\documentclass[12pt,a4paper]{article}
\usepackage[T1]{fontenc}
\usepackage[utf8]{inputenc}
\usepackage{lmodern}
\usepackage[margin=1in]{geometry}
\usepackage{microtype}
\usepackage{graphicx}
\usepackage[numbers]{natbib}
\usepackage{authblk}
\usepackage{amsfonts}
\usepackage{enumitem}
\usepackage{booktabs}
\usepackage{xcolor}
\usepackage{amsthm}
\usepackage{mathtools}
\usepackage{algorithm}
\usepackage{algorithmic}
\allowdisplaybreaks[4]
\newtheorem{remark}{Remark}[section]
\newtheorem{theorem}{Theorem}[section]
\newtheorem{lemma}{Lemma}[section]
\numberwithin{equation}{section}
\usepackage{subcaption}
\usepackage{pdflscape}
\usepackage{longtable}
\AddToHook{env/longtable/begin}{\setlength{\tabcolsep}{2pt}}
\usepackage{multirow}
\usepackage{array}
\usepackage{xurl}
\usepackage[hidelinks]{hyperref}
\newcommand\keywords[1]{\textbf{Keywords}: #1}

\title{An Unconditionally Stable Explicit
Robin–Robin Partitioned Scheme for
Fluid–Structure Interaction} 

\author[1]{Shihan Guo}
\author[2]{Ping Lin}
\author[3]{Yifan Wang}
\author[4]{Xiaohe Yue\thanks{Corresponding author: yuexiaohe666@gmail.com}}
\author[5]{Haibiao Zheng}

\affil[1]{
  Institute for Analysis and Numerics,
  Otto-von-Guericke-Universität Magdeburg,
  Magdeburg, 39106, Germany
}

\affil[2]{
  Division of Mathematics, University of Dundee,
  Dundee, DD1 4HN, United Kingdom
}

\affil[3]{
  Department of Mathematics and Statistics,
  Texas Tech University,
  Lubbock, TX, 79409, USA
}

\affil[4]{
  School of Mathematical Sciences,
  East China Normal University,
  Shanghai, 200241, China
}

\affil[5]{
  School of Mathematical Sciences,
  Ministry of Education Key Laboratory of Mathematics and Engineering
  Applications, Shanghai Key Laboratory of PMMP,
  East China Normal University,
  Shanghai, 200241, China
}

\date{}

\begin{document}
\maketitle

\begin{abstract}
We propose a parallelizable partitioned (loosely coupled) scheme for fluid
structure interaction (FSI) problems, specifically designed
\textcolor{black}{for efficient computation} in modern engineering simulations.
The FSI problem under consideration involves an incompressible viscous fluid,
governed by the Navier--Stokes equations, with a thick linear elastic structure.
The scheme adopts a Robin--Robin coupling condition, evaluating the right-hand
side of the Robin boundary terms at each time step solely from the previous-step
solutions. This explicit scheme allows the fluid and structure subproblems to be
solved entirely independently within each time step, eliminating the need for
staggered coupling or costly sub-iterations,
\textcolor{black}{and thus enabling their concurrent computation}.
\textcolor{black}{Various numerical experiments} demonstrate the stability,
accuracy, and \textcolor{black}{computational efficiency} of the proposed
approach, highlighting its potential for large-scale parallel FSI computations
in engineering applications.
\end{abstract}



\noindent\keywords{
Fluid-structure interaction, parallelizable partitioned scheme, unconditionally stable, ALE,  parallelization
}\par

\section{Introduction}
\label{sec1}
The study of fluid-structure interaction (FSI) has profound importance across disciplines such as biomedical engineering and environmental modeling \cite{ahajournal2024, vcanic2006modeling, doi:10.1137/20M1365144, wong2020shear}. These problems involve the tight interplay of
fluid and solid mechanics within time-evolving domains, often giving rise to nonlinear, multiphysics dynamics. Robust numerical methods for FSI are therefore essential, not only for advancing theoretical understanding but also for enabling computationally accurate and efficient large-scale simulations that support practical innovations ranging from biomedical device design to structural safety analysis, environmental resource management, etc. \cite{canic2024extendedonedimensionalreducedmodel, fluids7070222, bukavc2019fluid, bukavc2024computational}.

\textcolor{black}{
Numerical approaches to FSI are commonly classified as either monolithic or partitioned schemes. Monolithic methods solve the fluid and structural euqation simultaneously within a unified algebraic system \cite{failer2020parallel, richter2010finite, richter2013fully}, whereas partitioned methods retain separate fluid and structural subsolvers. The latter may employ either a fully implicit coupling procedure, generally requiring sub-iterations within each time step, or a loosely coupled strategy that advances the subproblems without such iterations. The loosely coupled splitting methods \cite{wang2018higher, fernandez2016convergence, bukac2016stability, burman2009stabilization, burman2014explicit, fernandez2015generalized, burman2022stability, Guzman2020, seboldt2021noniterative, gigante2021interface} are particularly attractive because of their modularity, their ability to reuse existing fluid and structural solvers, and their potential computational efficiency. By decomposing the coupled system into smaller subproblems, these methods can also reduce memory requirements and facilitate implementation on modern parallel computing architectures. Nevertheless, this advantage is not exclusive to partitioned methods: monolithic formulations equipped with carefully designed preconditioners can also achieve excellent computational efficiency and remain highly competitive \cite{gee2011truly}. 
Despite their advantages, loosely coupled methods remain challenging to design because stability and consistency must be maintained without fully resolving the interaction between the fluid and structure at each time step. This difficulty is particularly pronounced when the fluid and structural densities are comparable, in which case the added-mass effect can severely destabilize coupling procedure \cite{causin2005added, idelsohn2009fluid}.}

Considerable research has therefore focused on developing splitting strategies to overcome this instability from the added-mass effect. Here, we briefly review several contributions most relevant to the present work. Burman and Fern\`{a}ndez \cite{burman2009stabilization} proposed an explicit coupling scheme based on Nitsche's method and established its stability by introducing a pressure-stabilization term. However, this scheme exhibits a strong consistency error of order $\mathcal O(\Delta t/h)$. Building on this work, Burman and Fern\`{a}ndez subsequently investigated the connection between the Robin type conditions and the Nitsche based method, and compared them in \cite{burman2014explicit}. More recently, Robin–Robin partitioned schemes have gained significant attention for their ability to provide a stable interface treatment that effectively mitigates the added-mass effect and delivers reliable performance across a broad range of FSI applications \cite{durst2022recent}. Burman et al. proposed a loosely Robin-Robin scheme in \cite{burman2022stability}, and proved that this scheme is unconditionally stable with a $\mathcal{O}(\sqrt{T\Delta t})$ consistency error. They further gave error estimates in the energy norm of $\mathcal{O}(\sqrt{T}(\sqrt{\Delta t}+h)$ for the fully discrete system in \cite{Guzman2020}. Independently, Seboldt and Buka\v{c} \cite{seboldt2021noniterative} proposed an unconditionally stable, non-iterative Robin-Robin domain decomposition method, analyzing the convergence of the method and showing $\mathcal{O}(\sqrt{\Delta t})$ in time and optimal convergence in space. 
\textcolor{black}{
At the time semi-discrete level, the Robin--Robin schemes in
\cite{Guzman2020,seboldt2021noniterative} essentially coincide, while their
differences arise mainly at the fully discrete level and in the treatment of
the nonlinear problem. These sequential schemes are in fact first-order accurate but accuracy degrades
depending on the physical parameters. More recently, the nonlinear loosely
coupled Robin--Robin scheme was further analyzed in \cite{burman2025robin}, where a
nearly optimal temporal convergence result was established under additional
regularity assumptions. The analysis also shows that the splitting error
depends explicitly on the physical parameters.
}
\textcolor{black}{Beyond these developments, a common approach for improving the temporal accuracy of loosely coupled schemes is to introduce correction steps within each time step.
}

Despite their effectiveness in overcoming the added-mass effect, to the best of our knowledge, existing Robin–Robin partitioned approaches  for FSI are either inherently sequential, with one subsolver requiring updated interface data from the other within each time step, or iterative, requiring multiple subproblem solves to enforce the interface coupling.
\textcolor{black}{
These coupling procedures introduce sequential dependencies between the subsolvers and
additional computational overhead, which can limit the benefit of parallel execution, particularly when the fluid and structure workloads are strongly imbalanced.
}
In fact, parallelizable Robin-Robin schemes have been studied in many multi-physics domain problems such as Stokes-Darcy problems \cite{cao2014parallel} and fluid-poroelastic structure interaction (FPSI) problems \cite{guo2024fullyparallelizablelooselycoupled}.
In our recent work \cite{guo2024fullyparallelizablelooselycoupled}, we introduced an explicit, parallelizable Robin–Robin scheme for FPSI problems. \textcolor{black}{By constructing the interface conditions using only information available from the previous time step, the method fully decouples the coupled system into independent fluid and poroelastic-structure subproblems. Consequently, the two subproblems can be advanced concurrently at each time step without exchanging newly computed interface data.} This method prompts us to consider a question:
\begin{itemize}
    \item Can a fully parallelizable partitioned scheme based on Robin-Robin interface conditions be developed for FSI problems?
\end{itemize}

In the following sections, a positive answer would be provided. We propose in the present work a fully parallelizable, explicit Robin–Robin partitioned scheme for general FSI problems. The main contributions are threefold:
\begin{enumerate}[label=(\alph*)]
    \item Algorithmic development: We designed an explicit partitioned Robin–Robin FSI scheme that is unconditional stability and full parallelizable. A key advantage of the proposed scheme lies in its explicit handling and easy implementation of interface conditions, which allows the fluid and structure subproblems to be solved simultaneously without waiting for one another within each time step.
    \item Theoretical analysis: \textcolor{black}{We reformulate the original FSI problem as separate fluid and structural subproblems, coupled through explicit Robin-Robin interface conditions. We establish stability and consistency error estimates for the reformulated system and prove stability of the corresponding semi-discretized system.}
    {\color{black}
    \item We assess the proposed scheme using several numerical experiments, including a manufactured solution test, the Turek $\&$ Hron benchmark, a pressure wave driven blood flow simulation, and a 3D benchmark simulation to demonstrate accuracy, stability, and efficiency of the scheme.
    }
\end{enumerate}
By combining theoretical rigor with computational efficiency, the proposed approach contributes to the advancement of explicit partitioned FSI solvers and potentially enables large-scale multiphysics simulations in engineering and biomedical applications. 

The rest of this manuscript is organized as follows. In section \ref{section2}, we introduce the mathematical model for the FSI problem. Section \ref{section3} presents \textcolor{black}{the reformulated FSI system and} our explicit Robin-Robin partitioned scheme, giving a clear answer to the question above. \textcolor{black}{In section \ref{section4}, we first linearize the problem and establish the energy stability of the reformulated system. Next,  consistency error estimates between the original and reformulated system are studies. Finally, a discrete stability analysis for the semi-discretized system is provided.} Extensive numerical experiments in section 5 demonstrate the stability, accuracy, and \textcolor{black}{computational performance} of the proposed
approach.
{\color{black}

\begin{landscape}
\scriptsize
\setlength{\tabcolsep}{4pt}
\renewcommand{\arraystretch}{1.15}

\begin{longtable}{
>{\centering\arraybackslash}p{0.08\linewidth}
p{0.25\linewidth}
p{0.15\linewidth}
>{\centering\arraybackslash}p{0.08\linewidth}
p{0.25\linewidth}
p{0.15\linewidth}
}
\caption{Nomenclature and physical meaning of the main variables.}
\label{tab:nomenclature}\\

\toprule
\textbf{Symbol} & \textbf{Physical meaning} & \textbf{Domain / remark}
&
\textbf{Symbol} & \textbf{Physical meaning} & \textbf{Domain / remark}
\\
\midrule
\endfirsthead

\toprule
\textbf{Symbol} & \textbf{Physical meaning} & \textbf{Domain / remark}
&
\textbf{Symbol} & \textbf{Physical meaning} & \textbf{Domain / remark}
\\
\midrule
\endhead

\bottomrule
\endfoot

$\widehat{\Omega}_f$ 
& Reference fluid domain 
& Fixed reference configuration
&
$\Omega_f(t)$ 
& Current fluid domain 
& Time-dependent Eulerian configuration \\

$\widehat{\Omega}_s$ 
& Reference structure domain 
& Fixed reference configuration
&
$\Omega_s(t)$ 
& Current structure domain 
& Deformed structure configuration \\

$\widehat{\Gamma}$ 
& Fluid--structure interface in the reference configuration 
& Common boundary of $\widehat{\Omega}_f$ and $\widehat{\Omega}_s$
&
$\Gamma(t)$ 
& Fluid--structure interface in the current configuration 
& Time-dependent interface \\

$\widehat{\mathbf{T}}_f$ 
& ALE mapping for the fluid domain 
& Maps $\widehat{\Omega}_f$ to $\Omega_f(t)$
&
$\widehat{\mathbf{T}}_s$ 
& Lagrangian mapping for the structure domain 
& Maps $\widehat{\Omega}_s$ to $\Omega_s(t)$ \\

$\widehat{\mathbf{F}}_f$ 
& Fluid-domain deformation gradient 
& $\widehat{\mathbf{F}}_f=\widehat{\nabla}\widehat{\mathbf{T}}_f$
&
$\widehat{J}_f$ 
& Jacobian of the ALE mapping 
& $\widehat{J}_f=\det(\widehat{\mathbf{F}}_f)$ \\

$\widehat{\mathbf{F}}_s$ 
& Structure deformation gradient 
& $\widehat{\mathbf{F}}_s=\widehat{\nabla}\widehat{\mathbf{T}}_s$
&
$\widehat{J}_s$ 
& Jacobian of the structure mapping 
& $\widehat{J}_s=\det(\widehat{\mathbf{F}}_s)$ \\

$\mathbf{u}$ 
& Fluid velocity 
& Defined in $\Omega_f(t)$
&
$p$ 
& Fluid pressure 
& Defined in $\Omega_f(t)$ \\

$\widehat{\boldsymbol{\eta}}$ 
& Structure displacement 
& Defined in $\widehat{\Omega}_s$
&
$\widehat{\boldsymbol{\xi}}$ 
& Structure velocity 
& $\widehat{\boldsymbol{\xi}}=\partial_t\widehat{\boldsymbol{\eta}}$ \\

$\widehat{\boldsymbol{\eta}}_f$ 
& Fluid mesh displacement 
& ALE extension of the interface displacement
&
$\widehat{\mathbf{w}}$ 
& ALE mesh velocity 
& Usually $\widehat{\mathbf{w}}=\partial_t\widehat{\boldsymbol{\eta}}_f$ \\

$\boldsymbol{\sigma}_f$ 
& Fluid Cauchy stress tensor 
& $\boldsymbol{\sigma}_f=-p\mathbf{I}+2\mu_f\mathbf{D}(\mathbf{u})$
&
$\widehat{\boldsymbol{\sigma}}_s$ 
& Structure stress tensor 
& Saint-Venant--Kirchhoff model \\

$\mathbf{D}(\mathbf{u})$ 
& Strain-rate tensor of the fluid 
& $\mathbf{D}(\mathbf{u})=(\nabla\mathbf{u}+\nabla\mathbf{u}^{T})/2$
&
$\widehat{\mathbf{E}}$ 
& Green--Lagrange strain tensor 
& $\widehat{\mathbf{E}}=(\widehat{\mathbf{F}}_s^T\widehat{\mathbf{F}}_s-\mathbf{I})/2$ \\

$\rho_f$ 
& Fluid density 
& Positive physical parameter
&
$\rho_s$ 
& Structure density 
& Positive physical parameter \\

$\mu_f$ 
& Fluid dynamic viscosity 
& Positive physical parameter
&
$\mu_s$ 
& First Lamé parameter of the structure 
& Shear modulus \\

$\lambda_s$ 
& Second Lamé parameter of the structure 
& Related to compressibility
&
$\mathbf{f}$ 
& Fluid body force 
& Defined in the fluid domain \\

$\widehat{\mathbf{n}}_s$ 
& Outward unit normal vector to the reference structure domain 
& On $\partial\widehat{\Omega}_s$
&
$\widehat{\mathbf{n}}_f$ 
& Outward unit normal vector to the reference fluid domain 
& On $\partial\widehat{\Omega}_f$ \\

$\mathbf{n}_s$ 
& Outward unit normal vector to the current structure domain 
& On $\partial\Omega_s(t)$
&
$\mathbf{n}_f$ 
& Outward unit normal vector to the current fluid domain 
& On $\partial\Omega_f(t)$ \\

$\langle\cdot,\cdot\rangle_{\widehat{\Gamma}}$ 
& $L^2$ duality pairing or inner product on the interface 
& Over $\widehat{\Gamma}$
&
$\textcolor{black}{L}$ 
& \textcolor{black}{Robin parameter in both subproblems} 
& \textcolor{black}{Artificial coupling parameter} \\




$(\cdot,\cdot)_f$ 
& $L^2$ inner product over the fluid domain 
& Usually over $\widehat{\Omega}_f$
&
$(\cdot,\cdot)_s$ 
& $L^2$ inner product over the structure domain 
& Usually over $\widehat{\Omega}_s$ \\


\bottomrule
\end{longtable}

\normalsize
\end{landscape}
}

\section{Mathematical model}
\label{section2}
\subsection{Computational domains and mappings}
\begin{figure}[h]
  \centering
  \includegraphics[width=0.75\textwidth]{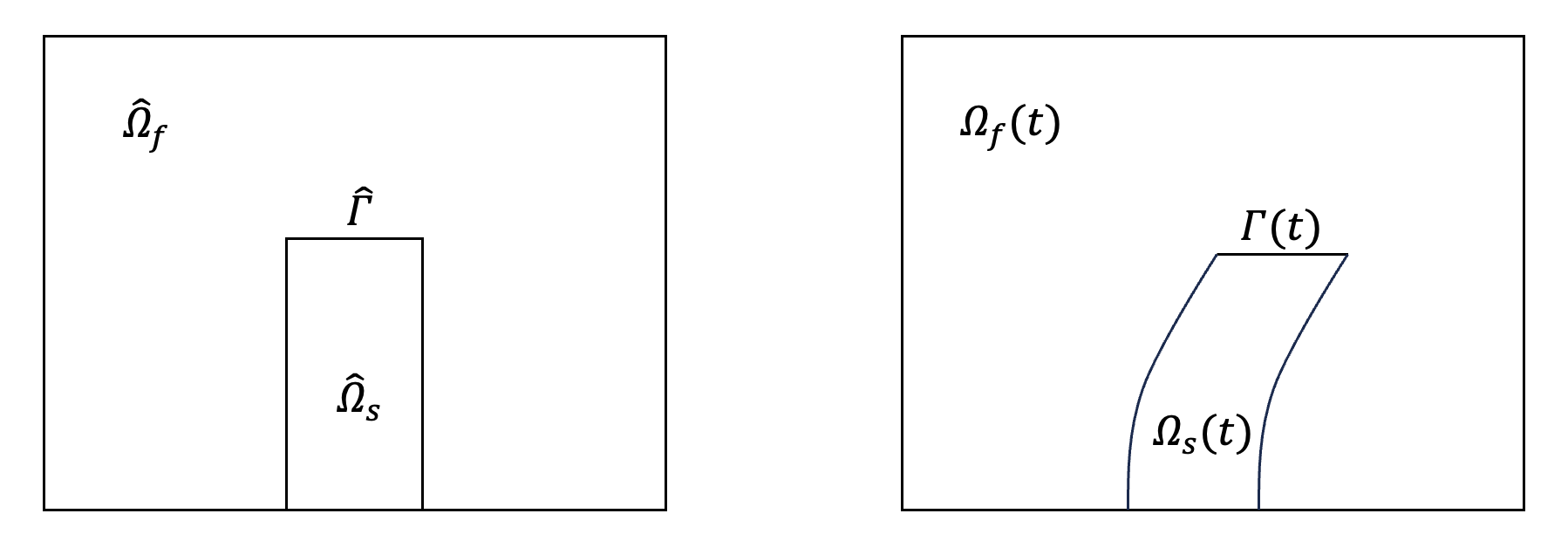}
  \caption{A sketch depicting the fluid–structure interaction domain: reference configuration (left) and physical configuration (right)}
  \label{fig:domain}
\end{figure}

As illustrated in Figure \ref{fig:domain}, let $\hat\Omega_s$ and $\hat\Omega_f$ denote the reference domains of the structure and the fluid, respectively. Owing to the deformability of the elastic structure, 
these domains evolve over time, giving rise to the time-dependent domains $\Omega_s(t)$ and $\Omega_f(t)$ in the Eulerian frame. The fluid structure interface is denoted by $\hat\Gamma$ in the reference domain and by $\Gamma(t)$ in the current domain. Throughout this study, we assume that both domains are regular, bounded regions in $\mathbb R^d$ with $d \in \{2, 3\}$.

Let $\hat{\boldsymbol\eta}:[0, T]\times\hat\Omega_s\to\mathbb{R}^d$ denote the displacement of the elastic material. To relate the reference configuration $\hat\Omega_s$ to the current domain $\Omega_s(t)$, we introduce the Lagrangian mapping as follows:
$$
\hat{\boldsymbol{T}}_s(\hat{\boldsymbol{x}}, t)=\hat{\boldsymbol{x}}+\hat{\boldsymbol{\eta}}(\hat{\boldsymbol{x}}, t), \quad \forall \hat{\boldsymbol{x}} \in \hat{\Omega}_s, ~t \in[0, T].
$$
We assume that $\hat{\boldsymbol{T}}_s$ is a $C^1$-diffeomorphism, ensuring a smooth and invertible correspondence between the reference and current domains. Consequently, for any scalar or vector function $\hat h:\hat{\Omega}_s\times[0, T]\to\mathbb R^i$, its counterpart in the Eulerian frame is defined as $h=\hat h\circ\hat{\boldsymbol{T}}_s^{-1}$ by the pullback:
$$
h:\Omega_s(t)\to\mathbb R^i, \quad h(\boldsymbol{x}, t) = \hat{h}\left(\hat{\boldsymbol{T}}_s^{-1}(\boldsymbol{x}, t), t\right).
$$
To track the evolution of the fluid domain over time, we further introduce the Arbitrary Lagrangian Eulerian (ALE) mapping $\hat{\boldsymbol{T}}_f: \hat{\Omega}_f \rightarrow \Omega_f(t)$, which provides a smooth, invertible transformation from the fixed reference domain to the current, physical configuration:
$$
\hat{\boldsymbol{T}}_f(\hat{\boldsymbol{x}},t)=\hat{\boldsymbol{x}}+\hat{\boldsymbol{\eta}}_f(\hat{\boldsymbol{x}}, t), \quad \forall \hat{\boldsymbol{x}} \in \hat{\Omega}_f,~ t \in[0, T] \text {. }
$$
Here, $\hat{\boldsymbol{\eta}}_f$ represents the fluid domain displacement, satisfying the interface condition $\hat{\boldsymbol \eta}=\hat{\boldsymbol \eta}_f$ on $\hat{\Gamma}$.
In constructing the ALE mappings, the displacement $\hat{\boldsymbol{\eta}}_f$ can be extended arbitrarily from the interface $\hat{\Gamma}$ into the interior of the fluid domain $\hat\Omega_f$. Typical strategies for such an extension include solving a harmonic extension or a biharmonic extension, see \cite{slyngstad2017verification} for further details.
For any scalar or vector function $f: \Omega_f(t) \times[0, T] \rightarrow \mathbb{R}^i$, we denote its representation in the reference domain as $\hat{f}=f \circ\hat{\boldsymbol T}_f$:
$$
\hat{f}: \hat{\Omega}_f \times[0, T] \rightarrow \mathbb{R}^i, \quad \hat{f}(\hat{\boldsymbol{x}}, t)=f\left(\boldsymbol{T}_f(\hat{\boldsymbol{x}}, t), t\right) .
$$
  Throughout this paper, quantities defined on the reference domain are denoted with a hat symbol $~\hat\cdot~$, while their counterparts in the current Eulerian configuration are written without a hat. 
\begin{remark}
In Section 2, the FSI model is first presented in its physical
configuration: the fluid equations are posed on the time-dependent Eulerian domain, whereas the thick structure is described on its fixed Lagrangian reference domain.
\textcolor{black}{In Section 3, the fluid equations are pulled back to a fixed reference domain through an ALE mapping. The main advantage of pulling the fluid equations back to the reference domain is that the fluid and structural subproblems will communicate through a fixed interface, thereby avoiding the transfer of interface data between moving and reference domains. Nevertheless, pulling the fluid equations back to the reference domain is not essential to the proposed coupling strategy. Alternatively, one may update the fluid mesh and solve the Navier--Stokes equations on the current fluid domain using the ALE time derivative, while retaining the mesh extension problem on the reference domain.}
\end{remark}

\subsection{The coupled Fluid-structure interaction problem}
We consider the fluid as an incompressible, viscous Newtonian fluid, and the structure as a Saint-Venant Kirchhoff material. Given the ALE mapping $\hat{\boldsymbol{T}}_f$, we denote $\hat{\boldsymbol{F}}_f=\hat\nabla\hat{\boldsymbol{T}}_f$ and $\hat{J}_f=\det(\hat{\boldsymbol{F}}_f)$ 
the deformation gradient and its associated Jacobian, respectively. Similarly, for the Lagrangian mapping $\hat{\boldsymbol{T}}_s$, we define $\hat{\boldsymbol{F}}_s=\hat\nabla\hat{\boldsymbol{T}}_s$ and $\hat{J}_s=\det(\hat{\boldsymbol{F}}_s)$ the deformation gradient and Jacobian of the Lagrangian map, respectively. 
The Saint-Venant Kirchhoff material constitutive relation is given by:
$$
\hat{\boldsymbol{\sigma}}_s=2\mu_s\hat{\boldsymbol{E}}+\lambda_s\text{tr}(\hat{\boldsymbol{E}})\boldsymbol{I},
$$
where $\hat{\boldsymbol{E}}=\frac{1}{2}(\hat{\boldsymbol{F}}_s^T\hat{\boldsymbol{F}}_s-\boldsymbol{I})$ is the Green-Lagrange strain tensor, $\boldsymbol{I}$ is the identity matrix, and $\mu_s$, $\lambda_s$ are Lamé parameters. 

The coupled fluid-structure interaction problem \cite{richter2017fluid} is to find fluid velocity $\boldsymbol{u}$ and pressure $p$, together with the structural displacement $\hat{\boldsymbol{\eta}}$ and velocity $\hat{\boldsymbol{\xi}}$, such that the following governing equations of the fluid, the structure, and the coupling conditions at the interface are satisfied:
\begin{equation}
  \begin{aligned}
    \rho_f(\frac{\partial\boldsymbol{u}}{\partial t} +(\boldsymbol{u}\cdot\nabla)\boldsymbol{u})-\nabla\cdot\boldsymbol{\sigma}_f&=\boldsymbol{f},\\
    \nabla\cdot\boldsymbol{u}&=0,\qquad&&\text{in}\ \Omega_f(t),\\
    \rho_s\frac{\partial\boldsymbol{\hat\xi}}{\partial t} - \hat{\nabla}\cdot(\hat{\boldsymbol{F}}_s\hat{\boldsymbol{\sigma}}_s)&=\hat{\boldsymbol{g}},\\
    \frac{d\hat{\boldsymbol{\eta}}}{d t}&=\hat{\boldsymbol{\xi}},\qquad&&\text{in}\ \hat\Omega_s,\\
    \boldsymbol{u}&=\boldsymbol{\xi},\qquad&&\text{on}\ \Gamma(t),\\
    \hat{\boldsymbol{n}}_f\cdot(\hat{J}_f\hat{\boldsymbol{\sigma}}_f\hat{\boldsymbol{F}}_f^{-T})&=\hat{\boldsymbol{n}}_f\cdot\hat{\boldsymbol{F}}_s\hat{\boldsymbol{\sigma}}_s,\qquad&&\text{on}\ \hat\Gamma.
  \end{aligned}
  \label{fsi_problem}
\end{equation}
Here, $\hat{\boldsymbol{n}}_f$ denotes the outward normal on $\partial\hat{\Omega}_f$. The fluid stress tensor is given by $\boldsymbol{\sigma}_f\left(\boldsymbol{u}, p\right)=-p \boldsymbol{I}+2 \mu_f \mathbb{D}(\boldsymbol{u})$, where $\mu_f$ is the fluid viscosity, and $\mathbb{D}(\boldsymbol{u})=$ $\left(\nabla \boldsymbol{u}+(\nabla \boldsymbol{u})^T\right) / 2$ is the strain rate tensor.

For simplicity, the remaining boundary segments of the fluid and structural domains are assumed to be subject either to homogeneous Dirichlet conditions or to natural stress-type conditions. These, together with appropriate initial conditions, are omitted from the formulation \eqref{fsi_problem} for brevity.

\begin{remark}
The well-posedness analysis of FSI problems is highly nontrivial, particularly for nonlinear models with moving interfaces. We refer the reader to~\cite{du2003analysis,du2004semidiscrete,coutand2005motion,coutand2006interaction} for representative results on the existence and uniqueness of solutions to both linearized and nonlinear FSI problems.
\end{remark}

\begin{remark}
In addition to the kinematic and dynamic coupling conditions, the FSI problem inherently satisfies a geometric compatibility condition, requiring the fluid and solid domains to remain in contact without separation or penetration. Therefore, the outward unit normals point in opposite directions on the interface, namely,
$
\hat{\boldsymbol{n}}_{f} = - \hat{\boldsymbol{n}}_{s}.
$
\end{remark}
\section{The Explicit Robin-Robin partitioned scheme}
\label{section3}
In this section, we present an explicit partitioned scheme for solving the FSI problem, designed to allow fully parallel computation. Within each time step, the original system is decoupled into a structural subproblem and a fluid subproblem, which are solved simultaneously without the need for sub-iterations. Compared with monolithic or strongly coupled schemes, the proposed approach notably reduces computational cost, while retaining stability and exhibiting \textcolor{black}{immunity} to added-mass effects.

\subsection{Robin-Robin interface conditions}
\label{continuous}
For simplicity, we consider a uniform partition of the time interval $[0, T]$ with step size $\Delta t$. Let $N\Delta t=T$ and $t_n=n\Delta t$. The variables with superscript $n+1$ represent the approximate solutions of corresponding exact solutions in sub-interval $[t_n, t_{n+1}]$. In each sub-interval $[t_n, t_{n+1}]$, we \textcolor{black}{reformulate the original FSI system into} the following two subproblems in the reference domain, subject to \textit{Robin interface conditions}:\\
\textbf{Fluid subproblem:} Given $\hat{J}_f^{n+1}, \left(\hat{\boldsymbol{F}}_f\right)^{n+1}, \hat{\boldsymbol{w}}^{n+1}, \hat{\boldsymbol{u}}^n$, $\hat{\boldsymbol{\xi}}^n$, $\hat{\mathcal F}^n = \left(\hat{J}_f\hat{\boldsymbol{\sigma}}_f\hat{\boldsymbol{F}}_f^{-T}\right)^{n}$ and $\hat{\mathcal{S}}^n=\left(\hat{\boldsymbol{F}}_s\hat{\boldsymbol{\sigma}}_s\right)^{n}$, find $\hat{\boldsymbol{u}}^{n+1}$ and $\hat{p}^{n+1}$ such that

\begin{align}
&\lefteqn{\rho_f\hat{J}_f^{n+1}\left[\partial_t\hat{\boldsymbol{u}}^{n+1}+ \left(\hat{\boldsymbol{F}}_f^{-1}\right)^{n+1}(\hat{\boldsymbol{u}}^{n+1} - \hat{\boldsymbol{w}}^{n+1})\cdot\hat\nabla\hat{\boldsymbol{u}}^{n+1}-\hat{\boldsymbol{f}}^{n+1}\right]=\hat\nabla\cdot\hat{\mathcal{F}}^{n+1},}\nonumber\\
&\hat\nabla\cdot\left[\hat{J}_f^{n+1}\left(\hat{\boldsymbol{F}}_f^{-1}\right)^{n+1}\hat{\boldsymbol{u}}^{n+1}\right]=0,&& \text{in }\hat{\Omega}_f\times[t_n, t_{n+1}],\nonumber\\
&\hat{\boldsymbol{u}}^{n+1}(\cdot,t_n)=\hat{\boldsymbol{u}}^n(\cdot,t_n), 
&& \text{in }\hat{\Omega}_f,\label{f_sub}\\
&\textcolor{black}{L}\hat{\boldsymbol{u}}^{n+1}+\hat{\mathcal F}^{n+1}\hat{\boldsymbol{n}}_f=\frac{\textcolor{black}{L}}{2}(\tilde{\boldsymbol{u}}^n+\tilde{\boldsymbol{\xi}}^n)+\frac{1}{2}(\tilde{\mathcal F}^{n}+\tilde{\mathcal{S}}^n)\hat{\boldsymbol{n}}_f,
&& \text{on }\hat{\Gamma}\times[t_n, t_{n+1}].\nonumber
\end{align}
Here, the ALE quantities at time level $n+1$, such as $\hat{J}_f^{n+1}$ and $\hat{\boldsymbol{w}}^{n+1}$ are treated as known data from \textcolor{black}{the} previous step, since the mesh motion has been updated first before the fluid step. Moreover, the equation $\hat{\boldsymbol{u}}^{n+1}(\cdot,t_n)=\hat{\boldsymbol{u}}^n(\cdot,t_n)$ means that the solution of $\hat{\boldsymbol{u}}^{n+1}$ at the beginning of time $t^n$ in the time interval $[t^n, t^{n+1}]$ must match with the final velocity obtained at the end of the previous step, namely $\boldsymbol{u}^{n}$ at time $n$.

\textbf{Structure subproblem:} Given $\hat{\boldsymbol{g}}^{n+1}$, $\hat{\boldsymbol{u}}^n, \hat{\boldsymbol{\eta}}^n, \hat{\boldsymbol{\xi}}^n, \hat{\mathcal{F}}^{n}$ and $\hat{\mathcal{S}}^{n}$, find $\hat{\boldsymbol{\eta}}^{n+1}$ and $\hat{\boldsymbol{\xi}}^{n+1}$ such that
  \begin{align}
    &\partial_t\hat{\boldsymbol{\eta}}^{n+1}=\hat{\boldsymbol{\xi}}^{n+1},\nonumber\\
    &\rho_s\partial_t\hat{\boldsymbol{\xi}}^{n+1} = \hat{\nabla}\cdot\hat{\mathcal{S}}^{n+1} + \hat{\boldsymbol{g}}^{n+1},&&\text{in }\hat{\Omega}_s\times[t_n,t_{n+1}],\nonumber\\
    &\hat{\boldsymbol{\xi}}^{n+1}(\cdot,t_n)=\hat{\boldsymbol{\xi}}^n(\cdot,t_n),\quad\hat{\boldsymbol{\eta}}^{n+1}(\cdot,t_n)=\hat{\boldsymbol{\eta}}^n(\cdot,t_n),
    && \text{in }\hat{\Omega}_s, \label{s_sub}\\
    &\textcolor{black}{L}\hat{\boldsymbol{\xi}}^{n+1}+\hat{\mathcal S}^{n+1}\hat{\boldsymbol{n}}_s=\frac{\textcolor{black}{L}}{2}(\tilde{\boldsymbol{u}}^n+\tilde{\boldsymbol{\xi}}^n)+\frac{1}{2}(\tilde{\mathcal F}^{n}+\tilde{\mathcal{S}}^n)\hat{\boldsymbol{n}}_s,
&& \text{on }\hat{\Gamma}\times[t_n, t_{n+1}].\nonumber
  \end{align}
\textcolor{black}{Here, $L$ denotes a positive Robin coupling parameter used in both subproblems.} The tilde superscript $\tilde{\cdot}^n$ ($n\ge1$) represents a time shifted quantity, which corresponds to the value of $\cdot^n$ shifted backward by a time step of $\Delta t$. For instance,
\begin{equation}
\label{tilde}
    \tilde{\boldsymbol{u}}^n=\hat{\boldsymbol{u}}^n(t-\Delta t), \qquad\tilde{\mathcal{F}}^n=\hat{\mathcal{F}}^n(t-\Delta t).
\end{equation}
So if $\hat{\boldsymbol{u}}^n$ is known on the interval $[t_{n-1}, t_n]$, then using $\tilde{\boldsymbol{u}}^n$ we can properly extend the definition of $\boldsymbol{u}^n$ in $[t_n, t_{n+1}]$. In particular, for initial time $n=0$, we set $\tilde\cdot^0(t)=\cdot(t_0)$ as an \textcolor{black}{identity}.

\begin{remark}
    It is worth noting that (\ref{f_sub}) and (\ref{s_sub}) are fully decoupled and can therefore be solved in parallel, since each subproblem relies solely on data from the previous time step. We may solve (\ref{f_sub}) and (\ref{s_sub}) instead of the original system \eqref{fsi_problem}. In section \ref{error estimate}, we will prove that, in a linearized case, solutions from reformulated problems converge to exact solutions in an appropriate norm as $\Delta t\rightarrow0$.
\end{remark}

Compared with the Robin--Robin loosely coupled splitting schemes proposed in~\cite{Guzman2020,seboldt2021noniterative}, an important feature of the present method is that the interface treatment is fully explicit on both sides. More precisely, the right-hand sides of the Robin interface conditions for both subproblems depend exclusively on data available from the previous time step. Therefore, the fluid and structure subproblems can be advanced concurrently within each time step, avoiding the sequential waiting in which one solver has to wait for the updated interface data produced by the other solver.

\begin{remark}
\textcolor{black}{The parameter $L$, , with unit of $Pa\cdot s/m$, is a numerical interface parameter introduced in the Robin--Robin splitting. It determines the relative
weight between the interface traction and the interface velocity/displacement
type quantities in the Robin transmission conditions. From this point of view,
$L$ may be interpreted as an interface impedance or relaxation parameter. 
The choice of $L$ may depend on the physical and geometric configuration, and affects the accuracy of the partitioned scheme. In the numerical stability analysis below, $L$ is assumed to be a positive constant, so that the interface terms generated by the Robin--Robin splitting can be controlled.
}
\end{remark}

\subsection{Adaptive ALE extension}
To extend the structure deformation $\hat{\boldsymbol{\eta}}$ from the interface $\hat\Gamma$ into the interior of the fluid domain, we adopt the adaptive ALE extension proposed by Masud and Hughes~\cite{masud1997space}, formulated as follows:
\begin{equation}
  \begin{aligned}
    &\hat\nabla\cdot\left(\left[1+\tau_m\right]\hat{\nabla}\hat{\boldsymbol{\eta}}_f\right) = 0 \qquad\text{ in } \hat{\Omega}_f,\\
    &\tau_m = \frac{1-\Delta_{min}/\Delta_{max}}{\Delta^e/\Delta_{max}}, \qquad e = 1, 2, \cdots, n_{el},\\
    &\hat{\boldsymbol{\eta}}_f = \hat{\boldsymbol{\eta}}(t) \qquad \text{on }\hat{\Gamma}, \qquad \quad\hat{\boldsymbol{\eta}}_f = \boldsymbol{0} \qquad \text{on }\partial\hat\Omega_f\setminus\hat\Gamma,
  \end{aligned}
  \label{adaptiveALE}
\end{equation}
where $\Delta^e$ denotes the area of the element under consideration; $n_{el}$ is the total number of elements in the fluid mesh; $\Delta_{min}$ and $\Delta_{max}$ represent the smallest and largest element areas, respectively. 
A weight function $\tau_m$ is introduced to impose a spatially varying stiffening effect within the computational domain.

\subsection{The splitting scheme}
In section \ref{continuous}, two Robin interface conditions are introduced to reformulate the original system in a continuous sense. In this part, we will apply the implicit Backward Euler scheme to obtain the semi-discretized form of (\ref{f_sub}) and (\ref{s_sub}). For notational convenience, we remark that the same notation $\cdot^n$ will be used to denote its discrete approximation $\cdot^n(t_n)$ at the time $t_n$.

Boundaries where Dirichlet conditions are imposed are denoted $\hat{\Sigma}_f^D$ for the fluid domain and $\hat{\Sigma}_s^D$ for the structure domain. For all $t \in [0,T]$, we define the following function spaces in the reference domain: 
\begin{align*}
\hat{\boldsymbol{V}}_{f} &= \{\hat{\boldsymbol{\boldsymbol{v}}}:\hat{\Omega}_{f}\rightarrow \mathbb{R}^{d}\ \big{|}\ \hat{\boldsymbol{v}}\in(H^{1}(\hat{\Omega}_f))^{d}, \hat{\boldsymbol{v}} = \boldsymbol{0}\ \text{on}\ \hat{\Sigma}_f^D\},\\%
\hat{Q}_{f} &= \{\hat{q}:\hat{\Omega}_{f}\rightarrow \mathbb{R}\ \big{|}  ~\hat{q}\in(L^{2}(\hat{\Omega}_f)) \},\\
\hat{\boldsymbol{V}}_{s} &= \{\hat{\boldsymbol{\zeta}}:\hat{\Omega}_{s}\rightarrow \mathbb{R}^{d}\ \big{|}\, \hat{\boldsymbol{\zeta}}\in(H^{1}(\hat{\Omega}_{s}))^{d}, \hat{\boldsymbol{\zeta}} = \boldsymbol{0}\ \text{on}\ \hat{\Sigma}_s^D \},
\end{align*} 
where $H^{1}$ denotes the usual Sobolev spaces. For notational convenience, we use $(\cdot,\cdot)_f$ and $(\cdot,\cdot)_s$ to represent integrals over the reference fluid and solid domains, respectively, and  $\left<\cdot,\cdot\right>_{\hat\Gamma}$ to denote the integrals over the reference fluid structure interface $\hat{\Gamma}$.

\textcolor{black}{By applying the implicit Backward Euler scheme to \eqref{f_sub} - \eqref{s_sub}}, we obtain the weak formulation of the proposed
Robin-Robin parallel loosely coupled splitting scheme \textcolor{black}{in semi-discretized form}, summarized in Algorithm \ref{semi}.

\begin{algorithm}[htbp!]
    \footnotesize
  \caption{\textcolor{black}{Parallel splitting scheme for  FSI system with moving interface}} 
  \label{semi} 
  Given $\hat{\boldsymbol{u}}^n$, $\hat{\boldsymbol{\eta}}^n$, $\hat{\boldsymbol{\xi}}^n$, $\left(\hat{J}_f\hat{\boldsymbol{\sigma}}_f\hat{\boldsymbol{F}}_f^{-T}\right)^{n}$ 
  and $\left(\hat{\boldsymbol{F}}_s\hat{\boldsymbol{\sigma}}_s\right)^{n}$, for $n=0, 1, 2, \cdots, N-1$, solve the following two sub-problems \textbf{in parallel.}

\noindent  \underline{\textbf{Fluid Subproblem:}}
\begin{enumerate}[label=\textbf{F\arabic*.}, leftmargin=2.2em, itemsep=0.6em]

 \item \textbf{ALE update.} 
 Solve for the ALE mapping $\boldsymbol{T}_f^{n+1}$ with $\tau_m$ defined in~(\ref{adaptiveALE}):
\begin{equation}
  \begin{aligned}
\hat\nabla\cdot\left(\left[1+\tau_m\right]\hat{\nabla}\hat{\boldsymbol{\eta}}_f^{n+1}\right) = 0 ~\text{in } \hat{\Omega}_f, 
\hat{\boldsymbol{\eta}}_f^{n+1} = \hat{\boldsymbol{\eta}}^{n} ~\text{on } \hat{\Gamma}, ~\hat{\boldsymbol{\eta}}_f^{n+1} = \boldsymbol{0}~\text{on }\partial\hat\Omega_f\setminus\hat\Gamma.
  \end{aligned}
  \label{alg_ALE}
\end{equation}
     Then, update $\hat{\boldsymbol{F}}^{n+1}_f$, $\hat{J}_f^{n+1}$ and $\boldsymbol{\omega}^{n+1}$ with $\boldsymbol{\omega}^{n+1}=d_t\hat{\boldsymbol{\eta}}_f^{n+1}$.
     
\item \textbf{Fluid solve.} 
{\color{black}
Solve the following problem for
$\hat{\boldsymbol{u}}^{n+1}$ and $\hat{p}^{n+1}$:
}
\begin{equation}
\begin{aligned} &\rho_f \left(\hat{J}_f^{n+1}\frac{\hat{\boldsymbol{u}}^{n+1}-\hat{\boldsymbol{u}}^n}{\Delta t} + \left(\hat{\boldsymbol{F}}_f^{-1}\right)^{n+1}(\hat{\boldsymbol{u}}^{n+1} - \hat{\boldsymbol{w}}^{n+1})\cdot\hat\nabla\hat{\boldsymbol{u}}^{n+1} , \hat{\boldsymbol{v}}\right)_f\\
      &+\left(\hat{J}_f^{n+1}\hat{\boldsymbol{\sigma}}_f^{n+1}\left(\hat{\boldsymbol{F}}_f^{-T}\right)^{n+1}, \hat\nabla\hat{\boldsymbol{v}}\right)_f+\left(\hat\nabla\cdot(\hat{J}_f^{n+1}\left(\hat{\boldsymbol{F}}_f^{-1}\right)^{n+1}\hat{\boldsymbol{u}}^{n+1}), \hat{q}\right)_f \\
      &+\left< \textcolor{black}{L}\hat{\boldsymbol{u}}^{n+1}, \hat{\boldsymbol{v}}\right>_{\hat{\Gamma}}=
      \left(\hat{\boldsymbol{f}}^{n+1},\hat{\boldsymbol{v}}\right)_f
      +
      \frac{\textcolor{black}{L}}{2}\left<\hat{\boldsymbol{u}}^n,\hat{\boldsymbol{v}}\right>_{\hat{\Gamma}}+\frac{\textcolor{black}{L}}{2}\left<\hat{\boldsymbol{\xi}}^n,\hat{\boldsymbol{v}}\right>_{\hat{\Gamma}}\\
      &+\frac{1}{2}\left<\left(\hat{\boldsymbol{F}}_s\hat{\boldsymbol{\sigma}}_s\right)^n\hat{\boldsymbol{n}}_f,\hat{\boldsymbol{v}}\right>_{\hat{\Gamma}}
      +\frac{1}{2}\left<\left(\hat{J}_f\hat{\boldsymbol{\sigma}}_f\hat{\boldsymbol{F}}_f^{-T}\right)^n\hat{\boldsymbol{n}}_f,\hat{\boldsymbol{v}}\right>_{\hat{\Gamma}}, \forall \left(\hat{\boldsymbol{v}},\hat{q}\right)\in\hat{\boldsymbol{V}}_f\times\hat{Q}_f.
     \end{aligned}
     \label{fluid}
    \end{equation}
  \item \textbf{Fluid-side interface update.} 
  Update the fluid-side interface traction-related quantity by the Robin condition:
  $\left(\hat{J}_f\hat{\boldsymbol{\sigma}}_f\hat{\boldsymbol{F}}_f^{-T}\right)^{n+1} 
=\frac{\textcolor{black}{L}}{2}(\hat{\boldsymbol{u}}^n+\hat{\boldsymbol{\xi}}^n)
+\frac{1}{2}((\hat{J}_f\hat{\boldsymbol{\sigma}}_f\hat{\boldsymbol{F}}_f^{-T})^{n} 
+ (\hat{\boldsymbol{F}}_{s}\hat{\boldsymbol{\sigma}}_{s})^{n} )\hat{\boldsymbol{n}}_f 
- \textcolor{black}{L}\hat{\boldsymbol{u}}^{n+1}$ .
\end{enumerate}

 \medskip
\noindent \underline{\textbf{Structure Subproblem:}}
 \begin{enumerate}[label=\textbf{S\arabic*.}, leftmargin=2.2em, itemsep=0.6em]

\item \textbf{Structure solve.} 
{\color{black}
Solve the following problem for
$\hat{\boldsymbol{\eta}}^{n+1}$ and
$\hat{\boldsymbol{\xi}}^{n+1}
=
(\hat{\boldsymbol{\eta}}^{n+1}
-\hat{\boldsymbol{\eta}}^{n})/\Delta t$:
}
\begin{equation}
     \begin{aligned}
      \label{A2}
      &\rho_s\left(\frac{\hat{\boldsymbol{\xi}}^{n+1}-\hat{\boldsymbol{\xi}}^n}{\Delta t}, \hat{\boldsymbol{\zeta}}\right)_s + \left(\hat{\boldsymbol{F}}_s^{n+1}\hat{\boldsymbol{\sigma}}_s^{n+1}, \hat{\nabla}\hat{\boldsymbol{\zeta}}\right)_s  + \left<\textcolor{black}{L}\hat{\boldsymbol{\xi}}^{n+1}, \hat{\boldsymbol{\zeta}}\right>_{\hat{\Gamma}}\\
      &=  \left(\hat{\boldsymbol{g}}^{n+1},\hat{\boldsymbol{\zeta}}\right)_s + \frac{\textcolor{black}{L}}{2}\left<\hat{\boldsymbol{u}}^n, \hat{\boldsymbol{\zeta}}\right>_{\hat{\Gamma}}+\frac{\textcolor{black}{L}}{2}\left<\hat{\boldsymbol{\xi}}^n,\hat{\boldsymbol{\zeta}}\right>_{\hat{\Gamma}}\\
      &+\frac{1}{2}\left<\left(\hat{\boldsymbol{F}}_s\hat{\boldsymbol{\sigma}}_s\right)^n\hat{\boldsymbol{n}}_s,\hat{\boldsymbol{\zeta}}\right>_{\hat{\Gamma}}+\frac{1}{2}\left<\left(\hat{J}_f\hat{\boldsymbol{\sigma}}_f\hat{\boldsymbol{F}}_f^{-T}\right)^n\hat{\boldsymbol{n}}_s, \hat{\boldsymbol{\zeta}}\right>_{\hat{\Gamma}}, \forall \hat{\boldsymbol{\zeta}}\in\hat{\boldsymbol{V}}_s.
    \end{aligned}
    \end{equation}
  \item \textbf{Structure-side interface update.} 
  Update the structure-side interface traction-related quantity by the Robin condition:
  $\left(\hat{\boldsymbol{F}}_s\hat{\boldsymbol{\sigma}}_s\right)^{n+1}
=\frac{\textcolor{black}{L}}{2}(\hat{\boldsymbol{u}}^n+\hat{\boldsymbol{\xi}}^n)
+\frac{1}{2}((\hat{J}_f\hat{\boldsymbol{\sigma}}_f\hat{\boldsymbol{F}}_f^{-T})^{n} 
+ (\hat{\boldsymbol{F}}_{s}\hat{\boldsymbol{\sigma}}_{s})^{n} )\hat{\boldsymbol{n}}_s 
- \textcolor{black}{L}\hat{\boldsymbol{\xi}}^{n+1}$.
\end{enumerate}
\textbf{Advancing to the next timestep.}
\end{algorithm}

{\color{black}
\subsection{Concurrency, load balance, and parallel efficiency}
\label{concurrency}
The proposed splitting removes the sequential dependency between the fluid and structure substeps and therefore permits concurrent execution. Its practical benefit depends on the workload and subsolvers. Let $T_f$ and $T_s$ denote the execution times of the fluid and structure subsolvers, respectively, under a given resource allocation. Both the sequential and concurrent implementations require the exchange of interface data between the two subproblems, although the underlying communication mechanisms may differ. We here assume that they are comparable and neglect them. The speedup attributable purely to concurrency is therefore bounded by
$$
\frac{T_{seq}}{T_{par}}=\frac{T_f + T_s}{\max\{T_f, T_s\}}\le 2,
$$
with equality for balanced subproblems.

At fixed resources, the situation is different from the idealized overlap argument above. Suppose that $P$ cores are available and that each subsolver exhibit ideal strong scaling, $T_f(P)=C_f/P$, $T_s(P)=C_s/P$. For the sequential scheme, $T_{seq}=(C_f+C_s)/P$ as two subproblems are consecutive. For the concurrent scheme, let $P_f+P_s=P$, so
$$
T_{par}=\max\{\frac{C_f}{P_f}, \frac{C_s}{P_s}\}.
$$
The minimum is attained when the two subsolvers are balanced, giving 
$$C_f/P_f=C_s/P_s=(C_f+C_s)/P,$$
and hence
$$\min_{P_f+P_s=P}T_{par}=T_{seq},$$
attained at $P_f = PC_f/(C_f+C_s)$. Under ideal scaling, concurrency alone does not reduce the execution time at a fixed total core count. 

In practical large-scale computations, however, the subsolvers generally exhibit diminishing strong-scaling efficiency as the number of cores increases. In this regime, concurrent execution provides an additional degree of freedom for resource allocation: instead of assigning all $P$ cores successively to subproblems that may already be approaching their strong-scaling limits, the available cores can be distributed between the fluid and structure solvers so that both operate in more efficient scaling regimes. The practical benefit of this resource allocation, and its departure from the ideal-scaling prediction above, are quantified in Section~\ref{sec:3d_concurrent}.
}

\section{\textcolor{black}{Numerical analysis}}
\label{section4}

{\color{black}
In this section, our analysis proceeds in three stages. First, we establish an energy stability estimate for the reformulated continuous
system (4.2) and (4.3). Next, before introducing any time discretization, we quantify the consistency error between the original FSI problem and this reformulated system. Finally, we discretize the reformulated system in time using the backward Euler method and show that the resulting time-discrete scheme is unconditionally stable.
}
\subsection{Linearized FSI problem}
 For the linearized FSI problem, both the fluid and structural domains are assumed to be fixed; therefore, we omit the notations involving hat that would otherwise indicate reference variables. Furthermore, we consider the fluid dynamics to be governed by the Stokes equations, while the structural response is described by the linear elasticity equations:
\begin{equation}
  \begin{aligned}
    \rho_f\frac{\partial\boldsymbol{u}}{\partial t}-\nabla\cdot\boldsymbol{\sigma}_f(\boldsymbol{u}, p)&=0,\\
    \nabla\cdot\boldsymbol{u}&=0,\qquad&&\text{in}\ \Omega_f,\\
    \rho_s\frac{\partial\boldsymbol{\xi}}{\partial t} - \nabla\cdot\boldsymbol{\sigma}_s(\boldsymbol{\eta})&=0,\\
    \boldsymbol{\sigma}_s(\boldsymbol{\eta})&=2\mu_s\mathbb D(\boldsymbol{\eta})+\lambda_s\nabla\cdot\boldsymbol{\eta}\boldsymbol{I},\\
    \frac{d{\boldsymbol{\eta}}}{d t}&=\boldsymbol{\xi},\qquad&&\text{in}\ \Omega_s,\\
    \boldsymbol{u}&=\boldsymbol{\xi},\quad\qquad&&\text{on}\ \Gamma,\\
    \boldsymbol{n}_f\cdot\boldsymbol{\sigma}_f&=\boldsymbol{n}_f\cdot\boldsymbol{\sigma}_s,\qquad&&\text{on}\ \Gamma.
  \end{aligned}
  \label{linear_fsi}
\end{equation}
Such assumptions are standard in the analysis of domain decomposition methods for fluid–structure interaction (FSI) problems \cite{Guzman2020, seboldt2021noniterative}. 
Under these simplifications, the fluid subproblem in sub-interval $[t_n, t_{n+1}]$ is formulated as follows: find $\boldsymbol{u}^{n+1}$ and $p^{n+1}$ such that
\begin{align}
    &\rho_f\partial_t\boldsymbol{u}^{n+1}-\nabla\cdot\boldsymbol{\sigma}_f^{n+1}=0.&&\text{in }\ \Omega_f\times[t_n, t_{n+1}],\nonumber\\
    &\nabla\cdot\boldsymbol{u}^{n+1}=0,&&\text{in }\ \Omega_f\times[t_n, t_{n+1}],\nonumber\\
    &\boldsymbol{u}^{n+1}(\cdot, t_n)=\boldsymbol{u}^n(\cdot, t_n),&&\text{in }\ \Omega_f,\label{linear_f_robin}\\
    &\textcolor{black}{L}\boldsymbol{u}^{n+1}+\boldsymbol{\sigma}_f^{n+1}\boldsymbol{n}_f=\frac{\textcolor{black}{L}}{2}\tilde{\boldsymbol{u}}^n+\frac{\textcolor{black}{L}}{2}\tilde{\boldsymbol{\xi}}^n+\frac{1}{2}\tilde{\boldsymbol{\sigma}}_f^n\boldsymbol{n}_f +\frac{1}{2}\tilde{\boldsymbol{\sigma}}_s^{n}\boldsymbol{n}_f, &&\text{on }\ \Gamma\times[t_n, t_{n+1}];\nonumber
  \end{align}
and the structure subproblem is: find $\boldsymbol{\eta}^{n+1}$ and $\boldsymbol{\xi}^{n+1}$ such that
\begin{align}
&\partial_t\boldsymbol{\eta}^{n+1}=\boldsymbol{\xi}^{n+1},&&\text{in }\Omega_s\times[t_n, t_{n+1}],\nonumber\\
  &\rho_s\partial_t\boldsymbol{\xi}^{n+1} - \nabla\cdot\boldsymbol{\sigma}_s^{n+1}=0,&&\text{in }\Omega_s\times[t_n, t_{n+1}],\nonumber\\
      &\boldsymbol{\xi}^{n+1}(\cdot,t_n)=\boldsymbol{\xi}^n(\cdot,t_n),\quad\boldsymbol{\eta}^{n+1}(\cdot,t_n)=\boldsymbol{\eta}^n(\cdot,t_n),
    && \text{in }\hat{\Omega}_s, \label{linear_s_robin}\\
  &\textcolor{black}{L}\boldsymbol{\xi}^{n+1}+\boldsymbol{\sigma}_s^{n+1}\boldsymbol{n}_s
  =\frac{\textcolor{black}{L}}{2}\tilde{\boldsymbol{u}}^n
  +\frac{\textcolor{black}{L}}{2}\tilde{\boldsymbol{\xi}}^n+\frac{1}{2}\tilde{\boldsymbol{\sigma}}_f^n\boldsymbol{n}_s +\frac{1}{2}\tilde{\boldsymbol{\sigma}}_s^{n}\boldsymbol{n}_s,&&\text{on }\Gamma\times[t_n, t_{n+1}].\nonumber
\end{align}

In section \ref{error estimate}, we analyze the energy stability and weak consistency of the reformulated problems \eqref{linear_f_robin} and \eqref{linear_s_robin} with respect to the original linearized problem \eqref{linear_fsi}. In section \ref{Stability analysis}, we further discretize problems \eqref{linear_f_robin} and \eqref{linear_s_robin} by the implicity Backward Euler scheme, and carry out a discrete stability analysis, showing that our scheme is unconditionally stable.

\begin{remark}
We assume that the reformulated problems \eqref{linear_f_robin} and \eqref{linear_s_robin} are well-posed and possess sufficient regularity. Specifically, we assume that $\boldsymbol{u}^{n+1}, \boldsymbol{\sigma}_f^{n+1}\boldsymbol{n}, \boldsymbol{\xi}^{n+1}$ and $\boldsymbol{\sigma}_s^{n+1}\boldsymbol{n}$ are all in $L^2(t_n, t_{n+1}; L^2(\Gamma))$. We note that, for instance, $\boldsymbol{\sigma}^{n+1}_s\boldsymbol{n}\in L^2(t_n, t_{n+1};L^2(\Gamma))$ implies $\boldsymbol{\xi}^{n+1}\in L^2(t_n,t_{n+1};H^{3/2}(\Omega_s))$. These regularity results are also used in \cite{burman2022stability}. 
\end{remark}

\subsection{\textcolor{black}{Stability and consistency error estimates}}
\label{error estimate}
In this section, we show that the splitting method with Robin-Robin interface conditions proposed above is \textcolor{black}{energy stable and} weakly consistent. Specifically, we will prove that the splitting error is $\sqrt{\Delta t}$. 

For simplicity of notation, let $\Vert\cdot\Vert_f$ denote $\Vert\cdot\Vert_{L^2(\Omega_f)}$, $\Vert\cdot\Vert_s$ denote $\Vert\cdot\Vert_{L^2(\Omega_s)}$ and $\Vert\cdot\Vert_\Gamma$ denote $\Vert\cdot\Vert_{L^2(\Gamma)}$. In the following, we also introduce the elastic energy of the structure defined by:
\begin{equation}
    \Vert\boldsymbol{\eta}\Vert^2_E=2\mu_s\Vert\mathbb D(\boldsymbol{\eta})\Vert^2_s+\lambda_s\Vert\nabla\cdot\boldsymbol{\eta}\Vert^2_s.
\label{energynorm}
\end{equation}
Besides, the polarization identity will be applied in the following analyses.
\begin{equation}
\label{polar}
    2(a-b)a = a^2 - b^2 + (a-b)^2.
\end{equation}

\color{black}
For $1\le n\le N$, introduce the following quantities:
\begin{equation*}
  \begin{aligned}
      \mathbb E^{n}&=\frac{\rho_s}{2}\Vert\boldsymbol{\xi}^n(t_n)\Vert^2_s+\frac{1}{2}\Vert\boldsymbol{\eta}^n(t_n)\Vert_E^2+\frac{\rho_f}{2}\Vert\boldsymbol{u}^n(t_n)\Vert^2_f,\\
      \mathbb H^n&=2\mu_f\int^{t_n}_{t_{n-1}}\Vert\mathbb D(\boldsymbol{u}^n)\Vert^2_f dt,\\
      \mathbb I^n&=\frac{L}{2}\int^{t_n}_{t_{n-1}}\Vert u^n\Vert^2_\Gamma +\Vert \boldsymbol{\xi}^n\Vert^2_\Gamma dt + \frac{1}{2L}\int^{t_n}_{t_{n-1}}\Vert \boldsymbol{\sigma}_f^{n}\boldsymbol{n}_f\Vert^2_\Gamma +\Vert \boldsymbol{\sigma}_s^{n}\boldsymbol{n}_s\Vert^2_\Gamma dt.
  \end{aligned}
\end{equation*}
The initial energy is
$$
\mathbb E^{0}=\frac{\rho_s}{2}\Vert\boldsymbol{\xi}^0(t_0)\Vert^2_s+\frac{1}{2}\Vert\boldsymbol{\eta}^0(t_0)\Vert_E^2+\frac{\rho_f}{2}\Vert\boldsymbol{u}^0(t_0)\Vert^2_f.
$$
Since $\tilde\cdot^0(t)=\cdot^0(t_0)$, set
$$
\mathbb I^n=\frac{L\Delta t}{2}\left(\Vert u^0(t_0)\Vert^2_\Gamma +\Vert \boldsymbol{\xi}^0(t_0)\Vert^2_\Gamma\right)+ \frac{\Delta t}{2L}\left(\Vert \boldsymbol{\sigma}_f^{0}(t_0)\boldsymbol{n}_f\Vert^2_\Gamma +\Vert \boldsymbol{\sigma}_s^{0}(t_0)\boldsymbol{n}_s\Vert^2_\Gamma\right).
$$
\begin{theorem}\label{cstability}
    Let $(\boldsymbol{u}^{n+1}, p^{n+1}, \boldsymbol{\eta}^{n+1}, \boldsymbol{\xi}^{n+1})$ be a sufficient regular solution of the reformulated systems \eqref{linear_f_robin}-\eqref{linear_s_robin} on $[t_n, t_{n+1}]$. Then, for every $1\le M\le N$,
    \begin{equation}
        \mathbb E^M + \sum_{n=1}^M\mathbb{H}^n + \mathbb{I}^M\le \mathbb{E}^0 + \mathbb{I}^0.
        \label{csta}
    \end{equation}
\end{theorem}
\begin{proof}
    See Appendix.
\end{proof}
\color{black}

Let $\mathbb{\boldsymbol{U}}, \mathbb{P}, \boldsymbol{\sigma}_\mathbb{{F}}, \mathbb{\boldsymbol{\Psi}}, \mathbb{\boldsymbol{\Xi}}$ and $\boldsymbol{\sigma}_\mathbb{{S}}$ be the exact solutions of the linearized system \eqref{linear_fsi}. We use the notation $\mathbb{\boldsymbol{U}}^{n+1}(t, \boldsymbol{x}) = \mathbb{\boldsymbol{U}}(t, \boldsymbol{x})$ for $t_{n}\leq t \leq t_{n+1}$ and $\boldsymbol{x}\in\boldsymbol{\Omega}$; similar notations are applied to other variables. We define the errors as follows:
$$
\boldsymbol{e}_{u}^{n} = \mathbb{\boldsymbol{U}}^{n} - \boldsymbol{u}^{n}, 
~\boldsymbol{e}_{f}^{n} = \boldsymbol{\sigma}_{\mathbb{F}}^{n} - \boldsymbol{\sigma}_{f}^{n},
$$
$$
\boldsymbol{e}_{\eta}^{n} = \mathbb{\boldsymbol{\Psi}}^{n} - \boldsymbol{\eta}^{n}, 
~\boldsymbol{e}_{\xi}^{n} = \mathbb{\boldsymbol{\Xi}}^{n} - \boldsymbol{\xi}^{n}, 
~\boldsymbol{e}_{s}^{n} = \boldsymbol{\sigma}_{\mathbb{S}}^{n} - \boldsymbol{\sigma}_{s}^{n}.
$$
Combining derivations of the original coupling conditions and Robin interface conditions with \eqref{tilde} yields the following equations:
    \begin{align}
        L\boldsymbol{e}_{u}^{n+1} + \boldsymbol{e}_{f}^{n+1}\boldsymbol{n}_{f}&=
        \frac{L}{2}\boldsymbol{e}_{u}^{n} + \frac{L}{2}\boldsymbol{e}_{\xi}^{n}
        + \frac{1}{2}\boldsymbol{e}_{f}^{n}\boldsymbol{n}_{f} 
        + \frac{1}{2}\boldsymbol{e}_{s}^{n}\boldsymbol{n}_{f} 
        + \boldsymbol{g}_{1}^{n+1}, \label{l1}\\
        L\boldsymbol{e}_{\xi}^{n+1} + \boldsymbol{e}_{s}^{n+1}\boldsymbol{n}_{s}&=
        \frac{L}{2}\boldsymbol{e}_{u}^{n} + \frac{L}{2}\boldsymbol{e}_{\xi}^{n}
        + \frac{1}{2}\boldsymbol{e}_{f}^{n}\boldsymbol{n}_{s} 
        + \frac{1}{2}\boldsymbol{e}_{s}^{n}\boldsymbol{n}_{s} 
        + \boldsymbol{g}_{2}^{n+1},\label{l2}
    \end{align}
where $\boldsymbol{g}_{1}^{n+1}$ and $\boldsymbol{g}_{2}^{n+1}$ are defined as follows:
\begin{equation*}
    \boldsymbol{g}_{1}^{n+1} 
=  \frac{L}{2}(\boldsymbol{U}^{n+1} - \tilde{\boldsymbol{U}}^{n} + \boldsymbol{\Xi}^{n+1} - \tilde{\boldsymbol{\Xi}}^{n}) 
+ \frac{1}{2}(\boldsymbol{\sigma}_{\mathbb{F}}^{n+1}\boldsymbol{n}_{f} - \tilde{\boldsymbol{\sigma}}_{\mathbb{F}}^{n}\boldsymbol{n}_{f}+\boldsymbol{\sigma}_{\mathbb{S}}^{n+1}\boldsymbol{n}_{f} - \tilde{\boldsymbol{\sigma}}_{\mathbb{S}}^{n}\boldsymbol{n}_{f}),
\end{equation*}
\begin{equation*}
\boldsymbol{g}_{2}^{n+1} 
=  \frac{L}{2}(\boldsymbol{U}^{n+1} - \tilde{\boldsymbol{U}}^{n} + \boldsymbol{\Xi}^{n+1} - \tilde{\boldsymbol{\Xi}}^{n}) 
+ \frac{1}{2}(\boldsymbol{\sigma}_{\mathbb{F}}^{n+1}\boldsymbol{n}_{s} - \tilde{\boldsymbol{\sigma}}_{\mathbb{F}}^{n}\boldsymbol{n}_{s}+\boldsymbol{\sigma}_{\mathbb{S}}^{n+1}\boldsymbol{n}_{s} - \tilde{\boldsymbol{\sigma}}_{\mathbb{S}}^{n}\boldsymbol{n}_{s}).
\end{equation*}

We assume that $\partial_t\boldsymbol{U}, \partial_t\boldsymbol{\sigma}_\mathbb{F}\boldsymbol{n}, \partial_t\boldsymbol{\Xi}$ and $\partial_t\boldsymbol{\sigma}_{\mathbb{S}}\boldsymbol{n}$ are in $L^2(0, T; L^2(\Gamma))$. The proof of the error estimates relies on the following lemma regarding local errors.

\begin{lemma}
\label{lemma}
For $\boldsymbol{g}_{1}^{n+1}$ and $\boldsymbol{g}_{2}^{n+1}$, we have the following estimates for $n\ge1$
    \begin{equation}
        \begin{aligned}
            \int_{t_{n}}^{t_{n+1}}\Vert\boldsymbol{g}_{1}^{n+1}(s)\Vert_{\Gamma}^{2}ds \leq \frac{(\Delta t)^{2}}{4}\int_{t_{n-1}}^{t_{n+1}}&\left(
            L^2\Vert\partial_{t}\mathbb{\boldsymbol{U}}\Vert_{\Gamma}^{2} 
            +L^2\Vert\partial_{t}\mathbb{\boldsymbol{\Xi}}\Vert_{\Gamma}^{2}\right.\\
            &\left.+\Vert\partial_{t}\boldsymbol{\sigma}_{\mathbb{F}}\boldsymbol{n}_f\Vert_{\Gamma}^{2} 
            +\Vert\partial_{t}\boldsymbol{\sigma}_{\mathbb{S}}\boldsymbol{n}_f\Vert_{\Gamma}^{2} 
            \right)ds,
        \end{aligned}
        \label{Lemma1}
    \end{equation}
        \begin{equation}
        \begin{aligned}
            \int_{t_{n}}^{t_{n+1}}\Vert\boldsymbol{g}_{2}^{n+1}(s)\Vert_{\Gamma}^{2}ds \leq \frac{(\Delta t)^{2}}{4}\int_{t_{n-1}}^{t_{n+1}}&\left(
            L^2\Vert\partial_{t}\mathbb{\boldsymbol{U}}\Vert_{\Gamma}^{2} 
            +L^2\Vert\partial_{t}\mathbb{\boldsymbol{\Xi}}\Vert_{\Gamma}^{2}\right.\\
            &\left.+\Vert\partial_{t}\boldsymbol{\sigma}_{\mathbb{F}}\boldsymbol{n}_s\Vert_{\Gamma}^{2} 
            +\Vert\partial_{t}\boldsymbol{\sigma}_{\mathbb{S}}\boldsymbol{n}_s\Vert_{\Gamma}^{2} 
            \right)ds.
        \end{aligned}
        \label{lamme2}
    \end{equation}
For $n=0$ we have
    \begin{equation}
        \begin{aligned}
            \int_{t_{0}}^{t_1}\Vert\boldsymbol{g}_{1}^1(s)\Vert_{\Gamma}^{2}ds \leq \frac{(\Delta t)^{2}}{4}\int_{t_0}^{t_1}&\left(
            L^2\Vert\partial_{t}\mathbb{\boldsymbol{U}}\Vert_{\Gamma}^{2} 
            +L^2\Vert\partial_{t}\mathbb{\boldsymbol{\Xi}}\Vert_{\Gamma}^{2}\right.\\
            &\left.+\Vert\partial_{t}\boldsymbol{\sigma}_{\mathbb{F}}\boldsymbol{n}_f\Vert_{\Gamma}^{2} 
            +\Vert\partial_{t}\boldsymbol{\sigma}_{\mathbb{S}}\boldsymbol{n}_f\Vert_{\Gamma}^{2} 
            \right)ds,
        \end{aligned}
    \end{equation}
        \begin{equation}
        \begin{aligned}
            \int_{t_0}^{t_1}\Vert\boldsymbol{g}_{2}^1(s)\Vert_{\Gamma}^{2}ds \leq \frac{(\Delta t)^{2}}{4}\int_{t_0}^{t_1}&\left(
            L^2\Vert\partial_{t}\mathbb{\boldsymbol{U}}\Vert_{\Gamma}^{2} 
            +L^2\Vert\partial_{t}\mathbb{\boldsymbol{\Xi}}\Vert_{\Gamma}^{2}\right.\\
            &\left.+\Vert\partial_{t}\boldsymbol{\sigma}_{\mathbb{F}}\boldsymbol{n}_s\Vert_{\Gamma}^{2} 
            +\Vert\partial_{t}\boldsymbol{\sigma}_{\mathbb{S}}\boldsymbol{n}_s\Vert_{\Gamma}^{2} 
            \right)ds.
        \end{aligned}
    \end{equation}
\end{lemma}
\begin{proof}
Only the proof of \eqref{Lemma1} is shown here, since the others proceed in a similar manner. For $t\in [t_{n-1}, t_{n}]$, we have:
\begin{equation*}
    \boldsymbol{U}^{n+1}-\tilde{\boldsymbol{U}}^n=\boldsymbol{U}(t)-\boldsymbol{U}(t-\Delta t) = \int^t_{t-\Delta t}\partial_t\boldsymbol{U}(s)ds.
\end{equation*}
By applying the Cauchy-Schwarz inequality, we have:
        \begin{align*}
        &\quad\int_{t_{n}}^{t_{n+1}}\Vert\mathbf{\boldsymbol{U}}^{n+1}(t) - \tilde{\boldsymbol{U}}^{n }(t)\Vert_{\Gamma}^{2}dt 
            =  \int_{t_{n}}^{t_{n+1}}\int_\Gamma\left(\mathbf{\boldsymbol{U}}(t,x) - \mathbf{\boldsymbol{U}}(t - \Delta t, x)\right)^2dxdt \\
            &= \int_{t_{n}}^{t_{n+1}}\int_\Gamma\left(\int_{t-\Delta t}^t\partial_t\boldsymbol{U}(s)ds\right)^2dxdt 
            \leq \int_{t_{n}}^{t_{n+1}}\int_\Gamma\Delta t\int_{t-\Delta t}^t\left(\partial_t\boldsymbol{U}(s)\right)^2ds dxdt\\
            &\leq \Delta t\int_\Gamma\int_{t_{n}}^{t_{n+1}}\int_{t_{n-1}}^{t_{n+1}}\left(\partial_t\boldsymbol{U}(s)\right)^2dsdxdt
            = (\Delta t)^{2}\int_{t_{n-1}}^{t_{n+1}}
           \Vert\partial_{t}\boldsymbol{U}(s)\Vert_{\Gamma}^{2}ds.
        \end{align*}
\end{proof}

We introduce the following quantities which will be used in our error estimates:
$$
\boldsymbol{E}^{n} = \frac{\rho_f}{2}\Vert\boldsymbol{e}_{u}^{n}(t_n)\Vert_{f}^{2} + \frac{\rho_s}{2}\Vert\boldsymbol{e}_{\xi}^{n}(t_n)\Vert_{s}^{2} + 
\frac{1}{2}\Vert\boldsymbol{e}_{\eta}^{n}(t_n)\Vert_{E}^{2},
$$
$$
\boldsymbol{T}^{n} = 2\mu_{f}\int_{t_{n-1}}^{t_{n}}\Vert\mathbb{D}(\boldsymbol{e}_{u}^{n}(s))\Vert_{f}^{2}ds,
$$
$$
\boldsymbol{S}^{n} = \frac{L}{2}\int_{t_{n-1}}^{t_n}(\Vert\boldsymbol{e}_{u}^{n}(s)\Vert_{\Gamma}^{2} + \Vert\boldsymbol{e}_{\xi}^{n}(s)\Vert_{\Gamma}^{2})ds
+ \frac{1}{2L}\int_{t_{n-1}}^{t_n}(\Vert\boldsymbol{e}_{f}^{n}(s)\boldsymbol{n}_{f}\Vert_{\Gamma}^{2} + \Vert\boldsymbol{e}_{s}^{n}(s)\boldsymbol{n}_{f}\Vert_{\Gamma}^{2})ds.
$$
Note that $\boldsymbol{E}^{0} = \boldsymbol{S}^{0} = 0$.

\begin{theorem}
\label{Error estimate}
Let $(\mathbb{\boldsymbol{U}}, \mathbb{P}, \boldsymbol{\sigma}_\mathbb{{F}}, \mathbb{\boldsymbol{\Psi}}, \mathbb{\boldsymbol{\Xi}}, \boldsymbol{\sigma}_\mathbb{{S}})$ solve the linearized system described in \eqref{linear_fsi} for $0 \leq t \leq T$. Furthermore, let $(\boldsymbol{u}^{n+1}, p^{n+1}, \boldsymbol{\sigma}_{f}^{n+1}, \boldsymbol{\eta}^{n+1}, \boldsymbol{\xi}^{n+1}, \boldsymbol{\sigma}_{s}^{n+1})$ denote the solution of the decoupled scheme given by \eqref{linear_f_robin} and \eqref{linear_s_robin} for $n = 0, 1, \cdots, N-1$. If $T = N\Delta t$ with $N \geq 1$, the following estimate holds:
\begin{equation}
    \begin{aligned}
        \boldsymbol{E}^{N} + \sum_{n=1}^{N}\boldsymbol{T}^n + \boldsymbol{S}^{N}
\leq CT\Delta t\left(
L\Vert\partial_{t}\mathbb{\boldsymbol{U}}\Vert_{L^2(0, T; L^2(\Gamma))}^{2}
+ L\Vert\partial_{t}\mathbb{\boldsymbol{\Xi}}\Vert_{L^2(0, T; L^2(\Gamma))}^{2}\right. \\
\left.+ \frac{1}{L}\Vert\partial_{t}\boldsymbol{\sigma}_{\mathbb{F}}\Vert_{L^2(0, T; L^2(\Gamma))}^{2}
+ \frac{1}{L}\Vert\partial_{t}\boldsymbol{\sigma}_{\mathbb{S}}\Vert_{L^2(0, T; L^2(\Gamma))}^{2}
\right),
    \end{aligned}
    \label{error result}
\end{equation}
where $C$ is a positive constant.
\end{theorem}

\begin{proof}
By multiplying \eqref{linear_fsi}, \eqref{linear_f_robin} and \eqref{linear_s_robin} by corresponding test functions $(\boldsymbol{e}_u, e_p, \boldsymbol{e}_\eta,\boldsymbol{e}_\xi)$, integrating over their respective domains, and simple calculations, one yields:
\begin{equation}
\frac{\rho_f}{2}\partial_{t}\Vert\boldsymbol{e}_{u}^{n+1}\Vert^{2}_{f} 
+ 2\mu_{f}\Vert\mathbb D(\boldsymbol{e}_u^{n+1}) \Vert^{2}_{f} 
+ \frac{\rho_s}{2}\partial_{t}\Vert\boldsymbol{e}_{\xi}^{n+1}\Vert^{2}_{s}
+ \frac{1}{2}\partial_{t}\Vert\boldsymbol{e_{\eta}}^{n+1}\Vert^{2}_{E} = J^{n+1},
\end{equation}
where
\begin{equation}
J^{n+1} = \int_{\Gamma}\boldsymbol{e}_{f}^{n+1}\boldsymbol{n}_{f}\cdot\boldsymbol{e}_{u}^{n+1} 
+ \int_{\Gamma}\boldsymbol{e}_{s}^{n+1}\boldsymbol{n}_{s}\cdot\boldsymbol{e}_{\xi}^{n+1} .
\label{J}
\end{equation}
By replacing \eqref{l1} and \eqref{l2} into \eqref{J}, we have:
\begin{equation}
    \begin{aligned}
        J^{n+1}
&= \int_{\Gamma}(\frac{L}{2}\boldsymbol{e}_{u}^{n} + 
\frac{L}{2}\boldsymbol{e}_{\xi}^{n}
 - L\boldsymbol{e}_{u}^{n+1} + \boldsymbol{g}_{1}^{n+1})
 \cdot\boldsymbol{e}_{u}^{n+1} \\
&+ \int_{\Gamma}(\frac{L}{2}\boldsymbol{e}_{u}^{n} + 
\frac{L}{2}\boldsymbol{e}_{\xi}^{n}
 - L\boldsymbol{e}_{\xi}^{n+1} +       \boldsymbol{g}_{2}^{n+1})\cdot\boldsymbol{e}_{\xi}^{n+1} \\
&+ \int_{\Gamma} (\frac{1}{2}\boldsymbol{e}_{f}^{n}\boldsymbol{n}_{f}
+ \frac{1}{2}\boldsymbol{e}_{s}^{n}\boldsymbol{n}_{f})
\cdot(\boldsymbol{e}_{u}^{n+1} - \boldsymbol{e}_{\xi}^{n+1} ) \\
&=: P_1 + P_2 + P_3.
    \end{aligned}
\end{equation}
Here, we denote three terms in the formulation of $J^{n+1}$ above as $P_1$, $P_2$ and $P_3$, respectively, for notational simplicity. Obviously from \eqref{polar},
\begin{equation}
    P_1 = \frac{L}{2}\left(\left\Vert\frac{\boldsymbol{e}_{u}^n }{2}
    +\frac{\boldsymbol{e}_{\xi}^n}{2} 
    + \frac{\boldsymbol{g}_{1}^{n+1}}{L}\right\Vert_{\Gamma}^{2} 
    - \left\Vert\boldsymbol{e}_{u}^{n+1} \right\Vert_{\Gamma}^{2} 
    - \left\Vert\frac{\boldsymbol{e}_{u}^n}{2}
    +\frac{\boldsymbol{e}_{\xi}^n}{2}
    + \frac{\boldsymbol{g}_{1}^{n+1}}{L}
    - \boldsymbol{e}_{u}^{n+1}\right\Vert_{\Gamma}^{2}\right),
    \label{p1}
\end{equation}
\begin{equation}
    P_2 = \frac{L}{2}\left(\left\Vert\frac{\boldsymbol{e}_{u}^n }{2}
    +\frac{\boldsymbol{e}_{\xi}^n}{2} 
    + \frac{\boldsymbol{g}_{2}^{n+1}}{L}\right\Vert_{\Gamma}^{2} 
    - \Vert\boldsymbol{e}_{\xi}^{n+1} \Vert_{\Gamma}^{2} 
    - \left\Vert\frac{\boldsymbol{e}_{u}^n}{2}
    +\frac{\boldsymbol{e}_{\xi}^n}{2}
    + \frac{\boldsymbol{g}_{2}^{n+1}}{L}
    - \boldsymbol{e}_{\xi}^{n+1}\right\Vert_{\Gamma}^{2}\right).
    \label{p2}
\end{equation}
To deal with term $P_3$, subtracting \eqref{l1} from \eqref{l2} yields:
$$
\boldsymbol{e}_{u}^{n+1} - \boldsymbol{e}_{\xi}^{n+1} = \frac{1}{L}(\boldsymbol{e}_{f}^{n}\boldsymbol{n}_{f} + \boldsymbol{e}_{s}^{n}\boldsymbol{n}_{f}
- \boldsymbol{e}_{f}^{n+1}\boldsymbol{n}_{f} -\boldsymbol{e}_{s}^{n+1}\boldsymbol{n}_{f}
+ \boldsymbol{g}_{1}^{n+1} - \boldsymbol{g}_{2}^{n+1}).
$$
Hence, by replacing the formulation above into $P_3$, and applying the polarization identity \eqref{polar} again, one obtains:
\begin{equation}
    \begin{aligned}
        P_3 
        &= \frac{L}{2}\left(\left\Vert\frac{\boldsymbol{e}_{f}^{n}\boldsymbol{n}_f}{2L}
         + \frac{\boldsymbol{e}_{s}^{n}\boldsymbol{n}_f}{2L}\right\Vert_{\Gamma}^{2}
         - \left\Vert\frac{\boldsymbol{e}_{f}^{n+1}\boldsymbol{n}_f}{L}\right\Vert_{\Gamma}^{2} 
         +\left\Vert\frac{\boldsymbol{e}_{f}^{n}\boldsymbol{n}_f}{2L} 
         + \frac{\boldsymbol{e}_{s}^{n}\boldsymbol{n}_f}{2L}
         -\frac{\boldsymbol{e}_{f}^{n+1}\boldsymbol{n}_f}{L}\right\Vert_{\Gamma}^{2}\right) \\
        &+\frac{L}{2}\left(\left\Vert\frac{\boldsymbol{e}_{f}^{n}\boldsymbol{n}_f}{2L}
         + \frac{\boldsymbol{e}_{s}^{n}\boldsymbol{n}_f}{2L}\right\Vert_{\Gamma}^{2}
         - \left\Vert\frac{\boldsymbol{e}_{s}^{n+1}\boldsymbol{n}_f}{L}\right\Vert_{\Gamma}^{2}
         + \left\Vert\frac{\boldsymbol{e}_{f}^{n}\boldsymbol{n}_f}{2L} 
         + \frac{\boldsymbol{e}_{s}^{n}\boldsymbol{n}_f}{2L}
         -\frac{\boldsymbol{e}_{s}^{n+1}\boldsymbol{n}_f}{L}\right\Vert_{\Gamma}^{2}\right) \\
        &+ \int_{\Gamma}\frac{1}{2}(\boldsymbol{e}_{f}^{n}\boldsymbol{n}_f +\boldsymbol{e}_{s}^{n}\boldsymbol{n}_f)\cdot\frac{1}{L}(\boldsymbol{g}_{1}^{n+1} - \boldsymbol{g}_{2}^{n+1}).
    \end{aligned}
    \label{p3}
\end{equation}
By adding \eqref{p1}, \eqref{p2}, and \eqref{p3} together, we observe that the third terms in each bracket cancel exactly due to \eqref{l1} and \eqref{l2}, thus yielding:
\begin{equation}
    \begin{aligned}
    &J^{n+1}
    = P_1 + P_2 + P_3 \\
    &= \frac{L}{2}\left(\left\Vert\frac{\boldsymbol{e}_{u}^n }{2}
    +\frac{\boldsymbol{e}_{\xi}^n}{2} 
    + \frac{\boldsymbol{g}_{1}^{n+1}}{L}\right\Vert_{\Gamma}^{2} 
    - \left\Vert\boldsymbol{e}_{u}^{n+1} \right\Vert_{\Gamma}^{2}
    +\left\Vert\frac{\boldsymbol{e}_{u}^n }{2}
    +\frac{\boldsymbol{e}_{\xi}^n}{2} 
    + \frac{\boldsymbol{g}_{2}^{n+1}}{L}\right\Vert_{\Gamma}^{2} 
    - \Vert\boldsymbol{e}_{\xi}^{n+1} \Vert_{\Gamma}^{2}\right)\\
    &+\frac{L}{2}\left(2\left\Vert\frac{\boldsymbol{e}_{f}^{n}\boldsymbol{n}_f}{2L}
    + \frac{\boldsymbol{e}_{s}^{n}\boldsymbol{n}_f}{2L}\right\Vert_{\Gamma}^{2}
    - \left\Vert\frac{\boldsymbol{e}_{f}^{n+1}\boldsymbol{n}_f}{L}\right\Vert_{\Gamma}^{2}
    - \left\Vert\frac{\boldsymbol{e}_{s}^{n+1}\boldsymbol{n}_f}{L}\right\Vert_{\Gamma}^{2}\right) \\
&+ \frac{1}{2L}\int_{\Gamma}(\boldsymbol{e}_{f}^{n}\boldsymbol{n}_{f} + 
\boldsymbol{e}_{s}^{n}\boldsymbol{n}_{f} )\cdot(\boldsymbol{g}_{1}^{n+1} - \boldsymbol{g}_{2}^{n+1}) \\
&=: Q_1 + Q_2 + Q_3. 
    \end{aligned}
    \end{equation}
By applying the Cauchy inequality and the Young's inequality, we have
    \begin{equation}
        \begin{aligned}
            Q_1 &= \frac{L}{2}\left( \frac{1}{2}\Vert\boldsymbol{e}_{u}^n + \boldsymbol{e}_{\xi}^{n}\Vert_{\Gamma}^{2}  
    + \frac{1}{L^2}(\Vert\boldsymbol{g}_{1}^{n+1}\Vert_{\Gamma}^{2} + \Vert\boldsymbol{g}_{2}^{n+1}\Vert_{\Gamma}^{2})
    + \frac{1}{L}\Vert(\boldsymbol{e}_{u}^{n} + \boldsymbol{e}_{\xi}^{n})\boldsymbol{g}_{1}^{n+1}\Vert_{\Gamma}\right.  \\
    &\left.+ \frac{1}{L}\Vert(\boldsymbol{e}_{u}^{n} + \boldsymbol{e}_{\xi}^{n})\boldsymbol{g}_{2}^{n+1}\Vert_{\Gamma} 
    - \Vert\boldsymbol{e}_{u}^{n+1}\Vert_{\Gamma}^{2} -\Vert\boldsymbol{e}_{\xi}^{n+1}\Vert_{\Gamma}^{2} \right) \\
    &\leq \frac{L}{2}\left(\Vert\boldsymbol{e}_{u}^{n}\Vert_{\Gamma}^{2} - \Vert\boldsymbol{e}_{u}^{n+1}\Vert_{\Gamma}^{2} + \Vert\boldsymbol{e}_{\xi}^{n}\Vert_{\Gamma}^{2} -\Vert\boldsymbol{e}_{\xi}^{n + 1}\Vert_{\Gamma}^{2} 
    + \frac{1}{L^2}\Vert\boldsymbol{g}_{1}^{n+1}\Vert_{\Gamma}^{2}
    + \frac{1}{L^2}\Vert\boldsymbol{g}_{2}^{n+1}\Vert_{\Gamma}^{2}\right. \\
    &\left.+ \delta\Vert\boldsymbol{e}_{u}^{n}\Vert_{\Gamma}^{2}
    + \delta\Vert\boldsymbol{e}_{\xi}^{n}\Vert_{\Gamma}^{2}
    + \frac{4}{\delta L^2}\Vert\boldsymbol{g}_{1}^{n+1}\Vert_{\Gamma}^{2}
    + \frac{4}{\delta L^2}\Vert\boldsymbol{g}_{2}^{n+1}\Vert_{\Gamma}^{2}
    \right),
        \end{aligned}
    \end{equation}
    \begin{equation}
        Q_2\le\frac{1}{2L}(\Vert\boldsymbol{e}_{f}^n\boldsymbol{n}_{f}\Vert_{\Gamma}^{2} + \Vert\boldsymbol{e}_{s}^n\boldsymbol{n}_{f}\Vert_{\Gamma}^{2})-\frac{1}{2L}\Vert\boldsymbol{e}_{f}^{n+1}\boldsymbol{n}_{f}\Vert_{\Gamma}^{2} - \frac{1}{2L}\Vert\boldsymbol{e}_{s}^{n+1}\boldsymbol{n}_{f}\Vert_{\Gamma}^{2},
    \end{equation}
    \begin{equation}
        Q_3 \le \frac{\delta}{4L}\Vert\boldsymbol{e}_{f}^{n}\boldsymbol{n}_{f} + \boldsymbol{e}_{s}^{n}\boldsymbol{n}_{f} \Vert_{\Gamma}^{2} + \frac{1}{4\delta L}\Vert\boldsymbol{g}_{1}^{n+1} -  \boldsymbol{g}_{2}^{n+1}\Vert_{\Gamma}^{2}.
    \end{equation}
Combining $Q_1, Q_2$ and $Q_3$ together, we obtain:
    \begin{equation}
        \begin{aligned}
            J^{n+1} &= Q_1 + Q_2 + Q_3 \\
            &\leq \frac{L}{2}\left[ (1 + \delta)\Vert\boldsymbol{e}_{u}^{n}\Vert_{\Gamma}^{2} +(1 + \delta) \Vert\boldsymbol{e}_{\xi}^{n}\Vert_{\Gamma}^{2} - \Vert\boldsymbol{e}_{u}^{n+1}\Vert_{\Gamma}^{2} - \Vert\boldsymbol{e}_{\xi}^{n+1}\Vert_{\Gamma}^{2}
            \right] \\
            &+ \frac{1}{2L}\left[(1 + \delta)\Vert\boldsymbol{e}_{f}^{n}\boldsymbol{n}_{f}\Vert_{\Gamma}^{2} + (1 + \delta)\Vert\boldsymbol{e}_{s}^{n}\boldsymbol{n}_{f}\Vert_{\Gamma}^{2}
            - \Vert\boldsymbol{e}_{f}^{n+1}\boldsymbol{n}_{f}\Vert_{\Gamma}^{2} - \Vert\boldsymbol{e}_{s}^{n+1}\boldsymbol{n}_{f}\Vert_{\Gamma}^{2}\right]\\
            &+ \frac{\delta + 5}{2\delta L}(\Vert\boldsymbol{g}_{1}^{n+1} \Vert_{\Gamma}^{2} + \Vert\boldsymbol{g}_{2}^{n+1} \Vert_{\Gamma}^{2}).
        \end{aligned}
    \end{equation}
We denote $\boldsymbol{G}^{n+1}$ as 
$$
\boldsymbol{G}^{n+1} = \frac{\delta + 5}{2\delta L}\int_{t_n}^{t_{n+1}}(\Vert\boldsymbol{g}_{1}^{n}(s)\Vert_{\Gamma}^{2}+ \Vert\boldsymbol{g}_{2}^{n}(s)\Vert_{\Gamma}^{2})ds.
$$
Taking the integral on $[t_{n}, t_{n+1}]$ we have:
\begin{equation}
        \boldsymbol{E}^{n+1} + \boldsymbol{T}^{n+1} + \boldsymbol{S}^{n+1}
\leq \boldsymbol{E}^{n} + \boldsymbol{S}^{n} + \delta\boldsymbol{S}^{n} + \boldsymbol{G}^{n+1}.
\end{equation}
It then follows directly that
\begin{equation}
        \boldsymbol{E}^{n+1} + \boldsymbol{T}^{n+1} + \boldsymbol{S}^{n+1}
\leq \boldsymbol{E}^{n} + \boldsymbol{S}^{n} + \delta\max_{0\le m\le N}\boldsymbol{S}^{m} + \boldsymbol{G}^{n+1}.
\end{equation}
By summing from $n = 0$ to $n = M - 1$ with $M\leq N$, and setting $\delta =\frac{\Delta t}{2T} $, one ultimately yields
\begin{equation}
    \begin{aligned}
         \boldsymbol{E}^{M} + \sum_{n=0}^{M-1}\boldsymbol{T}^{n+1} + \boldsymbol{S}^{M}
\leq  \frac{1}{2}\underset{0\leq m\leq N}{\max}\boldsymbol{S}^{m} + \sum_{n=0}^{M-1}\boldsymbol{G}^{n+1}.
    \end{aligned}
\end{equation}
Since this holds for any $1\leq M \leq N$, we have:
$$
\frac{1}{2}\underset{0\leq m\leq N}{\max}\boldsymbol{S}^{m} \leq \sum_{n=0}^{N-1}\boldsymbol{G}^{n+1}
$$
Hence, we have:
\begin{equation*}
    \begin{aligned}
        \boldsymbol{E}^{N} + \sum_{n=0}^{N-1}\boldsymbol{T}^{n+1} + \boldsymbol{S}^{N}
\leq  2\sum_{n=0}^{N-1}\boldsymbol{G}^{n+1}.
    \end{aligned}
\end{equation*}
Applying Lemma \ref{lemma} and that $\delta=\frac{\Delta t}{2T}$, we have the error estimates \eqref{error result}.
\end{proof}
\begin{remark}
    The sub-optimal order of temporal convergence is often obtained in {\color{black} loosely coupled splitting} schemes for the interaction between fluids and thick structures. In particular, sub-optimal order has been studied in \cite{Guzman2020, seboldt2021noniterative} for Robin-Robin methods, in \cite{fernandez2015generalized} for Robin-Neumann methods, and in \cite{burman2009stabilization} for the {\color{black}loosely coupled splitting} scheme based on Nitsche's method. For the interaction between fluids and thin structures (lower-dimensional model), optimal analytical results are studied in \cite{li2024optimal}. {\color{black} Recently, a loosely coupled splitting scheme for FSI problems was analyzed in \cite{burman2025robin}, and a nearly optimal convergence result was established under higher regularity and additional assumptions.}
\end{remark}

\begin{remark}
Regarding the value of the parameter $L$, the splitting error scales as $\sqrt{L + L^{-1}}\sqrt{\Delta t}$. Consequently, we anticipate that a smaller $\Delta t$ is required for larger $L$ to balance it, so that the convergence rate can be observed. More details will be shown in the numerical experiment section \ref{manufactured}.
\end{remark}
\subsection{\textcolor{black}{Discrete stability analysis}}
\label{Stability analysis}
In this section, we will discretize problems \eqref{linear_f_robin} and \eqref{linear_s_robin} using the same temporal discretization scheme, i.e. the Backward Euler scheme, described in Algorithm \ref{semi}.
\textcolor{black}{For the energy argument only, we temporarily write the parameter
$L$ as $L_1$ and $L_2$ in the fluid and structure subproblems, respectively;
the analysis below ultimately requires $L_1=L_2=L$.}
Specifically, the fluid subproblem is formulated as follows: find $\boldsymbol{u}^{n+1}$ and $p^{n+1}$ such that
\begin{subequations}
  \begin{align}
    &\rho_f\frac{\boldsymbol{u}^{n+1}-\boldsymbol{u}^n}{\Delta t}-\nabla\cdot\boldsymbol{\sigma}_f^{n+1}=0,&&\text{in}\ \Omega_f,\\
    &\nabla\cdot\boldsymbol{u}^{n+1}=0,&&\text{in}\ \Omega_f,\\
    &L_1\boldsymbol{u}^{n+1}+\boldsymbol{\sigma}_f^{n+1}\boldsymbol{n}_f=\frac{L_1}{2}\boldsymbol{u}^n+\frac{L_1}{2}\boldsymbol{\xi}^n+\frac{1}{2}\boldsymbol{\sigma}_f^n\boldsymbol{n}_f +\frac{1}{2}\boldsymbol{\sigma}_s^{n}\boldsymbol{n}_f, &&\text{on }\Gamma;\label{robin1}
  \end{align}
  \label{linear_f_robin_d}
\end{subequations}
and the structure subproblem is: find $\boldsymbol{\xi}^{n+1}=(\boldsymbol{\eta}^{n+1}-\boldsymbol{\eta}^n)/\Delta t$ such that
\begin{subequations}
\begin{align}
  &\rho_s\frac{\boldsymbol{\xi}^{n+1}-\boldsymbol{\xi}^n}{\Delta t} - \nabla\cdot\boldsymbol{\sigma}_s^{n+1}=0,&&\text{in }\Omega_s,\\
  &L_2\boldsymbol{\xi}^{n+1}+\boldsymbol{\sigma}_s^{n+1}\boldsymbol{n}_s=\frac{L_2}{2}\boldsymbol{u}^n+\frac{L_2}{2}\boldsymbol{\xi}^n+\frac{1}{2}\boldsymbol{\sigma}_f^n\boldsymbol{n}_s +\frac{1}{2}\boldsymbol{\sigma}_s^{n}\boldsymbol{n}_s,&&\text{on }\Gamma.\label{robin2}
\end{align}
\label{linear_s_robin_d}
\end{subequations}
Again, we remark that the same notation $\cdot^n$ will be used to denote its discrete approximation $\cdot^n(t_n)$ at the time $t_n$.
To state the stability result, the following quantities are introduced:
\begin{equation*}
  \begin{aligned}
      \mathcal E^n&=\frac{\rho_s}{2}\Vert\boldsymbol{\xi}^n\Vert^2_s+\frac{1}{2}\Vert\boldsymbol{\eta}^n\Vert_E^2+\frac{\rho_f}{2}\Vert\boldsymbol{u}^n\Vert^2_f,\\
      \mathcal D^n&=\mu_f\Delta t\Vert\mathbb D(\boldsymbol{u}^n)\Vert^2_f +\frac{\rho_s}{2}\Vert\boldsymbol{\xi}^n-\boldsymbol{\xi}^{n-1}\Vert^2_s +\frac{1}{2}\Vert\boldsymbol{\eta}^n-\boldsymbol{\eta}^{n-1}\Vert_E^2 + \frac{\rho_f}{2}\Vert\boldsymbol{u}^{n}-\boldsymbol{u}^{n-1}\Vert^2_f,\\
      \mathcal I^n&=\frac{L_1\Delta t}{2}\Vert\boldsymbol{u}^{n}\Vert_\Gamma^2 + \frac{L_2\Delta t}{2}\Vert\boldsymbol{\xi}^{n}\Vert_\Gamma^2+\frac{\Delta t}{2L_1}\Vert\boldsymbol{\sigma}_f^{n}\boldsymbol{n}_f\Vert_\Gamma^2+\frac{\Delta t}{2L_2}\Vert\boldsymbol{\sigma}_s^{n}\boldsymbol{n}_s\Vert_\Gamma^2.
  \end{aligned}
\end{equation*}

\begin{theorem}
\label{stability}
Let $(\boldsymbol{u}^{n+1}, p^{n+1}, \boldsymbol{\xi}^{n+1}, \boldsymbol{\eta}^{n+1})$ denote the solution of the decoupled scheme given by (\ref{linear_f_robin_d}) and (\ref{linear_s_robin_d}) for $0 \leq n \leq N-1$. The scheme is unconditionally stable when $L_1 = L_2$, and the following a priori energy estimate holds:
\begin{equation}
  \mathcal{E}^N + \sum_{n=1}^{N}\mathcal{D}^n + \mathcal{I}^N \le \mathcal{E}^0 + \mathcal{I}^0.
\end{equation}
\end{theorem}
\begin{proof}
  By taking $(\boldsymbol{v}, q, \boldsymbol{\zeta})=\Delta t(\boldsymbol{u}^{n+1}, p^{n+1}, \boldsymbol{\xi}^{n+1})$ in (\ref{linear_f_robin_d}) and (\ref{linear_s_robin_d}), adding the equations together, integrating by parts and applying (\ref{polar}), one yields:
\begin{equation}
\begin{aligned}
        &\frac{\rho_f}{2}\left(\Vert\boldsymbol{u}^{n+1}\Vert^2_f-\Vert\boldsymbol{u}^n\Vert^2_f+\Vert\boldsymbol{u}^{n+1}-\boldsymbol{u}^n\Vert^2_f\right) + 2\mu_f\Delta t\Vert\mathbb D(\boldsymbol{u}^{n+1})\Vert^2_f\\
        &+\frac{\rho_s}{2}\left(\Vert\boldsymbol{\xi}^{n+1}\Vert^2_s-\Vert\boldsymbol{\xi}^n\Vert^2_s+\Vert\boldsymbol{\xi}^{n+1}-\boldsymbol{\xi}^n\Vert^2_s\right)\\
        &+\mu_s\left(\Vert\mathbb D(\boldsymbol{\eta}^{n+1})\Vert^2_s-\Vert\mathbb D(\boldsymbol{\eta}^n)\Vert^2_s+\Vert\mathbb D(\boldsymbol{\eta}^{n+1}-\boldsymbol{\eta}^n)\Vert^2_s\right)\\
        &+\frac{\lambda_s}{2}\left(\Vert\nabla\cdot\boldsymbol{\eta}^{n+1}\Vert^2_s-\Vert\nabla\cdot\boldsymbol{\eta}^n\Vert^2_s+\Vert\nabla\cdot(\boldsymbol{\eta}^{n+1}-\boldsymbol{\eta}^n)\Vert^2_s\right)\\
        &+L_1\Delta t\left<\boldsymbol{u}^{n+1}-\frac{1}{2}\boldsymbol{\xi}^n-\frac{1}{2}\boldsymbol{u}^n,\boldsymbol{u}^{n+1}\right>_\Gamma+L_2\Delta t\left<\boldsymbol{\xi}^{n+1}-\frac{1}{2}\boldsymbol{\xi}^n-\frac{1}{2}\boldsymbol{u}^n,\boldsymbol{\xi}^{n+1}\right>_\Gamma\\
        &=\frac{\Delta t}{2}\left<\boldsymbol{\sigma}_f^n\boldsymbol{n}_f+\boldsymbol{\sigma}_s^n\boldsymbol{n}_f, \boldsymbol{u}^{n+1}\right>_\Gamma +\frac{\Delta t}{2}\left<\boldsymbol{\sigma}_f^n\boldsymbol{n}_s+\boldsymbol{\sigma}_s^n\boldsymbol{n}_s, \boldsymbol{\xi}^{n+1}\right>_\Gamma.
\end{aligned}
\label{weakform}
\end{equation}
For notational simplicity, let $\mathcal{F}^n=\boldsymbol{\sigma}_f^n\boldsymbol{n}_f$ and $\mathcal{S}^n=\boldsymbol{\sigma}_s^n\boldsymbol{n}_f$. Subtracting (\ref{robin2}) from (\ref{robin1}) yields:
\begin{equation}
  \boldsymbol{u}^{n+1}-\boldsymbol{\xi}^{n+1} = \left(\frac{1}{2L_1}+\frac{1}{2L_2}\right)\left(\mathcal{F}^n+\mathcal{S}^n\right)-\frac{\mathcal{F}^{n+1}}{L_1}-\frac{\mathcal{S}^{n+1}}{L_2}.
  \label{important}
\end{equation}
By applying (\ref{important}), the RHS of (\ref{weakform}) can be rewritten as:
\begin{equation}
  \begin{aligned}
      &\frac{\Delta t}{2}\left<\mathcal{F}^{n}+\mathcal{S}^n, \boldsymbol{u}^{n+1}-\boldsymbol{\xi}^{n+1}\right>_\Gamma\\
      &=\frac{\Delta t}{2}\left<\mathcal{F}^n+\mathcal{S}^n,\left(\frac{1}{2L_1}+\frac{1}{2L_2}\right)\left(\mathcal{F}^n+\mathcal{S}^n\right)-\frac{\mathcal{F}^{n+1}}{L_1}-\frac{\mathcal{S}^{n+1}}{L_2}\right>_\Gamma.
  \end{aligned}
  \label{t1}
\end{equation}
Using (\ref{polar}) again, the right-hand side above can be further rewritten as:
  \begin{align}
  &\frac{\Delta t}{2}\left<\mathcal{F}^n+\mathcal{S}^n, \frac{1}{2L_1}\left(\mathcal{F}^n+\mathcal{S}^n\right)-\frac{1}{L_1}\mathcal{F}^{n+1}\right>_\Gamma\nonumber\\
  +&\frac{\Delta t}{2}\left<\mathcal{F}^n+\mathcal{S}^n, \frac{1}{2L_2}\left(\mathcal{F}^n+\mathcal{S}^n\right)-\frac{1}{L_2}\mathcal{S}^{n+1}\right>_\Gamma\nonumber\\
  =&\frac{L_1\Delta t}{2}\left[\frac{\Vert\mathcal{F}^n+\mathcal{S}^n\Vert^2_\Gamma}{4L_1^2}-\frac{\Vert\mathcal{F}^{n+1}\Vert^2_\Gamma}{L_1^2}+\left\Vert\frac{\mathcal{F}^n+\mathcal{S}^n}{2L_1}-\frac{\mathcal{F}^{n+1}}{L_1}\right\Vert^2_\Gamma\right]\nonumber\\
  +&\frac{L_2\Delta t}{2}\left[\frac{\Vert\mathcal{F}^n+\mathcal{S}^n\Vert^2_\Gamma}{4L_2^2}-\frac{\Vert\mathcal{S}^{n+1}\Vert^2_\Gamma}{L_2^2}+\left\Vert\frac{\mathcal{F}^n+\mathcal{S}^n}{2L_2}-\frac{\mathcal{S}^{n+1}}{L_2}\right\Vert^2_\Gamma\right].\label{t3}
  \end{align}
Notice from (\ref{robin1}) and (\ref{robin2}) that the third terms in both square brackets are respectively equal to $\Vert\boldsymbol{u}^{n+1}-\frac{1}{2}\boldsymbol{\xi}^n-\frac{1}{2}\boldsymbol{u}^n\Vert^2_\Gamma$ and $\Vert\boldsymbol{\xi}^{n+1}-\frac{1}{2}\boldsymbol{\xi}^n-\frac{1}{2}\boldsymbol{u}^n\Vert^2_\Gamma$ based on the Robin boundary conditions defined in \eqref{robin1} and \eqref{robin2}. For the rest of terms above, by using Cauchy inequality we have:
\begin{equation}
  \begin{aligned}
    \Delta t\left(\frac{1}{8L_1}+\frac{1}{8L_2}\right)\Vert\mathcal{F}^n+\mathcal{S}^n\Vert^2_\Gamma-\frac{\Delta t}{2L_1}\Vert\mathcal{F}^{n+1}\Vert_\Gamma^2-\frac{\Delta t}{2L_2}\Vert\mathcal{S}^{n+1}\Vert_\Gamma^2\\
    \le\left(\frac{\Delta t}{4L_1}+\frac{\Delta t}{4L_2}\right)\left(\Vert\mathcal{F}^{n}\Vert_\Gamma^2+\Vert\mathcal{S}^{n}\Vert_\Gamma^2\right)-\frac{\Delta t}{2L_1}\Vert\mathcal{F}^{n+1}\Vert_\Gamma^2-\frac{\Delta t}{2L_2}\Vert\mathcal{S}^{n+1}\Vert_\Gamma^2.
  \end{aligned}
  \label{t2}
\end{equation}
To derive unconditionally stable energy estimates, we here require:
\begin{equation}
  \begin{aligned}
  \frac{1}{4L_1}+\frac{1}{4L_2}\le\frac{1}{2L_1}&\qquad\text{and }\qquad\frac{1}{4L_1}+\frac{1}{4L_2}\le\frac{1}{2L_2},
\end{aligned}
\end{equation}
which implies that
$$
\Rightarrow L_1=L_2.
$$
Namely, the variables from adjacent time steps cancel out with each other during summation from 0 to $N-1$.

For the last two boundary integral terms on the left-hand side of (\ref{weakform}), by using (\ref{polar}) one yields:
\begin{align*}
  &L_1\Delta t\left<\boldsymbol{u}^{n+1}-\frac{1}{2}\boldsymbol{\xi}^n-\frac{1}{2}\boldsymbol{u}^n,\boldsymbol{u}^{n+1}\right>_\Gamma+L_2\Delta t\left<\boldsymbol{\xi}^{n+1}-\frac{1}{2}\boldsymbol{\xi}^n-\frac{1}{2}\boldsymbol{u}^n,\boldsymbol{\xi}^{n+1}\right>_\Gamma\\
  &=\frac{L_1\Delta t}{2}\left[\Vert\boldsymbol{u}^{n+1}\Vert^2_\Gamma-\frac{1}{4}\Vert\boldsymbol{u}^n+\boldsymbol{\xi}^n\Vert^2_\Gamma+\Vert\boldsymbol{u}^{n+1}-\frac{1}{2}\boldsymbol{\xi}^n-\frac{1}{2}\boldsymbol{u}^n\Vert^2_\Gamma\right]\\
  &+\frac{L_2\Delta t}{2}\left[\Vert\boldsymbol{\xi}^{n+1}\Vert^2_\Gamma-\frac{1}{4}\Vert\boldsymbol{u}^n+\boldsymbol{\xi}^n\Vert^2_\Gamma+\Vert\boldsymbol{\xi}^{n+1}-\frac{1}{2}\boldsymbol{\xi}^n-\frac{1}{2}\boldsymbol{u}^n\Vert^2_\Gamma\right].
\end{align*}
Notice that the last two terms in both square brackets above would be canceled by corresponding terms in (\ref{t3}). For the rest of the terms above, by using the Cauchy inequality and $L_1=L_2$, one yields:
\begin{equation}
\begin{aligned}
  &\Delta t\left(\frac{L_1}{2}\Vert\boldsymbol{u}^{n+1}\Vert_\Gamma^2-\frac{L_1}{8}\Vert\boldsymbol{u}^n+\boldsymbol{\xi}^n\Vert^2_\Gamma+\frac{L_2}{2}\Vert\boldsymbol{\xi}\Vert_\Gamma^2-\frac{L_2}{8}\Vert\boldsymbol{u}^{n}+\boldsymbol{\xi}^n\Vert_\Gamma^2\right)\\
  &\ge\frac{L_1\Delta t}{2}\left(\Vert\boldsymbol{u}^{n+1}\Vert_\Gamma^2-\Vert\boldsymbol{u}^n\Vert_\Gamma^2\right)+\frac{L_2\Delta t}{2}\left(\Vert\boldsymbol{\xi}^{n+1}\Vert_\Gamma^2-\Vert\boldsymbol{\xi}^n\Vert_\Gamma^2\right).
\end{aligned}
\label{t4}
\end{equation}
Finally, combining all estimates above, we have:
$$
\mathcal{E}^{n+1} + \mathcal{D}^{n+1} + \mathcal{I}^{n+1} \le \mathcal{E}^{n} + \mathcal{I}^n.
$$
The result follows after summing the inequalities over $n$ from 0 to $N-1$.
\end{proof}

\begin{remark}
The condition $L_1=L_2$ is only a sufficient condition in our analysis. 
{\color{black}For $L_1\neq L_2$, stability is not established by the present analysis and
remains an open question.}

From another point of view,  (\ref{robin1})  and (\ref{robin2}) can be rewritten as
\begin{subequations}
        \begin{align}
            L_{1}\Delta t\frac{\boldsymbol{u}^{n+1} - \boldsymbol{u}^{n}}{\Delta t} + \boldsymbol{\sigma}_{f}^{n+1}\boldsymbol{n}_{f} = \frac{L_{1}}{2}(\boldsymbol{\xi}^{n} - \boldsymbol{u}^{n}) + \frac{1}{2}\boldsymbol{\sigma}_{f}^{n}\boldsymbol{n}_{f} + \frac{1}{2}\boldsymbol{\sigma}_{s}^{n}\boldsymbol{n}_{f} ,\label{robin11}
            \\
             L_{2}\Delta t\frac{\boldsymbol{\xi}^{n+1} - \boldsymbol{\xi}^{n}}{\Delta t} + \boldsymbol{\sigma}_{s}^{n+1}\boldsymbol{n}_{s} = \frac{L_{2}}{2}(\boldsymbol{u}^{n} - \boldsymbol{\xi}^{n}) + \frac{1}{2}\boldsymbol{\sigma}_{f}^{n}\boldsymbol{n}_{s} + \frac{1}{2}\boldsymbol{\sigma}_{s}^{n}\boldsymbol{n}_{s} \label{robin22} .
        \end{align}
\end{subequations}
Equations (\ref{robin11})  and (\ref{robin22}) can be considered as  approximation for the continuous equations defined on the interface $\Gamma$:
$$
L_{1}\Delta t\partial_{t}\boldsymbol{u}^{n}  
= \frac{L_{1}}{2}(\boldsymbol{\xi}^{n} - \boldsymbol{u}^{n}) 
+ (\boldsymbol{\sigma}_{f}^{n}\boldsymbol{n}_{f} 
-  \boldsymbol{\sigma}_{f}^{n+1}\boldsymbol{n}_{f} ),
$$
$$
L_{2}\Delta t\partial_{t}\boldsymbol{\xi}^{n} 
= \frac{L_{2}}{2}(\boldsymbol{u}^{n} - \boldsymbol{\xi}^{n}) 
+ (\boldsymbol{\sigma}_{s}^{n}\boldsymbol{n}_{s} 
- \boldsymbol{\sigma}_{s}^{n+1}\boldsymbol{n}_{s} ).
$$
When $L_{1} = L_{2}=L$, subtracting the second equation from the first one, we obtain:
\begin{equation}
L\Delta t\partial_{t}(\boldsymbol{u}^{n}  - \boldsymbol{\xi}^{n} ) 
+ L(\boldsymbol{u}^{n}  - \boldsymbol{\xi}^{n} )
= 2(\boldsymbol{\sigma}_{f}^{n}\boldsymbol{n}_{f} 
- \boldsymbol{\sigma}_{f}^{n+1}\boldsymbol{n}_{f})
= \mathcal{O}(\Delta t).
\end{equation}
Roughly speaking, if  $\boldsymbol{u} - \boldsymbol{\xi} = 0$ initially at $t = t_{0}$, it provides an approximation of the original interface condition $\boldsymbol{u} = \boldsymbol{\xi}$ at any time $t$ by $L\Delta t\partial_{t}(\boldsymbol{u}  - \boldsymbol{\xi}) + L(\boldsymbol{u}  - \boldsymbol{\xi})= \mathcal{O}(\Delta t)$ then $\boldsymbol{u} - \boldsymbol{\xi} =  - e^{-(t - t_{0})/\Delta t} + C_{1}\mathcal{O}(\Delta t) / L$. Obviously, this approximation error will not increase with time, and this is due to the term $\Delta t\partial_{t}(\boldsymbol{u} - \boldsymbol{\xi})$ introduced from our Robin-Robin interface conditions when $L_{1} = L_{2}$.
\end{remark}

\section{Numerical \textcolor{black}{experiments}}
\label{section5}

This section presents a series of numerical examples that illustrate the stability 
and accuracy of the proposed method. The spatial discretization is performed using 
the finite element method with Taylor--Hood ($P_2-P_1$) elements for the fluid variables 
and ($P_2-P_2$) elements for the structure variables.
\textcolor{black}{
All computations are implemented in FEniCS~\cite{AlnaesEtal2014,AlnaesEtal2015},
with parallel execution based on \textit{mpi4py} and the built-in MPI support of FEniCS.
The nonlinear problems in Sections~5.2--5.4 are solved by
Newton's method with absolute and relative tolerances of $10^{-8}$, and the resulting
linear systems are solved directly using MUMPS without preconditioning.
For comparisons with strong coupling, we use the monolithic FSI solver
\textit{turtleFSI}~\cite{bergersen2020turtlefsi}, implemented in the same FEniCS framework and using MUMPS without preconditioning as well. The Robin--Robin and
Strong-Coupling computations use the same mesh, finite element spaces, time
discretization, and nonlinear-solver tolerances for each comparison.
For wall-clock measurements, each computation is repeated
three times, and the median time is reported.
Implementations of the numerical benchmarks are available at
\url{https://github.com/hnkoshi/A-fully-parallelizable-Robin-Robin-scheme-for-FSI.git}.
}

\subsection{Manufactured solution test}
\label{manufactured}

\begin{figure}[htbp]
    \centering

    \begin{subfigure}[b]{0.40\textwidth}
        \centering
        \includegraphics[width=\linewidth]{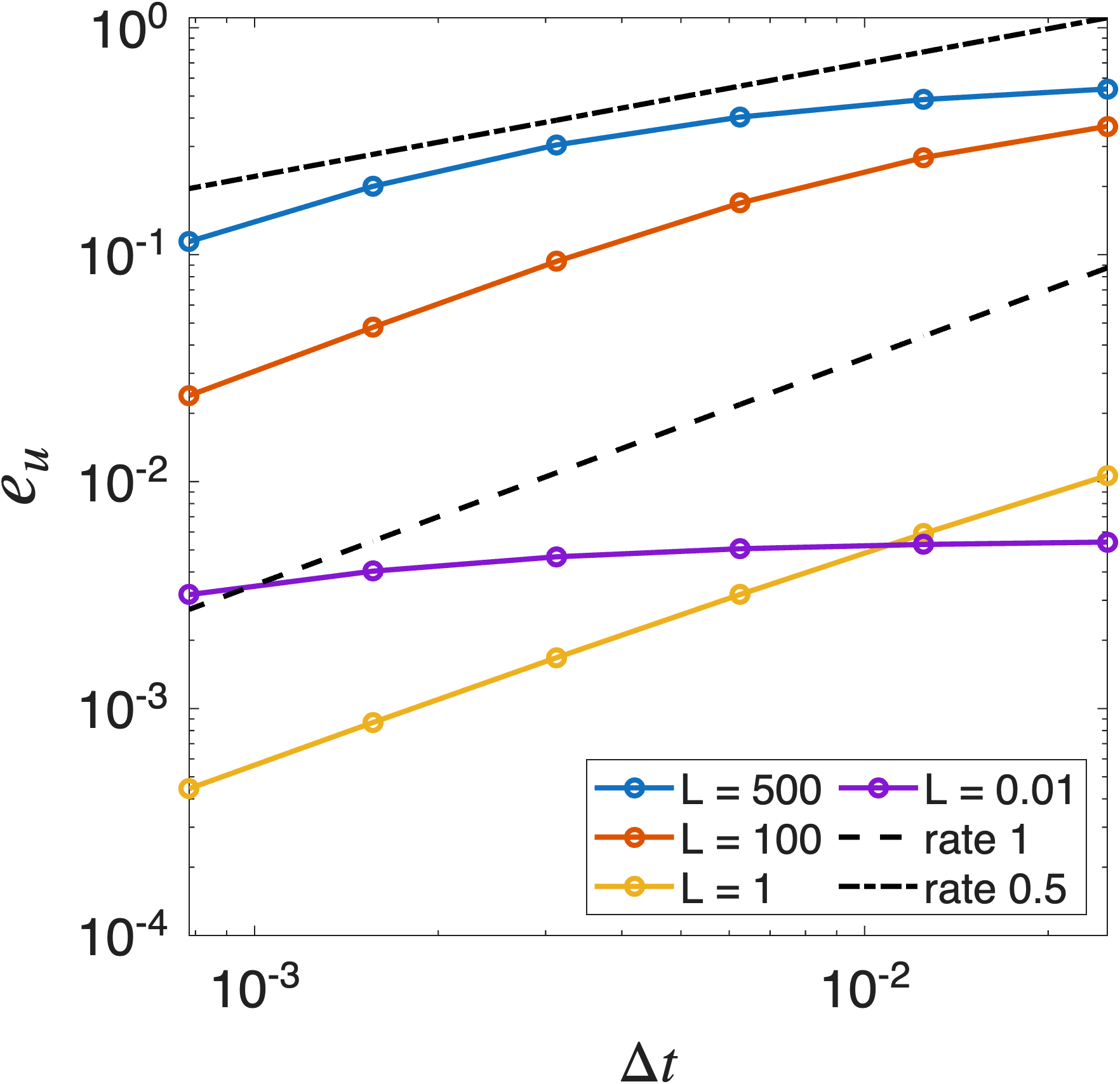}
        \label{fig_ex2:subfig1}
    \end{subfigure}
    \hspace{0.04\textwidth}
    \begin{subfigure}[b]{0.40\textwidth}
        \centering
        \includegraphics[width=\linewidth]{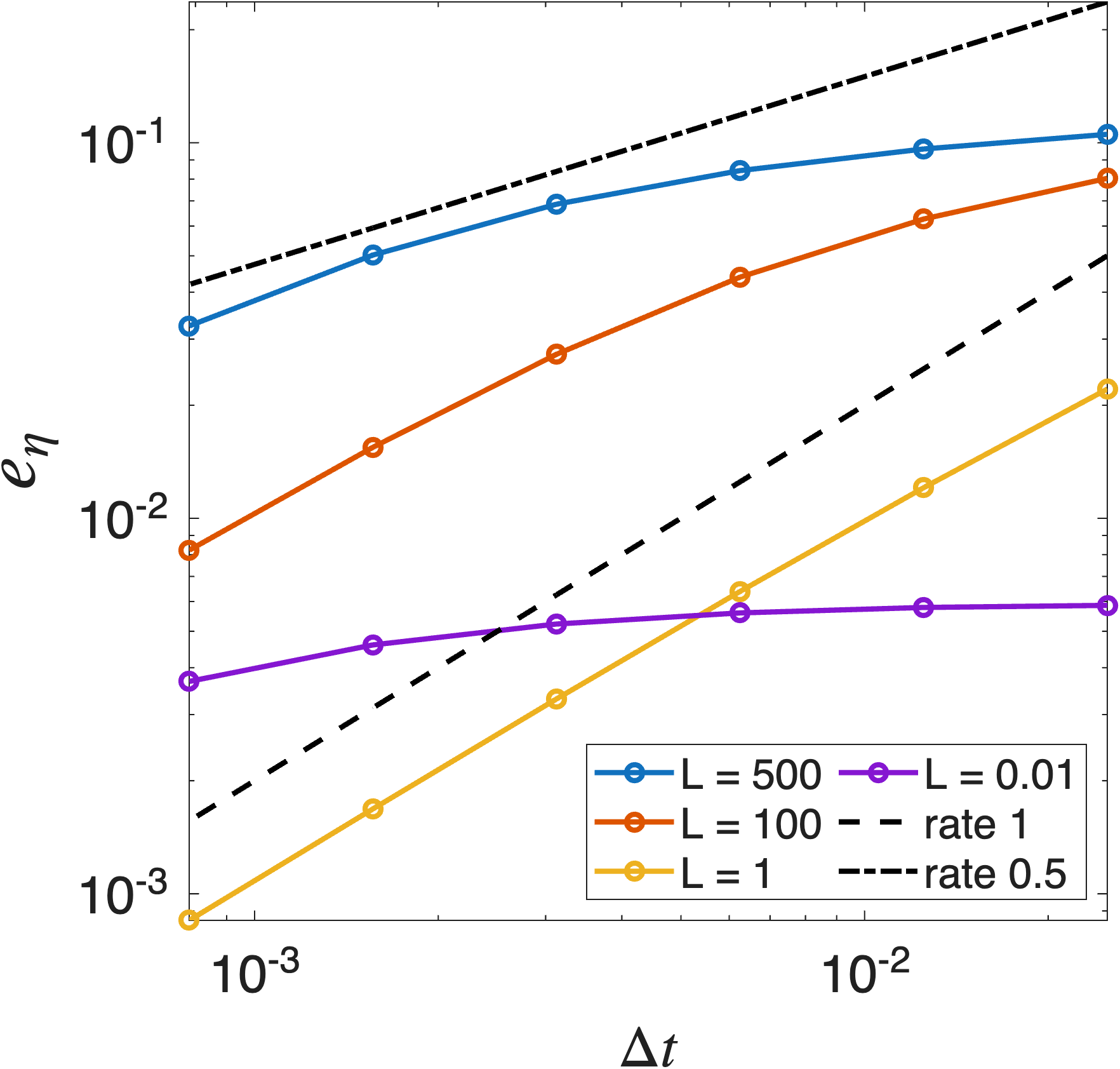}
        \label{fig_ex2:subfig2}
    \end{subfigure}

    \vspace{0.10cm}

    \begin{subfigure}[b]{0.40\textwidth}
        \centering
        \includegraphics[width=\linewidth]{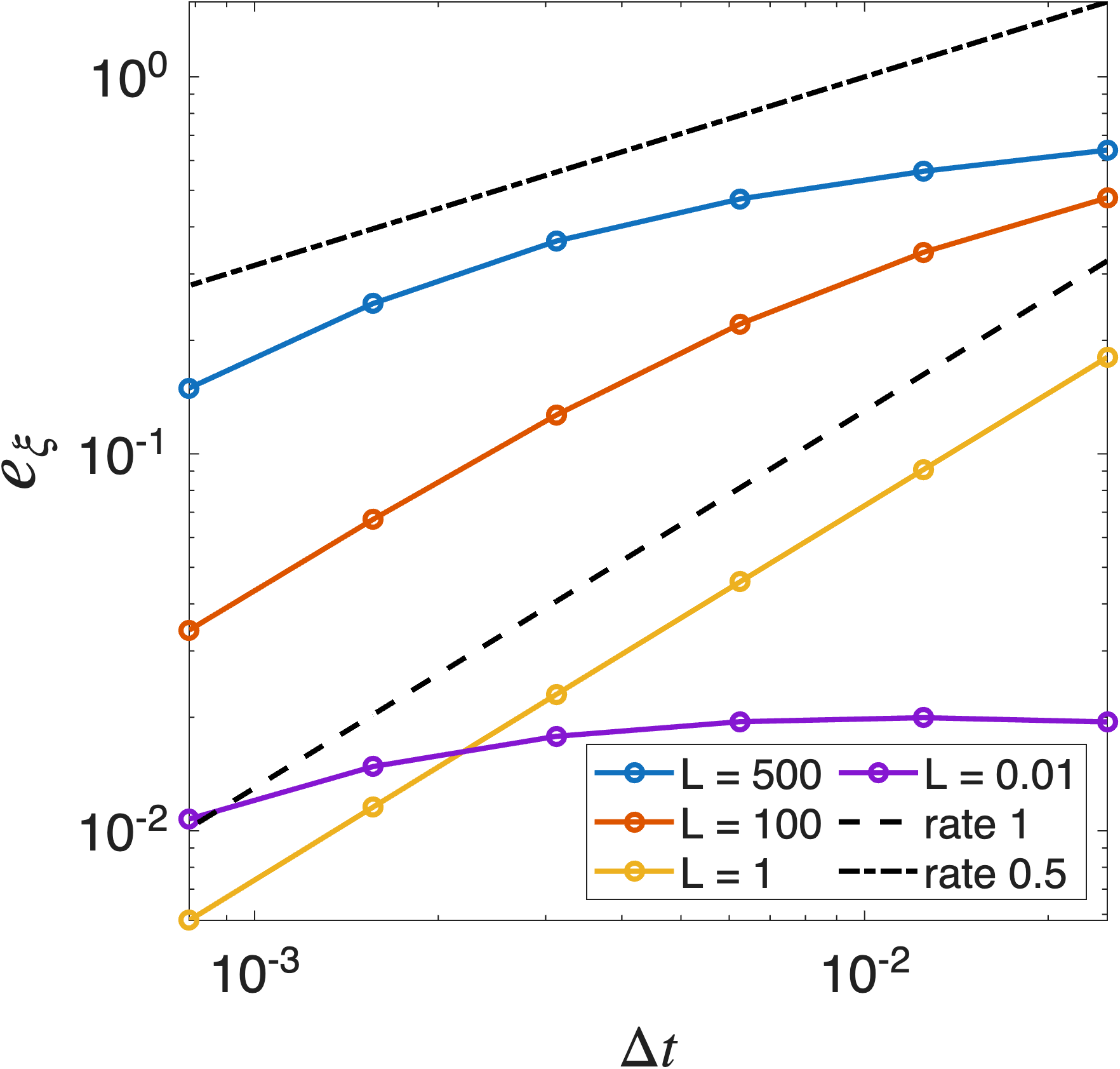}
        \label{fig_ex2:subfig4}
    \end{subfigure}
    \hspace{0.04\textwidth}
    \begin{subfigure}[b]{0.40\textwidth}
        \centering
        \includegraphics[width=\linewidth]{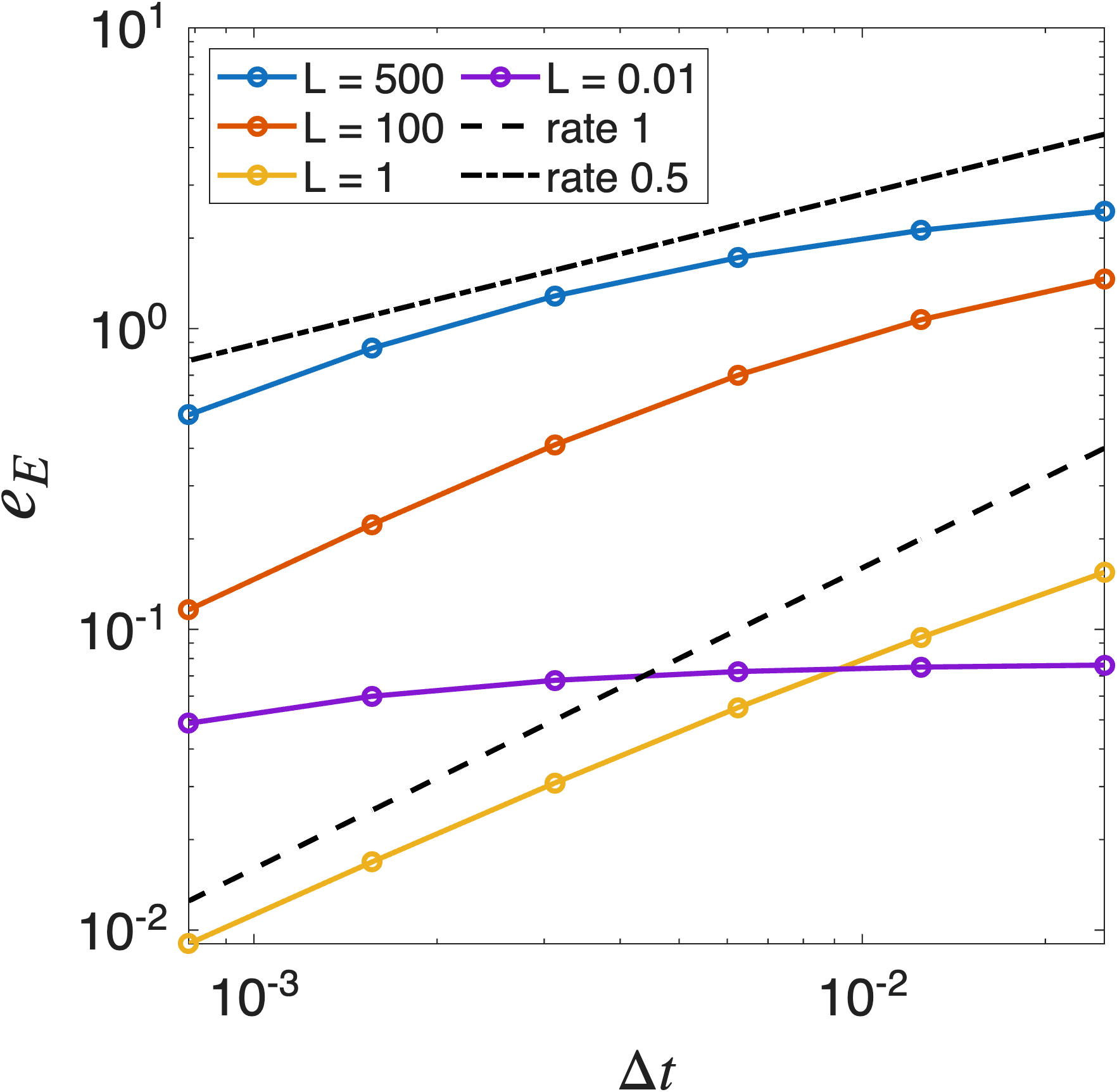}
        \label{fig_ex2:subfig5}
    \end{subfigure}

    \caption{
    Convergence plots for $\boldsymbol{e}_u$ (top-left),
    $\boldsymbol{e}_{\eta}$ (top-right),
    $\boldsymbol{e}_{\xi}$ (bottom-left), and
    $\boldsymbol{e}_E$ (bottom-right) with different Robin parameters
    $L$ at the final time $T=0.5$.
    }
    \label{fig_EX2:five-subfigures}
\end{figure}

The precision of the proposed decoupled scheme is assessed by comparison with a manufactured closed form solution of the linearized system \eqref{linear_fsi}. The following manufactured solution is defined in the composite domain consisting of a fluid region $\Omega_{f}=[0, 1]\times[0, 1]$ and a structure region $\Omega_{s}=[0, 1]\times[-1, 0]$, which are separated by the fluid-structure interface $\Gamma$.
\begin{align*}
    \boldsymbol{\eta} & = \begin{bmatrix}
    \sin(\pi t)(\cos(y) - 3x) \\
    \sin(\pi t)(y + 1)
    \end{bmatrix}, &
    \boldsymbol{\xi} & = \frac{\partial \boldsymbol{\eta}_{1}}{\partial t},\\
    \boldsymbol{u} & = \begin{bmatrix}
    \pi\cos(\pi t)(\cos(y) - 3x) \\
    \pi\cos(\pi t)(y + 1)
    \end{bmatrix}, &
    p & = 2\pi\cos(\pi t).
\end{align*}
 In particular, the model parameters are chosen as
\(\rho_{f} = \mu_{f} = \rho_{s} = \mu_{s} = \lambda_{s} = 1.0\). It should be noted that the prescribed fluid velocity \(\boldsymbol{u}\) does not satisfy the divergence-free condition. Therefore, the mass conservation equation is modified by replacing it with \(\nabla \cdot \boldsymbol{u}_{1} = f_m\). Dirichlet boundary conditions are imposed on all external boundaries, and the simulation is carried out until the final time \(T=0.5\). We define the following error norms for relevant variables
\begin{align*}
    &\boldsymbol{\mathbf{e}}_{u} = \Vert\boldsymbol{U}(T) - \boldsymbol{u}^N\Vert_f,\quad
    &&\boldsymbol{\mathbf{e}}_{\eta} = \Vert\mathbb{\boldsymbol{\Psi}}(T) - \boldsymbol{\eta}^N\Vert_s ,\\
    &\boldsymbol{\mathbf{e}}_{\xi} = \Vert\mathbb{\boldsymbol{\Xi}}(T) - \boldsymbol{\xi}^N\Vert_s,\quad
    &&\boldsymbol{\mathbf{e}}_{E} = \Vert\mathbb{\boldsymbol{\Psi}}(T) - \boldsymbol{\eta}^N\Vert_{E}.
\end{align*} 
The time step sizes $\Delta t$ are set as $0.1/n, \text{where}~ n = 4, 8, 16, 32, 64, 128$, while the mesh size is fixed, namely, $\Delta h=0.01$. The resulting numerical errors for different values of $L$ are reported in Figure \ref{fig_EX2:five-subfigures}.

{\color{black}
We observe that, for $L=1$ and $L=100$, the convergence rates of all
error measures are close to $\mathcal{O}(\Delta t)$. This indicates that the
theoretical $\mathcal{O}(\sqrt{\Delta t})$ estimate is a rigorous upper bound
and is not necessarily sharp for these parameter values. In contrast, for the
more extreme choices $L=0.01$ and $L=500$, the convergence behavior
deteriorates toward the theoretical half-order rate. In particular, for
$L=500$, the observed rate gradually improves as the time step is refined.
This behavior is consistent with Theorem~\ref{Error estimate}, where the
splitting error contains the factor
\[
\sqrt{L+L^{-1}}\sqrt{\Delta t}.
\]
Therefore, both excessively small and excessively large values of $L$ can
amplify the splitting-error constant, suggesting that extreme choices of the
Robin parameter should be avoided in practice.
}


\subsection{The Turek \& Hron's problem}
\label{benchmark}

\begin{figure}
    \centering
    \includegraphics[width=0.8\linewidth]{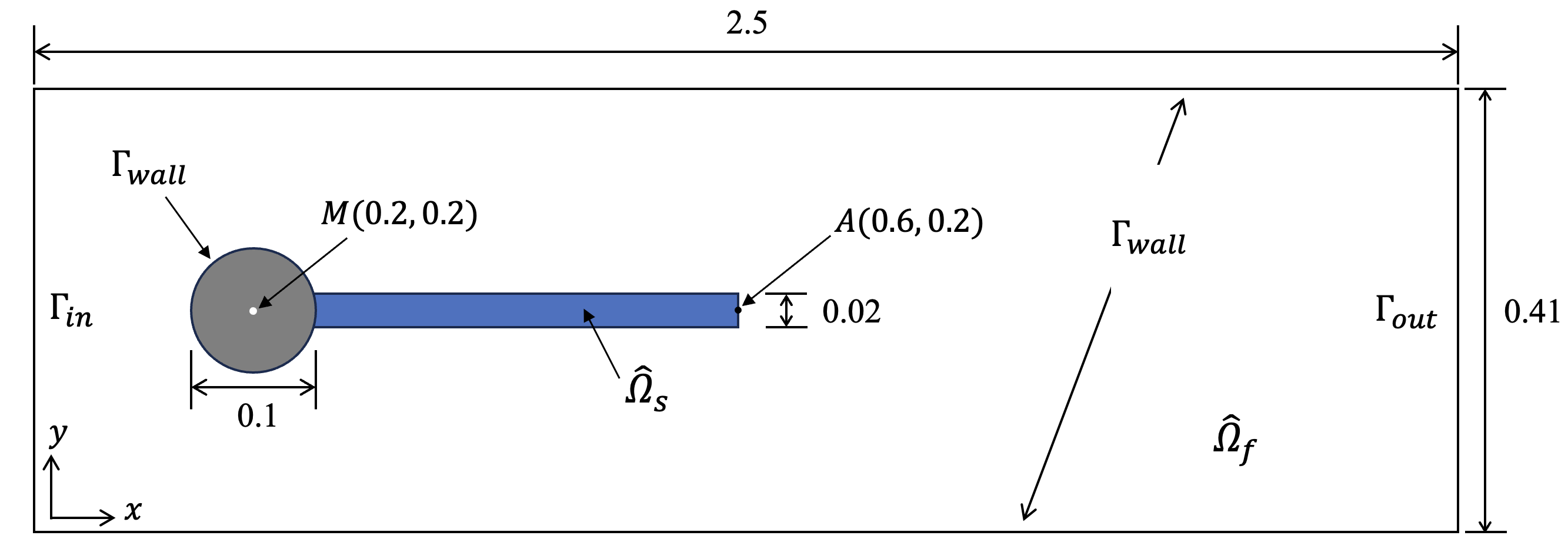}
    \caption{Geometrical setup of the Turek \& Hron's benchmark FSI problem.}
    \label{configuration}
\end{figure}

In this example, we consider the well-known Turek \& Hron's benchmark problem \cite{turek2010numerical} to validate the performance of the proposed scheme (as reported in Algorithm \ref{semi}). The benchmark configuration, as shown in Figure \ref{configuration}, consists of an elastic beam (in blue) attached to a rigid cylinder (in gray) placed inside a flow channel. The corresponding physical parameters are reported in Table \ref{tab:2dparameters}. At the channel inlet, a parabolic velocity profile is imposed, with the flow field smoothly ramped from $\boldsymbol{v}=0$ at $t=0$ s to its maximum velocity at $t=2s$:
\begin{equation*}
        \boldsymbol{v}_\text{in}(y, t)= \left\{
    \begin{aligned}
     &\frac{1-\cos(\pi t/2)}{2}\cdot\frac{12y(0.41-y)}{0.41^2} \text{m}\text{ s}^{-1}&&t<2,\\
     &\frac{12y(0.41-y)}{0.41^2} \text{m}\text{ s}^{-1}&&t\ge2.
    \end{aligned}
    \right.
\end{equation*}
\textcolor{black}{
The no-slip condition
$\boldsymbol{u}=\boldsymbol{0}$ is imposed on $\Gamma_{\mathrm{wall}}$,
while $\boldsymbol{\sigma}_f\boldsymbol{n}_f=\boldsymbol{0}$ is imposed on
$\Gamma_{\mathrm{out}}$. The elastic beam is clamped at its attachment to
the rigid cylinder, where $\boldsymbol{\eta}=\boldsymbol{0}$.
}

\begin{table}[h!]
    \centering
    \begin{tabular}{l c | l c}
        \toprule
        \textbf{Parameter} & \textbf{Value} & \textbf{Parameter} & \textbf{Value} \\ 
        \midrule
        $\rho_s$ & 1000 kg$\text{ m}^{-3}$    &  $\rho_f$ & 1000 kg$\text{ m}^{-3}$      \\ 
        $\mu_s$  & $2 \times 10^6 \text{ kg m}^{-1}\text{ s}^{-2}$ & $\mu_f$  & 1 kg $\text{m}^{-1}\text{ s}^{-1}$      \\ 
        $\lambda_s$ & $8 \times 10^6 \text{ kg m}^{-1}\text{ s}^{-2}$ & $T$      & 10 s\\
        \bottomrule
    \end{tabular}
    \caption{List of parameters used in Sec.\ref{benchmark}.}
    \label{tab:2dparameters}
\end{table}
\color{black}


We fix the Robin parameter at $L=2300$ and consider different time-step sizes $\Delta t$. 
The quantities of interest are the $x$- and $y$-displacements of the beam tip 
$A=(0.6,0.2)$. Over the time interval $I=[9,10]$, their mean values and amplitudes are defined as
\begin{equation*}
    \frac{\max_{t\in I} \boldsymbol{\eta}_i + \min_{t\in I}\boldsymbol{\eta}_i}{2}\pm\frac{\max_{t\in I} \boldsymbol{\eta}_i - \min_{t\in I}\boldsymbol{\eta}_i}{2}.
    \label{interest}
\end{equation*}

The Robin--Robin and Strong-Coupling results for different time-step sizes, together with published benchmark values~\cite{turek2010numerical, failer2020parallel}, are summarized in Table~\ref{tab:th}. Both methods agree well with the benchmark data, and the agreement improves as the time step is refined.

\begin{table}[htbp]
\centering
\caption{Results of the Turek--Hron benchmark with $L=2300$ using different time-step sizes $\Delta t = 1.0\times10^{-3}s$ and $\Delta t = 2.5\times10^{-4}s$ within the interval $I = [9, 10]$.}
\begin{tabular}{lcc}
\toprule
 & $\boldsymbol{\eta}_x \times 10^{-3}$ & $\boldsymbol{\eta}_y \times 10^{-3}$ \\
\midrule
Robin--Robin ($\Delta t = 1.0\times10^{-3}$)
& $-2.81 \pm 2.66$ & $1.57 \pm 34.60$ \\

Strong-Coupling ($\Delta t = 1.0\times10^{-3}$)
& $-2.79 \pm 2.64$ & $1.46 \pm 34.29$ \\

Robin--Robin ($\Delta t = 2.5\times10^{-4}$)
& $-2.88 \pm 2.73$ & $1.48 \pm 34.88$ \\

Strong-Coupling ($\Delta t = 2.5\times10^{-4}$)
& $-2.84 \pm 2.68$ & $1.44 \pm 34.64$ \\

Turek et al.~\cite{turek2010numerical}
& $-2.86 \pm 2.70$ & $1.45 \pm 34.93$ \\

Failer \& Richter~\cite{failer2020parallel}
& $-2.89 \pm 2.73$ & $1.47 \pm 35.03$ \\
\bottomrule
\end{tabular}
\label{tab:th}
\end{table}

\begin{figure}[h!]
    \centering
\begin{subfigure}[b]{0.75\textwidth}
        \centering
        \includegraphics[width=\textwidth]{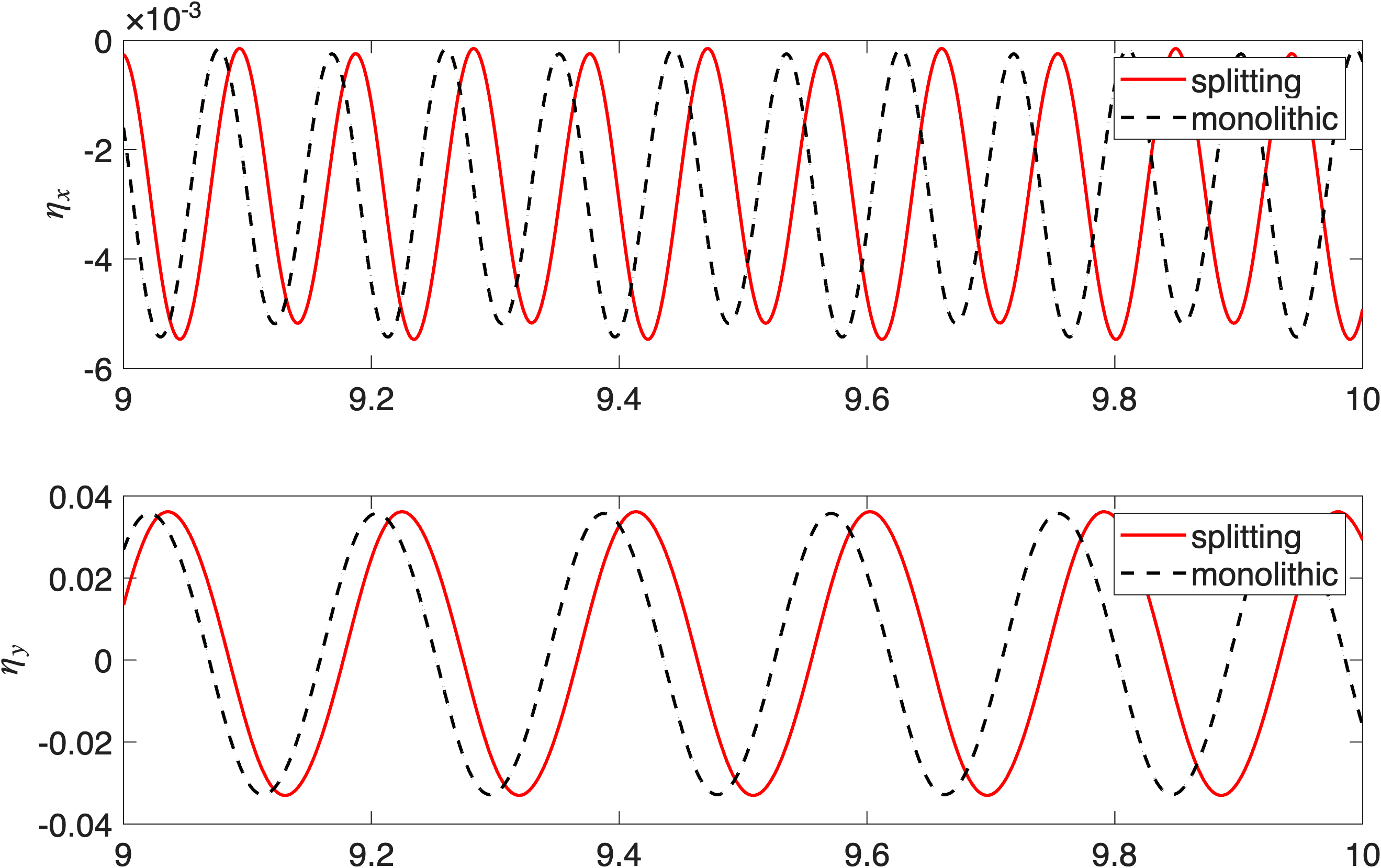}
        \caption{$\Delta t = 1\times10^{-3}$}
    \end{subfigure}
    
\begin{subfigure}[b]{0.75\textwidth}
        \centering
        \includegraphics[width=\textwidth]{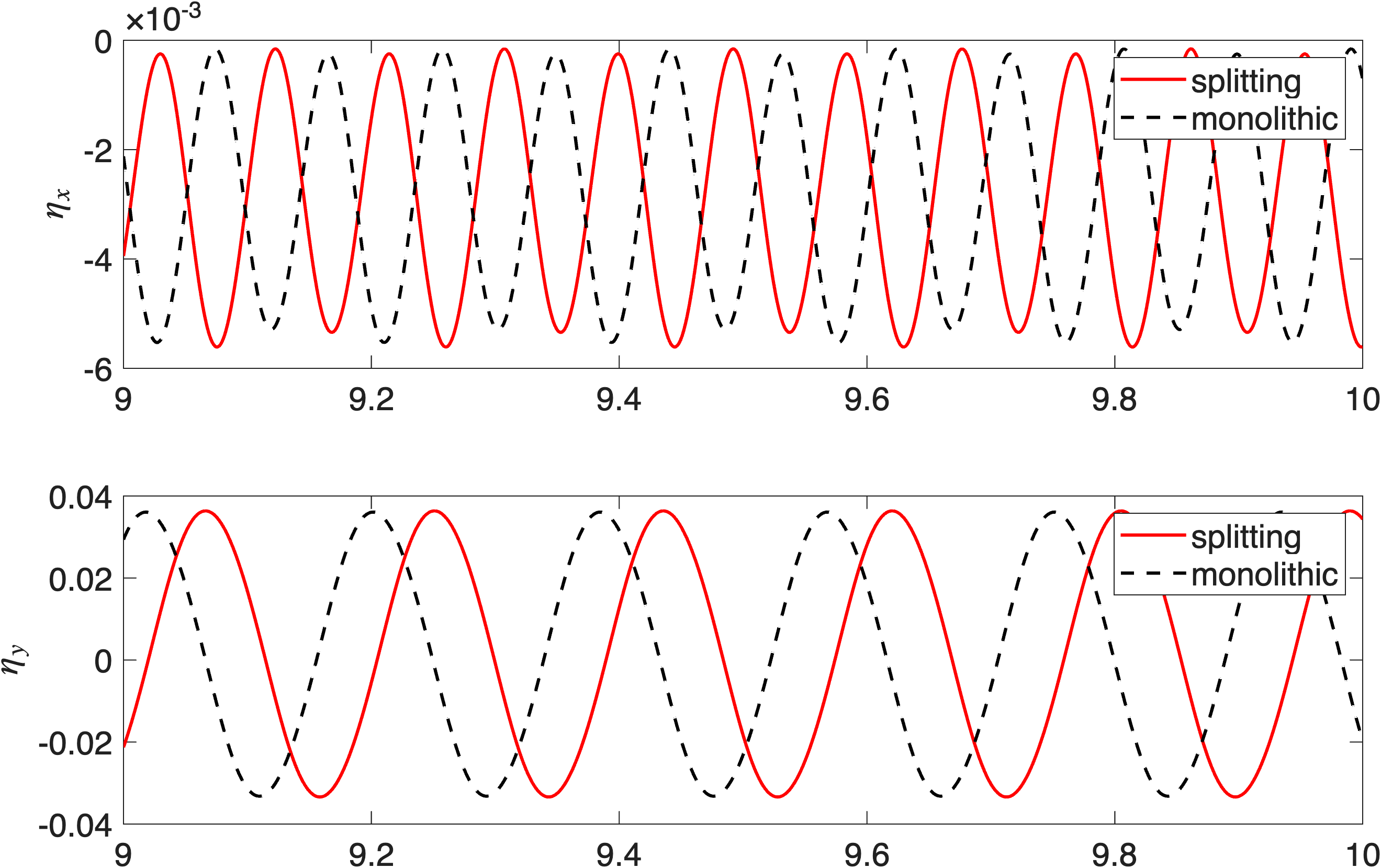}
        \caption{$\Delta t = 2.5\times10^{-4}$}
    \end{subfigure}
    \caption{The time history of horizontal (top) and vertical (bottom) deflections over the interval $I=[9, 10] ~\text{s}$, with the Robin--Robin (splitting) and Strong-Coupling (monolithic) results superimposed for comparison, are shown for (a) $\Delta t=1\times10^{-3}$ and (b) $\Delta t = 2.5\times10^{-4}$, respectively.}
    \label{fig:periodic}
\end{figure}

Figure~\ref{fig:periodic} compares the tip-displacement histories obtained
with the Robin--Robin and Strong-Coupling methods for two time-step sizes.
The two methods show close agreement in oscillation amplitude, while a small
phase difference remains in the periodic regime. This difference may be
attributed to the splitting error introduced during the initial transient,
which leads to a phase shift after the solution reaches the periodic state.

\color{black}

{\color{black}
\subsection{Pressure-wave propagation benchmark}
\label{pressure_wave}

To address the biomedical application, we consider the pressure-wave
propagation benchmark for blood flow in a compliant vessel, following \cite{bukavc2014modular}.
As illustrated in Figure~\ref{fig:pressure_wave_domain}, the reference fluid domain is $\widehat{\Omega}_f=(0,\ell)\times(0,R)$, 
while the structural domain is
$\widehat{\Omega}_s=(0,\ell)\times(R,R+h_s)$, where
$\ell=6~\mathrm{cm}$, $R=0.5~\mathrm{cm}$, and
$h_s=0.1~\mathrm{cm}$. Only the upper half of the vessel is modeled, with
the lower boundary $\Gamma_{f}^{\mathrm{sym}}$ corresponding to the axis of symmetry.

\begin{figure}[htbp]
    \centering
    \includegraphics[width=0.65\textwidth]{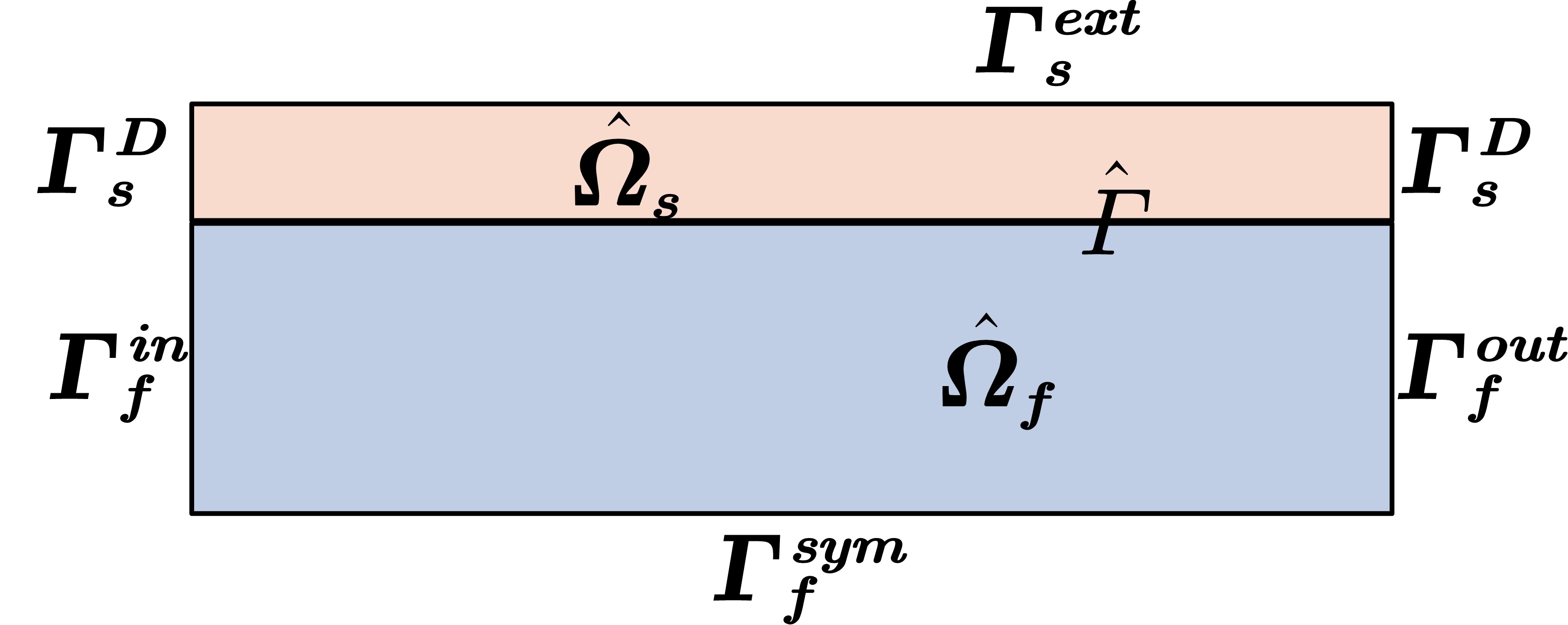}
    \caption{Reference configuration of the pressure-wave propagation
    benchmark.}
    \label{fig:pressure_wave_domain}
\end{figure}

The flow is driven by a time-dependent pressure pulse applied at the inlet
$\Gamma_{f}^{\mathrm{in}}$. We
prescribe
$\boldsymbol{\sigma}_f\boldsymbol{n}_f
=-p_{\mathrm{in}}(t)\boldsymbol{n}_f$ on $\Gamma_{f}^{\mathrm{in}}$, whereas a zero normal traction is imposed at the outlet $\Gamma_{f}^{\mathrm{out}}$. The inlet pressure is given by
\begin{equation*}
p_{\mathrm{in}}(t)=
\begin{cases}
\dfrac{p_{\max}}{2}
\left[
1-\cos\left(\dfrac{2\pi t}{t_{\max}}\right)
\right],
& 0\leq t\leq t_{\max},\\[2mm]
0, & t>t_{\max},
\end{cases}
\end{equation*}
where $p_{\max}=1.333\times10^{4}~\mathrm{dyne/cm^2}$ and
$t_{\max}=3\times10^{-3}~\mathrm{s}$. On the symmetry boundary
$\Gamma_{f}^{\mathrm{sym}}$, the conditions $\boldsymbol{u}_{f,y}=0$ and $\partial \boldsymbol{u}_{f,x}/\partial y=0$ are imposed. The structure is fixed at the inlet and outlet ends. On the exterior structural boundary $\Gamma_s^{\mathrm{ext}}$,
the axial displacement is constrained, namely, $\boldsymbol{\eta}_x=0$. All physical quantities are given in the CGS system. The fluid physical parameters are taken as $\rho_f=1.0~\mathrm{g/cm^3}$ and
$\mu_f=0.035~\mathrm{g/(cm\,s)}$. For the vessel wall, we
set $\rho_s=1.1~\mathrm{g/cm^3}$,
$\mu_s=5.75\times10^{5}~\mathrm{dyne/cm^2}$, and
$\lambda_s=1.7\times10^{6}~\mathrm{dyne/cm^2}$.
A zeroth-order recoil term $c_s\boldsymbol{\eta}$ is included in the structural momentum equation with $c_s=4.0\times10^{6}~\mathrm{dyne/cm^4}$.


The simulations are performed over $0\leq t\leq0.012~\mathrm{s}$, with the
Robin parameter fixed at $L=500$. We consider the time-step sizes
$\Delta t=1.0\times10^{-4}$, $5.0\times10^{-5}$, $2.5\times10^{-5}$, and
$1.25\times10^{-5}~\mathrm{s}$. The mesh size is $h=0.02~\mathrm{cm}$,
with 17500 triangular elements in the fluid domain and 3606 in the structure
domain. A highly resolved Strong-Coupling solution with
$\Delta t=1.0\times10^{-6}~\mathrm{s}$ and
$h=3.125\times10^{-3}~\mathrm{cm}$ is used as the reference solution.



To assess the effect of correction subiterations, Robin-sub0 denotes the
scheme without correction, while Robin-sub$k$ ($k = 1, 2$) denotes Robin-sub0 followed
by $k$ correction steps using updated interface data.
Figure~\ref{fig:pw_solution_comparison} shows that one correction substantially
improves the agreement with the reference solution, while two corrections
yield displacement and pressure profiles close to those of the Strong-Coupling
method.

\begin{figure}[H]
    \centering
    \captionsetup[subfigure]{skip=0pt}

    \begin{subfigure}[t]{0.43\textwidth}
        \centering
        \includegraphics[width=\linewidth]
        {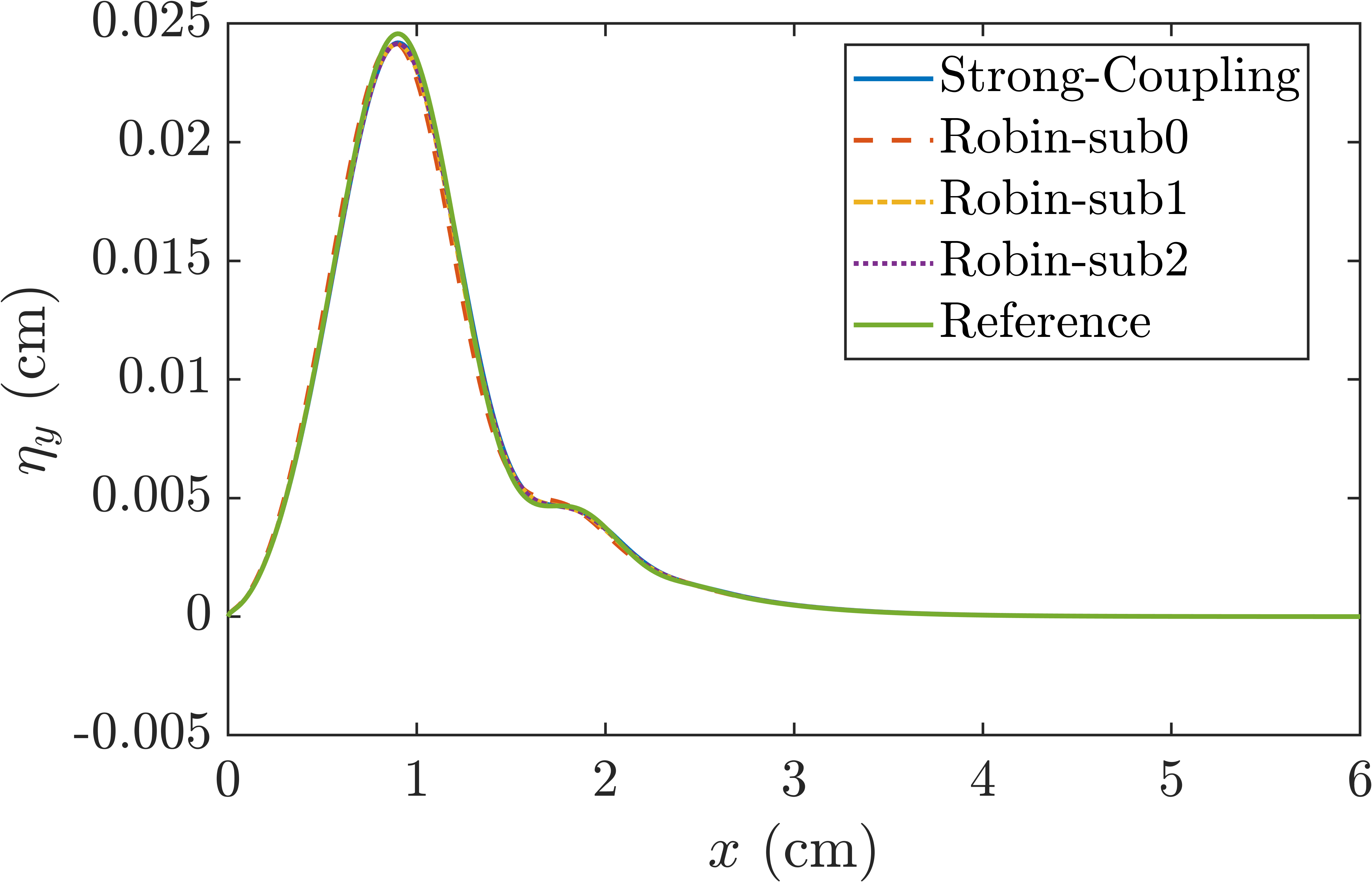}
        \caption{$y$-displacement at $t=4\,\mathrm{ms}$.}
        \label{fig:pw_disp_4ms}
    \end{subfigure}
    \hspace{0.04\textwidth}
    \begin{subfigure}[t]{0.43\textwidth}
        \centering
        \includegraphics[width=\linewidth]
        {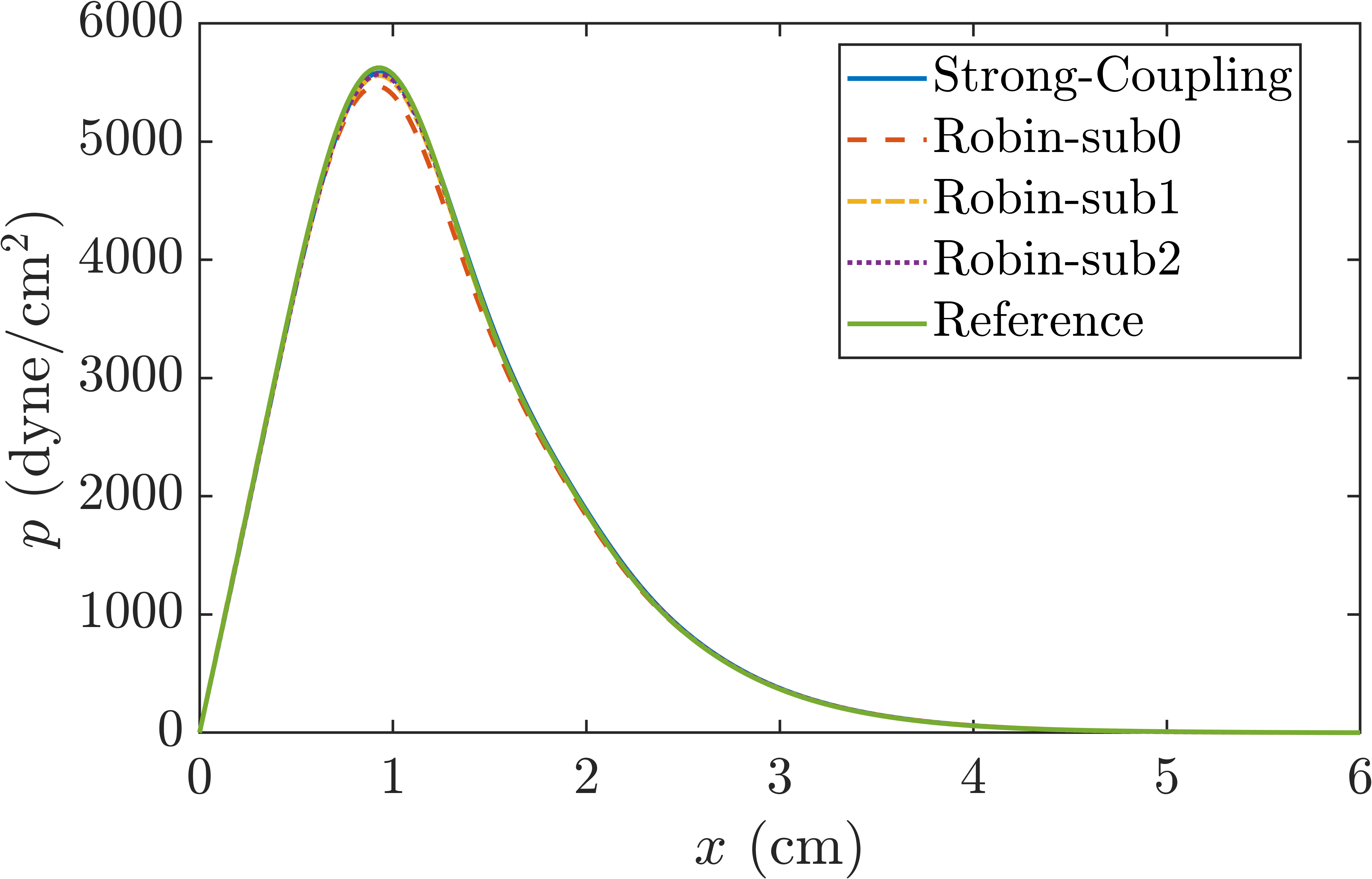}
        \caption{Pressure at $t=4\,\mathrm{ms}$.}
        \label{fig:pw_pressure_4ms}
    \end{subfigure}

    \vspace{-0.12cm}

    \begin{subfigure}[t]{0.43\textwidth}
        \centering
        \includegraphics[width=\linewidth]
        {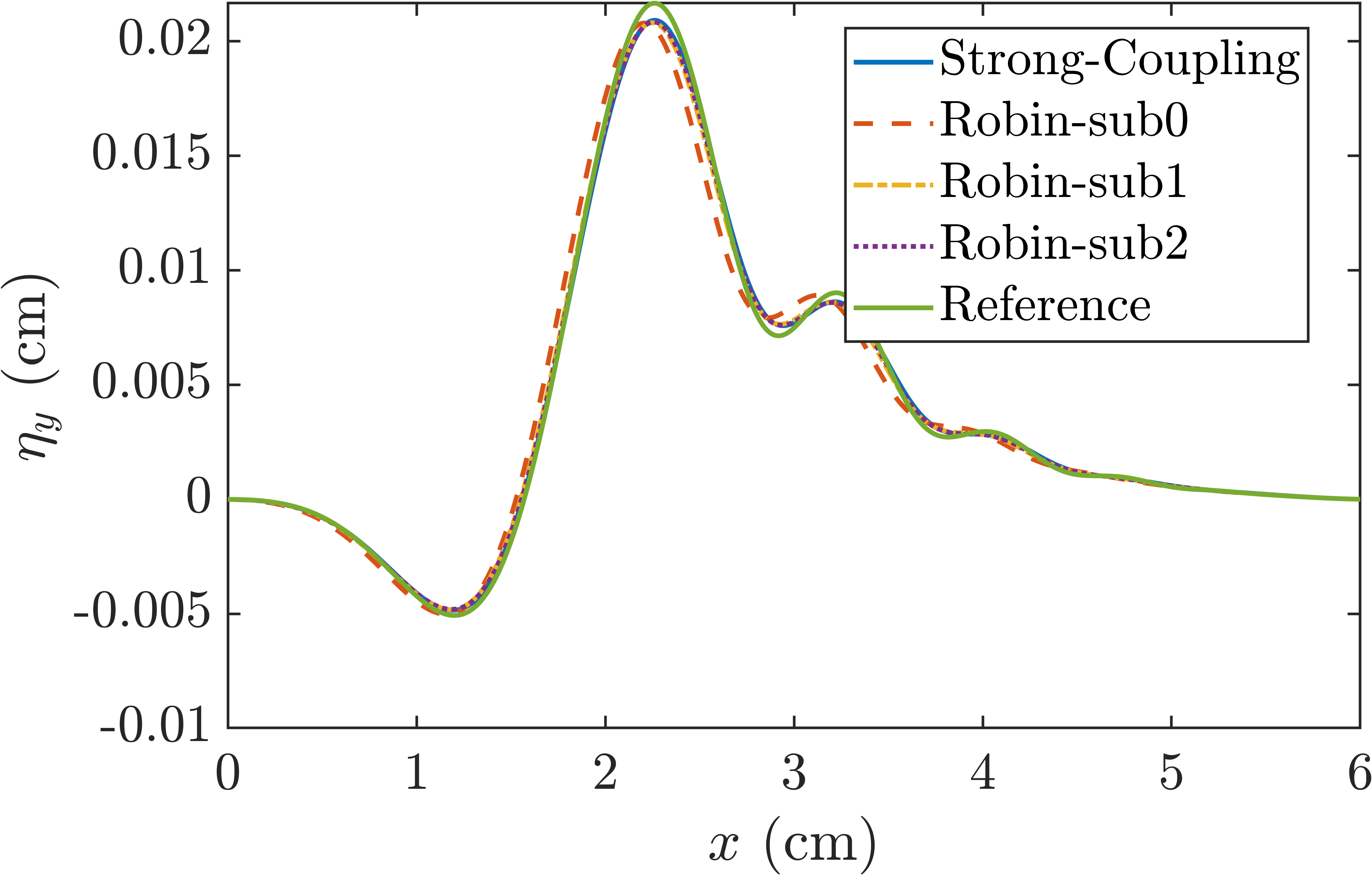}
        \caption{$y$-displacement at $t=8\,\mathrm{ms}$.}
        \label{fig:pw_disp_8ms}
    \end{subfigure}
    \hspace{0.04\textwidth}
    \begin{subfigure}[t]{0.43\textwidth}
        \centering
        \includegraphics[width=\linewidth]
        {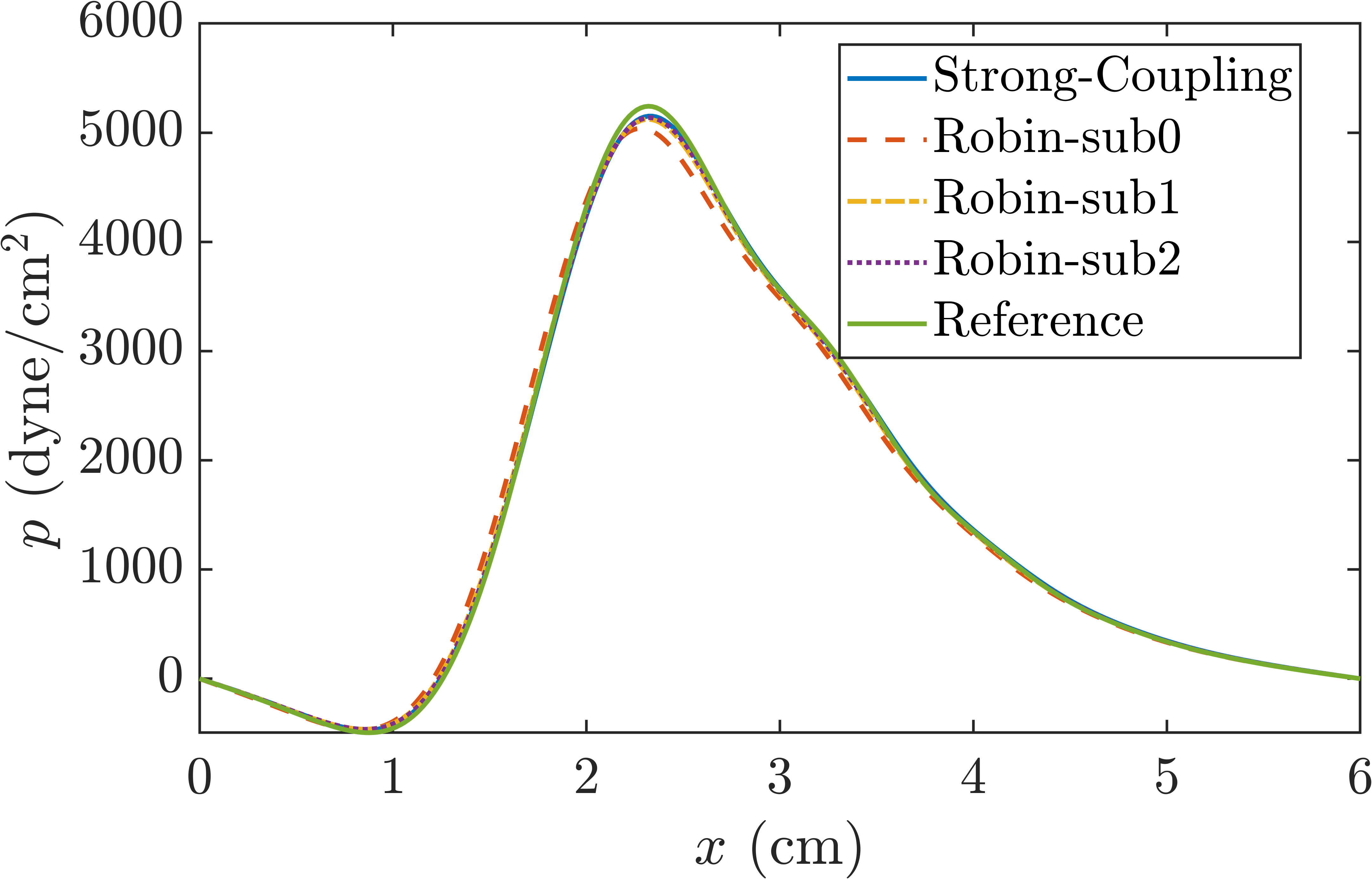}
        \caption{Pressure at $t=8\,\mathrm{ms}$.}
        \label{fig:pw_pressure_8ms}
    \end{subfigure}

    \vspace{-0.12cm}

    \begin{subfigure}[t]{0.43\textwidth}
        \centering
        \includegraphics[width=\linewidth]
        {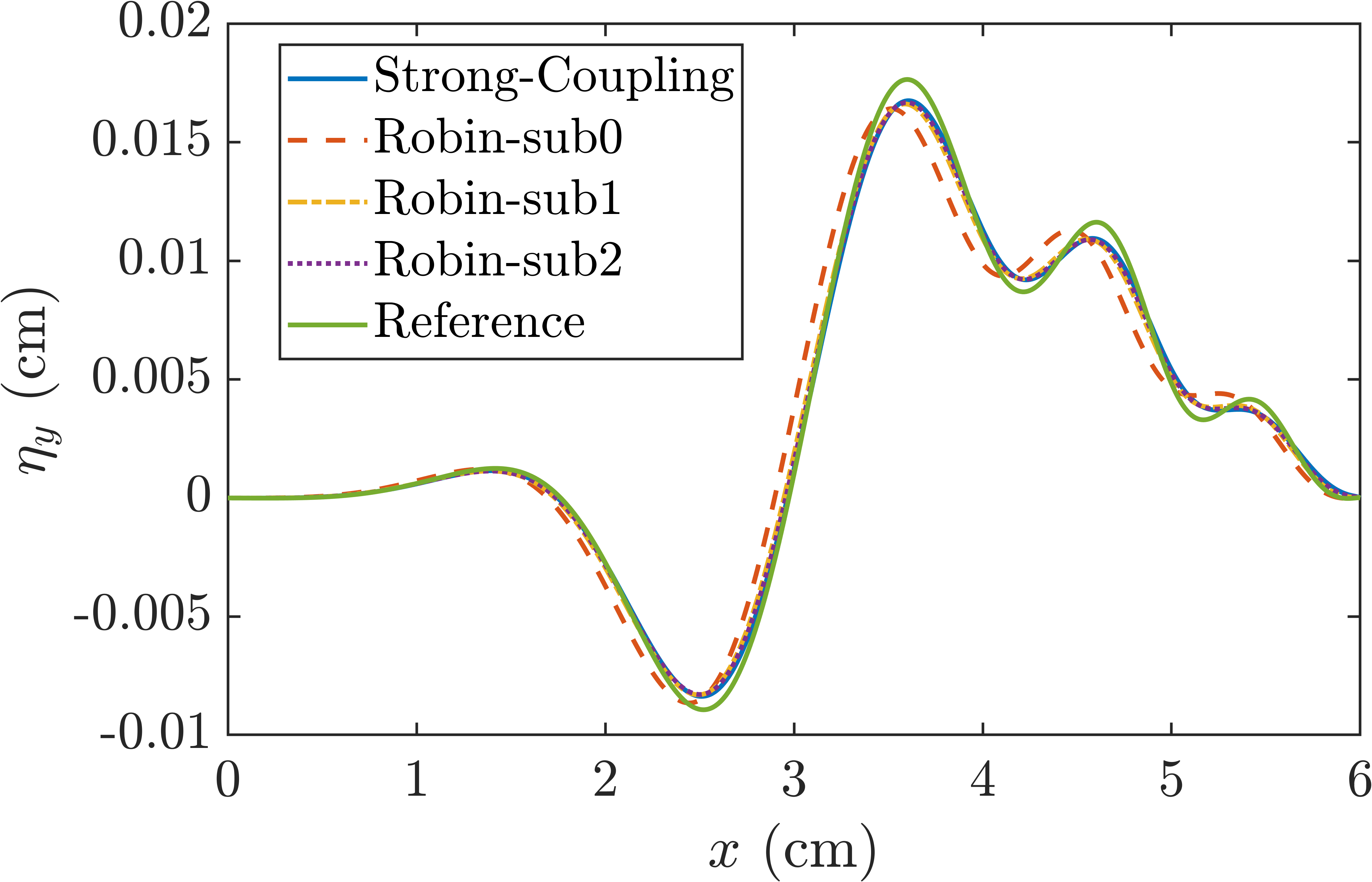}
        \caption{$y$-displacement at $t=12\,\mathrm{ms}$.}
        \label{fig:pw_disp_12ms}
    \end{subfigure}
    \hspace{0.04\textwidth}
    \begin{subfigure}[t]{0.43\textwidth}
        \centering
        \includegraphics[width=\linewidth]
        {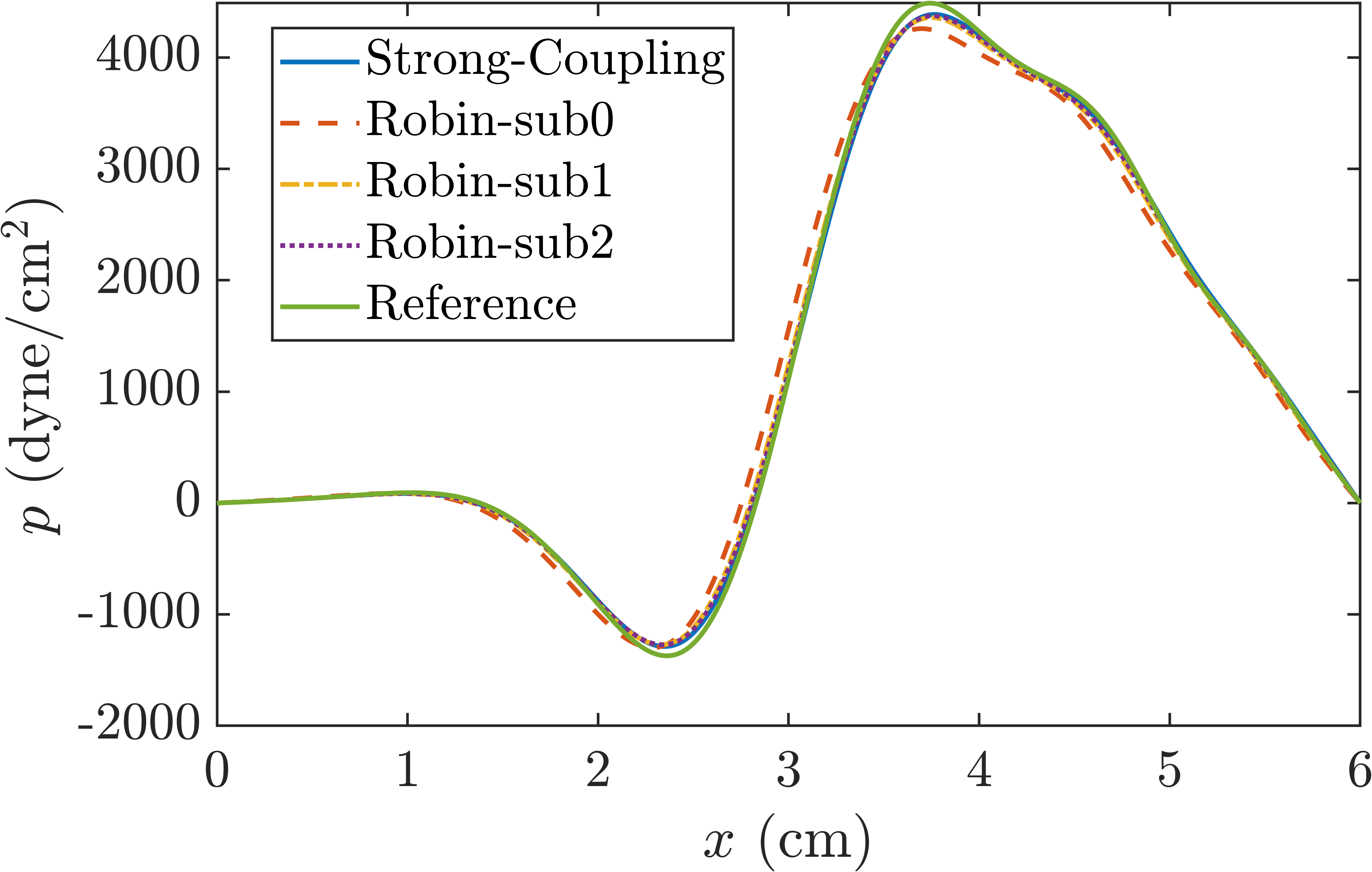}
        \caption{Pressure at $t=12\,\mathrm{ms}$.}
        \label{fig:pw_pressure_12ms}
    \end{subfigure}

    \caption{
    Comparison of the Strong-Coupling and Robin-Robin solutions with the reference solution
    for the pressure-wave benchmark at $t=4$, $8$, and $12\,\mathrm{ms}$.
    The left and right columns show the interface $y$-displacement and the
    pressure along the symmetry boundary, respectively.
    }
    \label{fig:pw_solution_comparison}
\end{figure}

Tables~\ref{tab:pw_error_disp} and~\ref{tab:pw_error_pressure} report the
relative $L^2$ errors at $t=12\,\mathrm{ms}$, defined by
$e_{\eta}=\|\eta_y-\eta_{y,\mathrm{ref}}\|_{L^2(\widehat{\Gamma})}
/\|\eta_{y,\mathrm{ref}}\|_{L^2(\widehat{\Gamma})}$ and
$e_p=\|p-p_{\mathrm{ref}}\|_{L^2(\Gamma_f^{\mathrm{sym}})}
/\|p_{\mathrm{ref}}\|_{L^2(\Gamma_f^{\mathrm{sym}})}$.
The observed rates approach first order as $\Delta t$ decreases.
The first correction substantially reduces the splitting error, while
Robin-sub2 reaches an accuracy close to the Strong-Coupling method.
Without correction, Robin-sub0 requires approximately a four-times smaller
time step to achieve comparable accuracy.

\begin{table}[htbp]
\centering
\caption{Relative $L^2$ errors and temporal convergence rates of the interface
$y$-displacement at $t=12\,\mathrm{ms}$.}
\label{tab:pw_error_disp}

\small
\renewcommand{\arraystretch}{1.15}
\setlength{\tabcolsep}{3.5pt}

\begin{tabular}{c cc cc cc cc}
\hline
& \multicolumn{2}{c}{Strong-Coupling}
& \multicolumn{2}{c}{Robin-sub0}
& \multicolumn{2}{c}{Robin-sub1}
& \multicolumn{2}{c}{Robin-sub2} \\
\cline{2-9}

$\Delta t$
& Err. (\%) & Ord.
& Err. (\%) & Ord.
& Err. (\%) & Ord.
& Err. (\%) & Ord. \\
\hline

$1.00\times10^{-4}$
& 27.20 & --
& 66.09 & --
& 35.93 & --
& 30.42 & -- \\

$5.00\times10^{-5}$
& 16.77 & 0.698
& 45.49 & 0.539
& 21.63 & 0.732
& 18.50 & 0.718 \\

$2.50\times10^{-5}$
& 9.67 & 0.795
& 28.10 & 0.695
& 12.28 & 0.817
& 10.63 & 0.799 \\

$1.25\times10^{-5}$
& 5.27 & 0.876
& 15.74 & 0.836
& 6.63 & 0.889
& 5.80 & 0.874 \\

\hline
\end{tabular}
\end{table}

\begin{table}[htbp]
\centering
\caption{Relative $L^2$ errors and temporal convergence rates of the pressure
along the symmetry boundary at $t=12\,\mathrm{ms}$.}
\label{tab:pw_error_pressure}
\small
\renewcommand{\arraystretch}{1.15}
\setlength{\tabcolsep}{3.5pt}

\begin{tabular}{c cc cc cc cc}
\hline
& \multicolumn{2}{c}{Strong-Coupling}
& \multicolumn{2}{c}{Robin-sub0}
& \multicolumn{2}{c}{Robin-sub1}
& \multicolumn{2}{c}{Robin-sub2} \\
\cline{2-9}

$\Delta t$
& Err. (\%) & Ord.
& Err. (\%) & Ord.
& Err. (\%) & Ord.
& Err. (\%) & Ord. \\
\hline

$1.00\times10^{-4}$
& 15.02 & --
& 42.28 & --
& 22.12 & --
& 17.97 & -- \\

$5.00\times10^{-5}$
& 8.28 & 0.859
& 25.90 & 0.707
& 12.00 & 0.882
& 9.77 & 0.879 \\

$2.50\times10^{-5}$
& 4.37 & 0.922
& 14.53 & 0.833
& 6.30 & 0.930
& 5.14 & 0.925 \\

$1.25\times10^{-5}$
& 2.25 & 0.961
& 7.76 & 0.906
& 3.24 & 0.959
& 2.65 & 0.960 \\

\hline
\end{tabular}
\end{table}

Table~\ref{tab:pw_wall_time} reports the wall-clock time, and
Figure~\ref{fig:pw_work_precision} shows the corresponding work--precision
comparison. All computations were performed on a dual-socket Intel Xeon Gold 6330 CPU. All runs use eight cores: eight
for the Strong-Coupling method, and a $6+2$ fluid--structure split for the
Robin--Robin method. For Robin-sub0, a smaller time step is required to reach an accuracy comparable
to the Strong-Coupling method, which reduces its cost advantage. In contrast,
the correction subiterations improve the work--precision balance. In particular,
Robin-sub2 achieves nearly the same accuracy as the Strong-Coupling method while
requiring only about $71$--$73\%$ of its wall-clock time. Furthermore, we noticed that when upholding a comparable level of accuracy, Robin-sub1 outperforms Robin-sub2 in terms of computational efficiency, as shown in Figure.~\ref{fig:pw_work_precision}. Therefore, in practice, Robin-sub1 would be a more efficient choice while maintaining a comparable level of accuracy.

\begin{table}[H]
\centering
\caption{Wall-clock time (s) for the Strong-Coupling and Robin-Robin methods with different
time-step sizes. The values in parentheses denote the percentage of the wall-clock time
relative to the Strong-Coupling method at the same time-step size.}
\label{tab:pw_wall_time}
\small
\renewcommand{\arraystretch}{1.15}
\setlength{\tabcolsep}{3.5pt}
\begin{tabular}{c|cccc}
\hline
$\Delta t$
& Strong-Coupling
& Robin-sub0
& Robin-sub1
& Robin-sub2
\\
\hline
$1.0\times10^{-4}$
& 417.52
& 106.44 (25.5\%)
& 204.53 (49.0\%)
& 304.51 (72.9\%)
\\

$5.0\times10^{-5}$
& 830.00
& 200.25 (24.1\%)
& 400.28 (48.2\%)
& 595.70 (71.8\%)
\\

$2.5\times10^{-5}$
& 1666.24
& 399.03 (23.9\%)
& 798.15 (47.9\%)
& 1195.25 (71.7\%)
\\

$1.25\times10^{-5}$
& 3308.15
& 796.64 (24.1\%)
& 1574.48 (47.6\%)
& 2357.94 (71.3\%)
\\
\hline
\end{tabular}
\end{table}

\begin{figure}[H]
    \centering
    \includegraphics[width=0.95\textwidth]{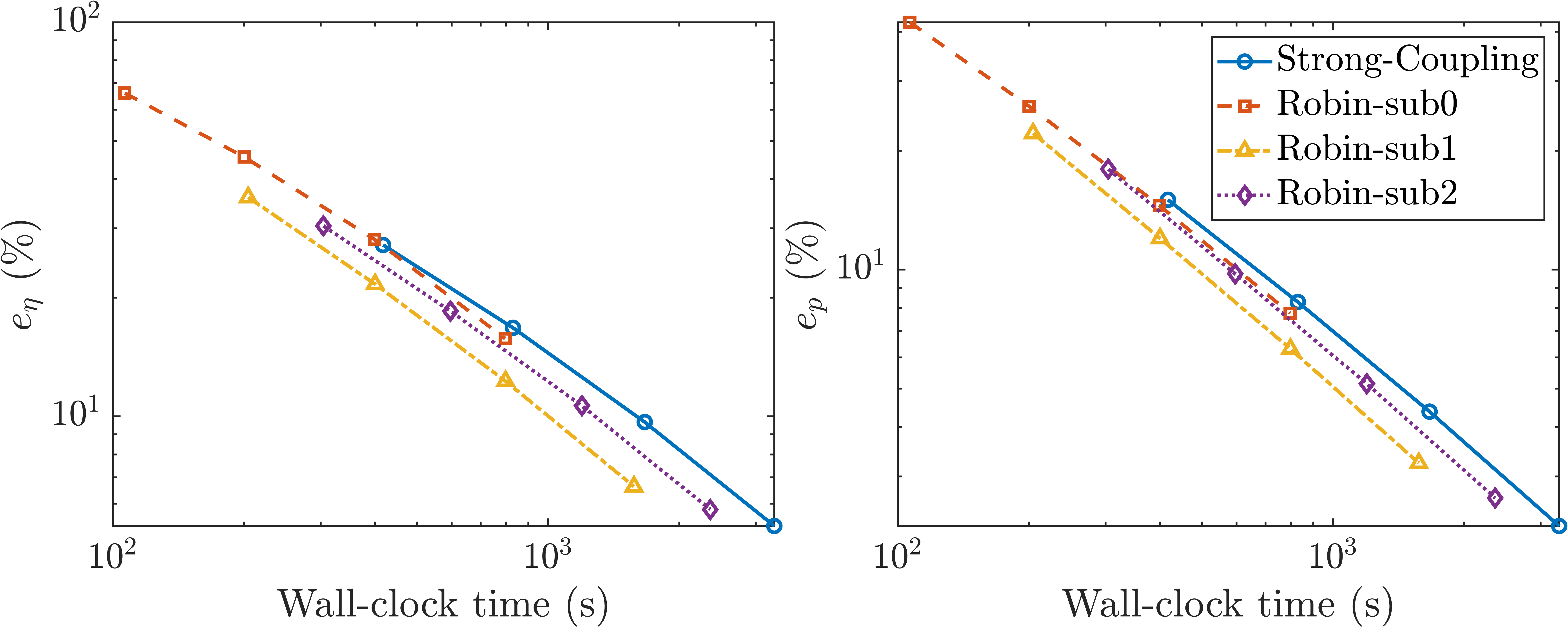}
    \caption{Work--precision comparison for the pressure-wave benchmark.}
    \label{fig:pw_work_precision}
\end{figure}

}

{\color{black}
\subsection{3D benchmark test}
\label{3dnumerical}

We further consider the three-dimensional FSI benchmark proposed in
\cite{failer2020parallel}, which can be viewed as a three-dimensional
extension of the Turek--Hron benchmark in Section~\ref{benchmark}.
The configuration is shown in Figure~\ref{fig:3dbenchmark}. The geometry
and boundary conditions are the natural three-dimensional counterparts
of those considered in Section~\ref{benchmark}. At the inlet, a bi-parabolic velocity profile is prescribed and smoothly
ramped from rest during the first second:
\begin{equation*}
\boldsymbol{u}_{\mathrm{in}}(y,z,t)
=
\begin{cases}
\displaystyle
\left(
\frac{1}{2}-\frac{1}{2}\cos(\pi t)
\right)
1.75
\frac{36y(0.41-y)z(0.41-z)}{0.41^4}\,
\boldsymbol{e}_x,
& 0\leq t<1,
\\[2mm]
\displaystyle
1.75
\frac{36y(0.41-y)z(0.41-z)}{0.41^4}\,
\boldsymbol{e}_x,
& t\geq1.
\end{cases}
\end{equation*}
A stress-free condition is imposed at the outlet $x=2.8$, while no-slip
conditions are imposed on the remaining rigid walls. The physical parameters are the same as those used in
Section~\ref{benchmark}, yielding $\mathrm{Re}=175$. We monitor the $x$-, $y$-, and
$z$-components of the displacement at the beam tip
$B=(0.9,0.2,0.3)$. The mean values and amplitudes are evaluated over
the periodic interval $I=[9,10]$. Unless otherwise stated, $L=4500$.
The time step is $\Delta t=10^{-3}\,\mathrm{s}$ on $[0,2]$ and
$\Delta t=5\times10^{-4}\,\mathrm{s}$ thereafter.
\begin{figure}[htbp]
    \centering
    \includegraphics[width=0.85\linewidth]
    {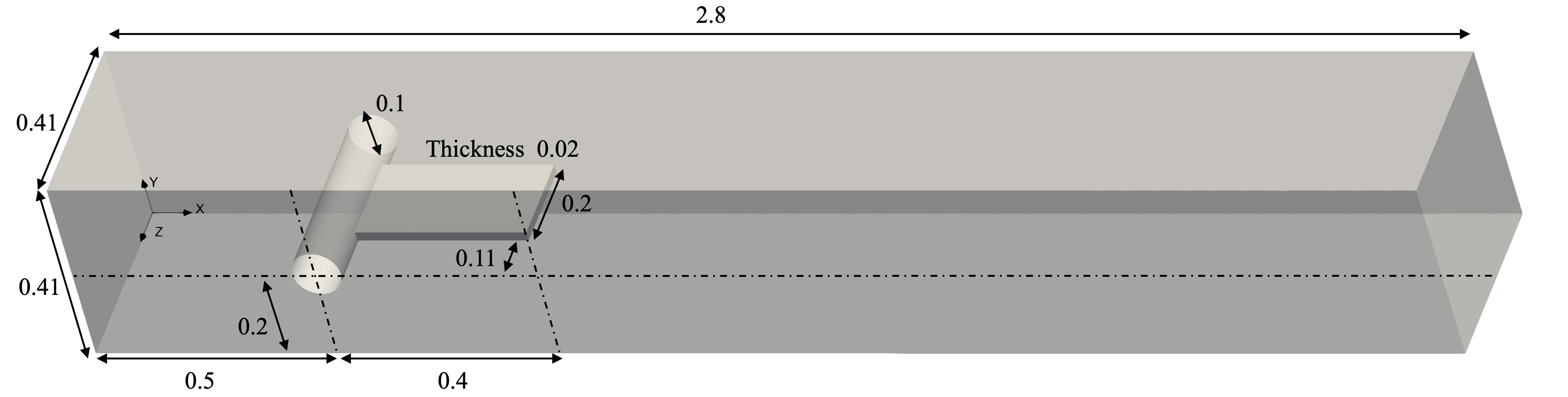}
    \caption{Geometrical configuration of the three-dimensional FSI benchmark.}
    \label{fig:3dbenchmark}
\end{figure}

We first compare the concurrent and sequential Robin--Robin
implementations and examine the effect of workload imbalance. We then
validate the three-dimensional response and compare the memory requirements.

\subsubsection{Concurrent implementation and load balance}
\label{sec:3d_concurrent}

To isolate the benefit of concurrency, we compare concurrent and sequential
implementations of the same Robin--Robin scheme. Both use the same
discretization, time step, Robin parameter, solver, and total number of
MPI processes. In the sequential implementation, the fluid and structure subproblems are
solved successively using all available MPI processes. In the concurrent
implementation, the processes are divided between the two subproblems,
The benefit of concurrency depends on the workload balance between
the two subproblems.


To investigate this effect, three mesh configurations are considered.
The fluid mesh is fixed while the structure mesh is progressively reduced.
All computations use a total of $32$ MPI processes. For each configuration,
we report the end-to-end wall-clock time over $100$ time steps, including
initialization, assembly, subproblem solution, and interface communication.
These contributions are not timed separately because they may overlap in the
concurrent implementation. The results are shown in
Table~\ref{tab:concurrent_sequential}. For the nearly balanced case, the
concurrent implementation reduces the wall-clock time from $4486$~s to
$2555$~s, corresponding to a reduction of approximately $43\%$. As the
workload becomes increasingly fluid-dominated, the reduction decreases to
approximately $22\%$ and $18\%$ for $N_f/N_s=4.0$ and $10.0$, respectively.

\begin{table}[htbp]
\centering
\caption{
Wall-clock time for the concurrent and sequential Robin--Robin
implementations using $32$ MPI processes.
}
\label{tab:concurrent_sequential}
\small
\begin{tabular}{c|c c c c|c c c}
\hline
Case
& $N_f$
& $N_s$
& $N_f/N_s$
& $(P_f,P_s)$
& Con. (s)
& Seq. (s)
& Reduction
\\
\hline
1
& \multirow{3}{*}{169,740}
& 125,624
& 1.4
& $(16,16)$
& 2555
& 4486
& 43.0\%
\\
2
&
& 42,618
& 4.0
& $(28,4)$
& 1484
& 1908
& 22.2\%
\\
3
&
& 16,914
& 10.0
& $(31,1)$
& 1346
& 1636
& 17.7\%
\\
\hline
\end{tabular}
\end{table}

These results confirm that the benefit of concurrent execution is largest
for a balanced workload and decreases as the imbalance increases.
Nevertheless, a reduction in wall-clock time is observed in all three
cases. The actual gain also depends on the strong-scaling behavior of
the individual subsolvers and the communication overhead.

\subsubsection{Benchmark validation and computational resources}
\label{sec:3d_validation}

We next assess the three-dimensional response of the Robin--Robin method.
The results are compared with the two monolithic computations reported in
\cite{failer2020parallel}. Figure~\ref{fig:3d_tip_displacement} shows the
tip-displacement histories over the periodic interval $I=[9,10]$. The corresponding mean values and amplitudes are compared with the
published results with a good agreement, as shown in Table~\ref{tab:3d_benchmark_validation}.

\begin{figure}[h!]
    \centering

    \begin{subfigure}[b]{0.45\textwidth}
        \centering
        \includegraphics[width=\linewidth]
        {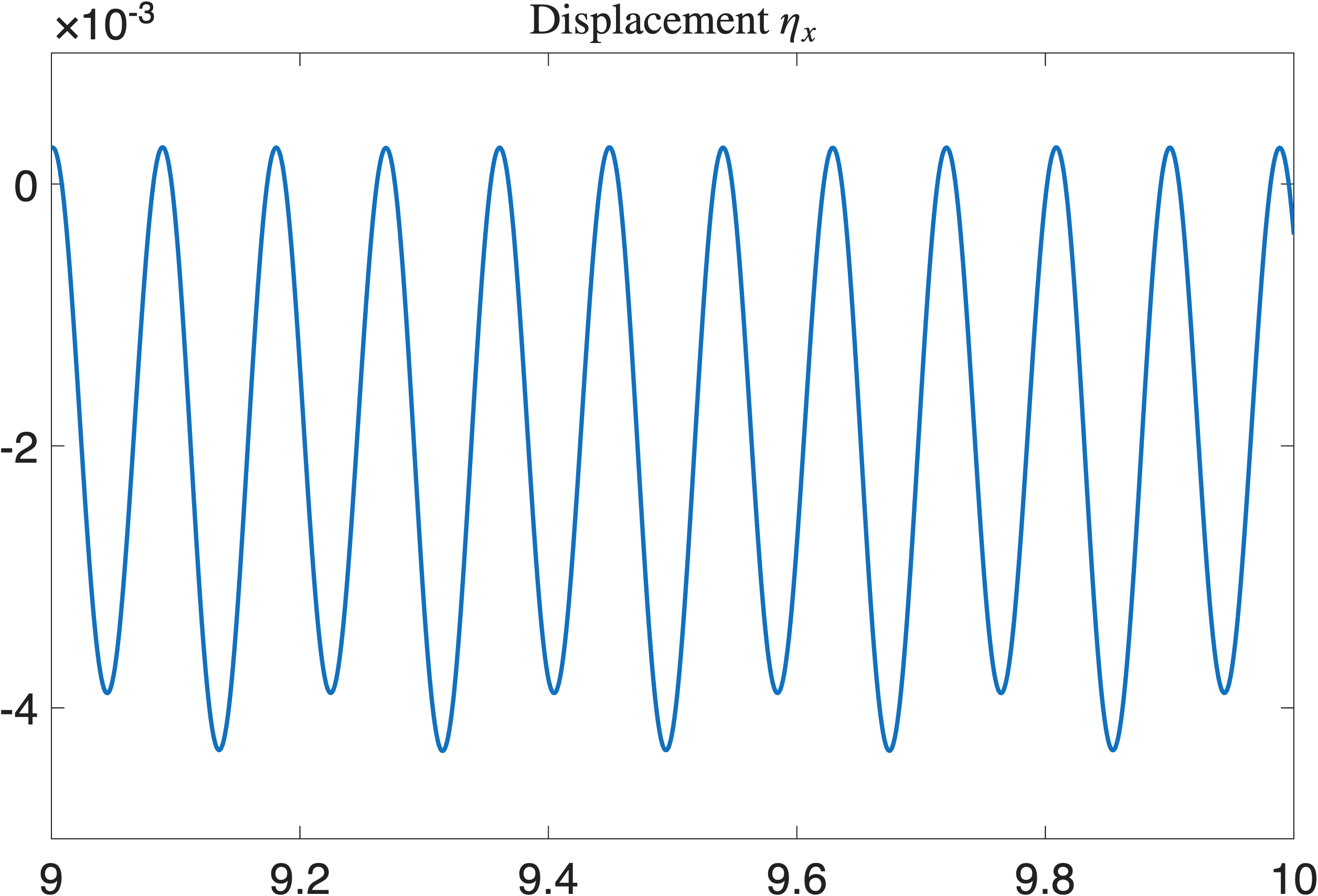}
    \end{subfigure}
    \hfill
    \begin{subfigure}[b]{0.45\textwidth}
        \centering
        \includegraphics[width=\linewidth]
        {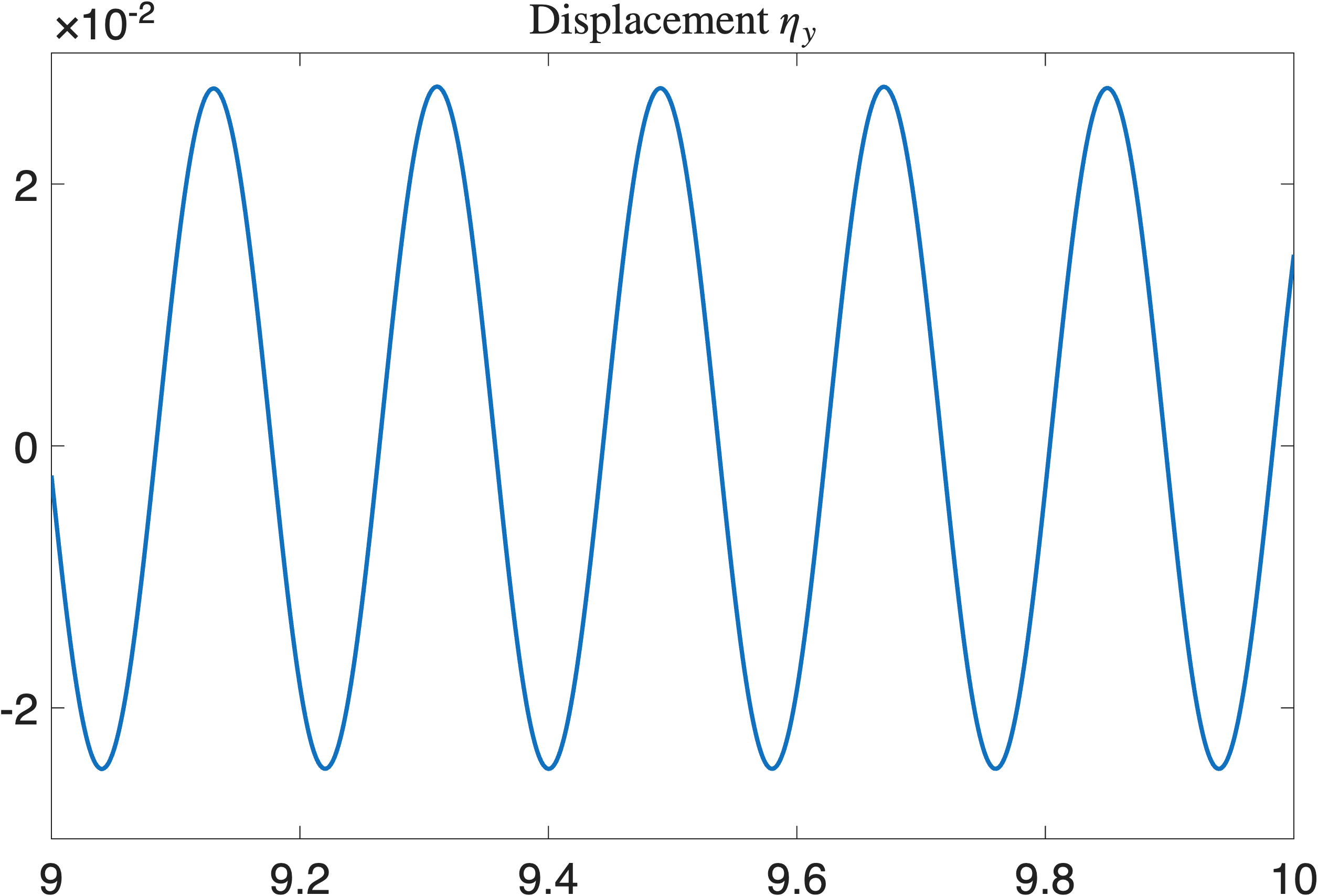}
    \end{subfigure}

    \vspace{0.08cm}

    \begin{subfigure}[b]{0.45\textwidth}
        \centering
        \includegraphics[width=\linewidth]
        {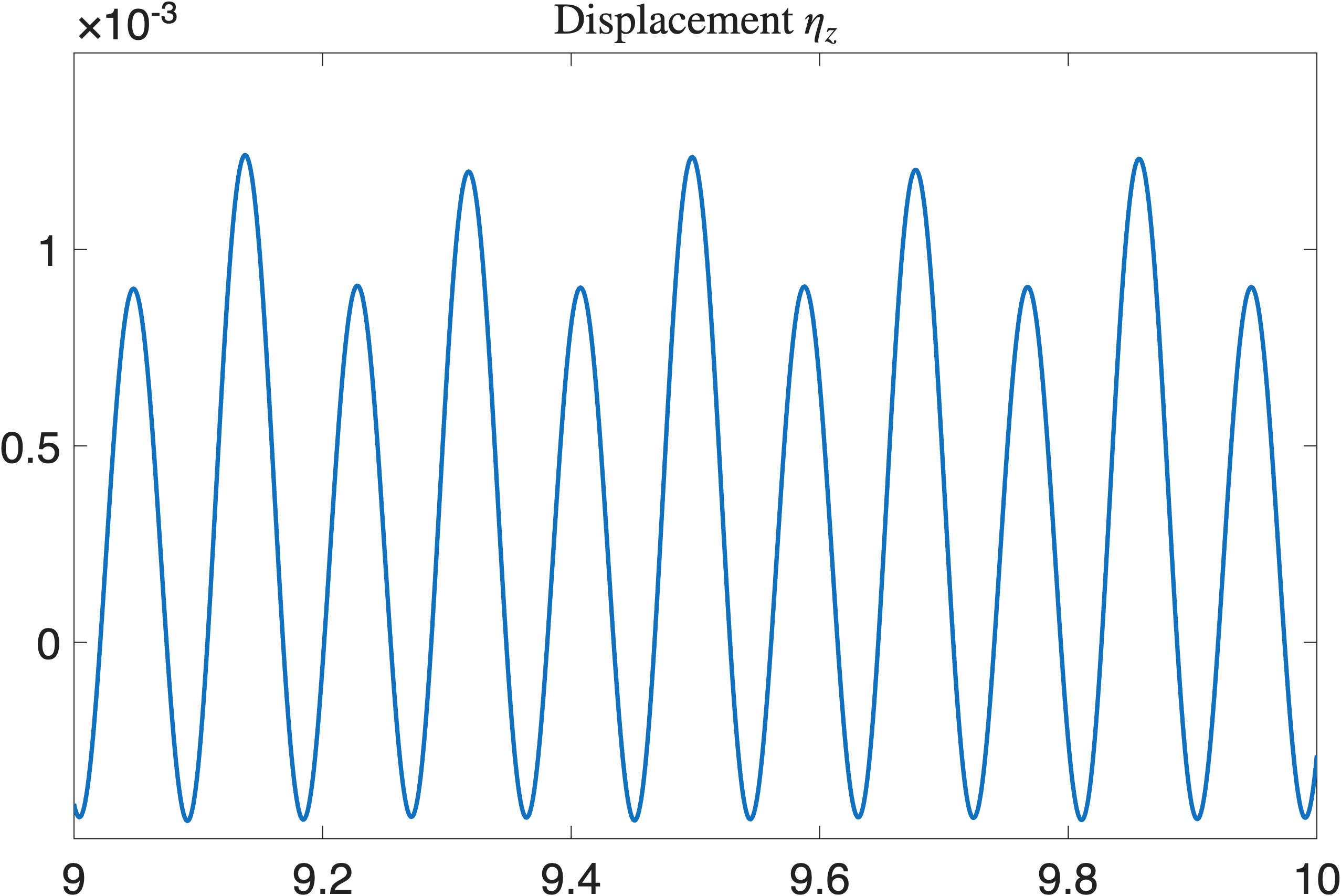}
    \end{subfigure}
    \hfill
    \begin{subfigure}[b]{0.45\textwidth}
        \centering
        \includegraphics[width=\linewidth]
        {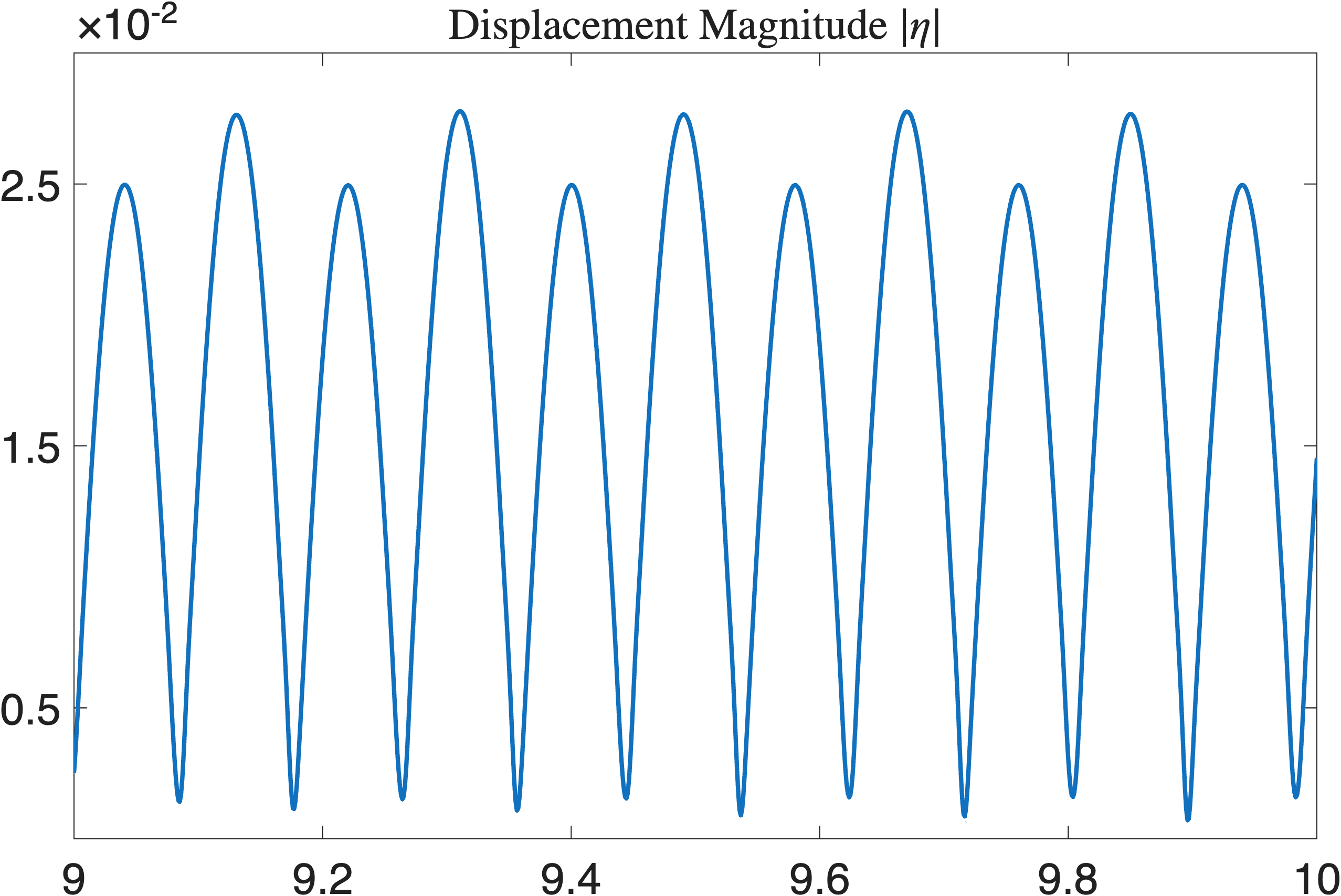}
    \end{subfigure}

    \caption{
    Tip displacement of the elastic beam over $I=[9,10]$.
    }
    \label{fig:3d_tip_displacement}
\end{figure}

\begin{table}[H]
\centering
\caption{
Comparison with the 3D FSI benchmark results in
\cite{failer2020parallel} over $I=[9,10]$.
}
\label{tab:3d_benchmark_validation}
\small
\renewcommand{\arraystretch}{1.15}
\setlength{\tabcolsep}{4pt}
\begin{tabular}{lcccc}
\toprule
&
$\boldsymbol{\eta}_x\times10^{-3}$
&
$\boldsymbol{\eta}_y\times10^{-3}$
&
$\boldsymbol{\eta}_z\times10^{-3}$
&
DoFs
\\
\midrule
Reference 1
& $-2.650 \pm 2.881$
& $1.566 \pm 29.224$
& $0.435 \pm 1.091$
& 463,988
\\
Robin--Robin
& $-2.024 \pm 2.305$
& $1.386 \pm 26.049$
& $0.393 \pm 0.847$
& 1,580,632
\\
Reference 2
& $-2.143 \pm 2.383$
& $2.699 \pm 25.594$
& $0.486 \pm 0.877$
& 3,531,304
\\
\bottomrule
\end{tabular}
\end{table}

The Robin--Robin results reproduce the characteristic periodic response
and are generally within the range of the two published monolithic
computations. The differences between the two reference results also
indicate the sensitivity of this benchmark to spatial discretization.
To establish reliable reference values and construct a comprehensive dataset, future work will involve performing these computations on systematically refined meshes.

We further examine the computational cost of the Robin--Robin
implementation. Figure~\ref{fig:3d_computational_cost}(a) shows its
strong-scaling behavior as the number of MPI processes increases from
$1$ to $32$. The reported wall-clock time is the mean computational
time per time step averaged over $100$ time steps. The wall-clock time
decreases consistently with increasing core count. The $32$-core
monolithic result is included as a reference. The peak memory consumption is compared in Figure~\ref{fig:3d_computational_cost}(b). The Robin--Robin and monolithic computations have comparable problem sizes of approximately
$2.7$ million degrees of freedom and are executed on the same
AMD EPYC 7662 system. At $32$ MPI processes, the Robin--Robin implementation requires less
wall-clock time than the monolithic computation.
The measured speedup from 1 to 32 cores is approximately $14$ times.
Moreover, the monolithic
computation requires approximately $656$ GB of peak memory, about four
times that of the Robin--Robin implementation. The reduced memory
requirement results from solving separate fluid and structure subsystems
instead of a single large monolithic system.

\begin{figure}[htbp]
    \centering

    \begin{subfigure}[b]{0.47\textwidth}
        \centering
        \includegraphics[width=\linewidth]
        {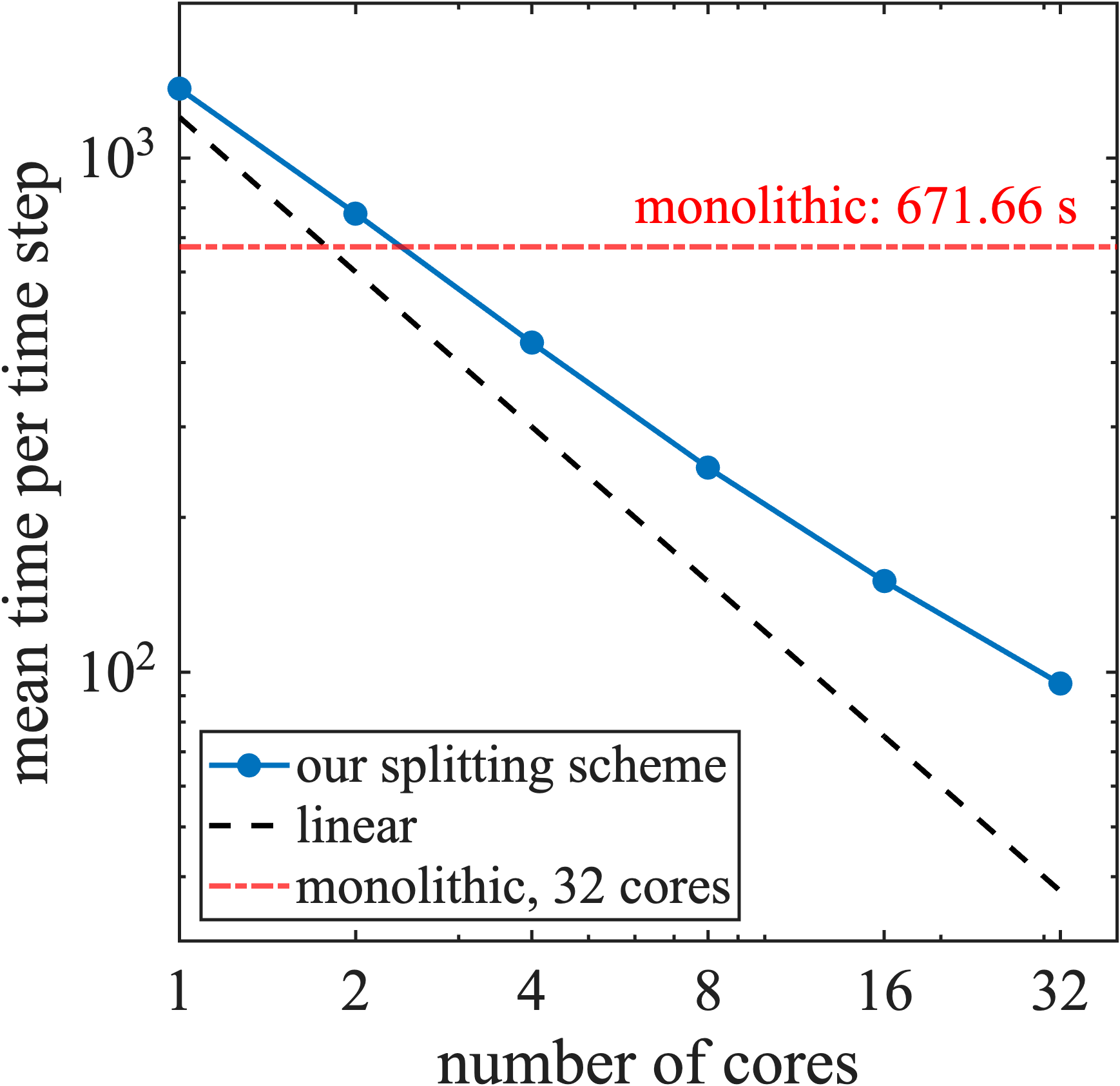}
        \caption{Wall-clock time.}
        \label{fig:3d_parallel_time}
    \end{subfigure}
    \hfill
    \begin{subfigure}[b]{0.47\textwidth}
        \centering
        \includegraphics[width=\linewidth]
        {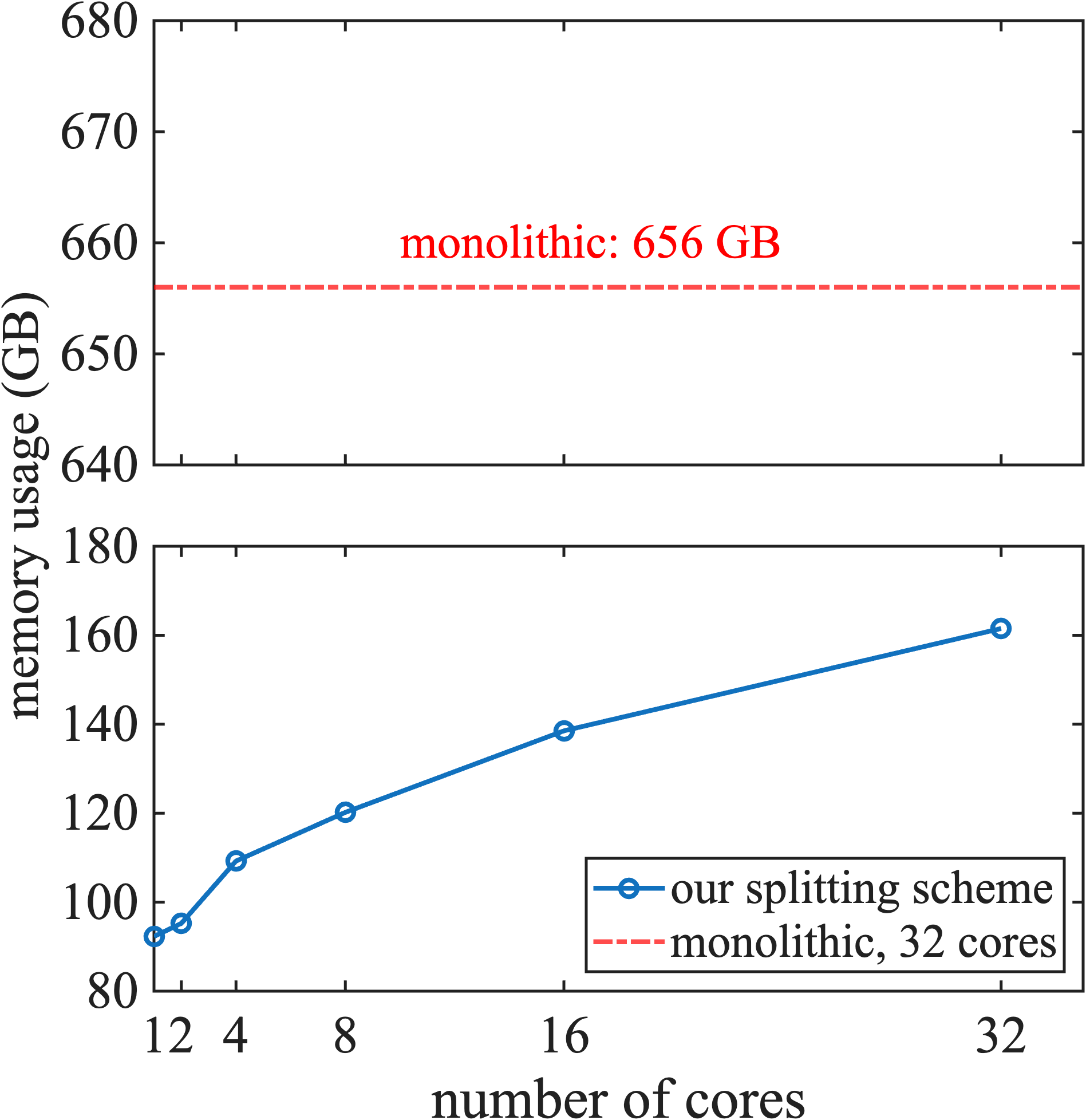}
        \caption{Peak memory consumption.}
        \label{fig:3d_memory}
    \end{subfigure}

    \caption{
    Computational cost of the Robin--Robin and monolithic implementations
    for the three-dimensional benchmark.
    }
    \label{fig:3d_computational_cost}
\end{figure}


}

\section{Conclusion}
\label{section6}
We have presented a fully parallelizable, partitioned Robin–Robin scheme for fluid–structure interaction (FSI) problems governed by the incompressible Navier–Stokes equations coupled with a linearly elastic structure. The method achieves unconditional stability while decoupling the fluid and structure subproblems in a manner that allows their simultaneous solution, rather than relying on iterative procedures. 
This property \textcolor{black}{enables concurrent execution of the two subproblems}
and \textcolor{black}{can reduce the computational cost},
particularly for large-scale 3D simulations.
Theoretical analysis established half-order temporal consistency error and unconditional stability, while numerical experiments verified the accuracy, efficiency, and robustness of the approach across both benchmark problems and realistic test cases. Beyond its demonstrated performance in classical FSI benchmarks, the proposed method has the potential to extend explicit partitioned strategies to a broader class of multiphysics systems. In particular, its flexibility suggests applicability to coupled parabolic and parabolic hyperbolic problems. 

\textcolor{black}{
It is worth noting that our numerical studies indicate that different values
of the parameter $L$ may be required for different problem configurations to
maintain accurate numerical solutions. Both the splitting-error constant and
a suitable choice of $L$ may also depend on the physical parameters of the
problem. This dependence motivates further investigation into a systematic
parameter-selection strategy. In particular, developing a theoretical
criterion for choosing $L$ represents an important direction for future work
and may further enhance the accuracy, efficiency, and general applicability
of the proposed method.
}

Taken together, the theoretical and numerical results demonstrate the potential of the proposed Robin-Robin scheme as a scalable, reliable, and computationally efficient framework for large-scale FSI and multiphysics simulations.

\section{Acknowledgment}
This work has been supported in part by the following grants or awards: NSF award (DMS-2247001), NSF of China (No.12471406, No.12371388), EPSRC Grant (EP/Y024974/1), Science and Technology Commission of Shanghai Municipality (No.22DZ2229014), and ERDF within the International Max Planck Research School for Systems and Process Engineering for a Sustainable Chemical Production (IMPRS SysProSus). The authors would like to thank Prof. Thomas Richter for helpful discussions and valuable suggestions. \textcolor{black}{The authors would also like to thank the anonymous reviewers for prompting the discussion presented in Section \ref{concurrency}.}

\color{black}
\appendix
\section{Proof of Theorem\ref{cstability}}
\begin{proof}
Test the fluid equations with $(\boldsymbol{u}^{n+1}, p^{n+1})$ and the structure equation with $\boldsymbol{\xi}^{n+1}=\partial_t\boldsymbol{\eta}^{n+1}$. After integration by parts, the pressure contribution vanishes since $\nabla\cdot\boldsymbol{u}^{n+1}=0$, while
\begin{equation}
(\boldsymbol{\sigma}_s(\boldsymbol{\eta}^{n+1}), \nabla\boldsymbol{\xi}^{n+1})_s=\frac{1}{2}\frac{d}{dt}\Vert\boldsymbol{\eta}^{n+1}\Vert_E^2.
\end{equation}
For notational simplicity, let $\mathcal{F}^n=\boldsymbol{\sigma}_f^n\boldsymbol{n}_f$ and $\mathcal{S}^n=\boldsymbol{\sigma}_s^n\boldsymbol{n}_f$. Thanks to $\boldsymbol{n}_s=-\boldsymbol{n}_f$, the tested equations give
\begin{equation}
\begin{aligned}
&\frac{d}{dt}\mathbb E^{n+1}(t)+2\mu_f\Vert\mathbb{D}(\boldsymbol{u}^{n+1})\Vert^2_f+L\left<\boldsymbol{u}^{n+1}-\frac{1}{2}\tilde{\boldsymbol{\xi}}^n-\frac{1}{2}\tilde{\boldsymbol{u}}^n,\boldsymbol{u}^{n+1}\right>_\Gamma\\
&+L\left<\boldsymbol{\xi}^{n+1}-\frac{1}{2}\tilde{\boldsymbol{\xi}}^n-\frac{1}{2}\tilde{\boldsymbol{u}}^n,\boldsymbol{\xi}^{n+1}\right>_\Gamma=\frac{1}{2}\left<\tilde{\mathcal{F}}^n+\tilde{\mathcal{S}}^n, \boldsymbol{u}^{n+1}-\boldsymbol{\xi}^{n+1}\right>_\Gamma.
\end{aligned}
\label{app:main}
\end{equation}
Notice that two Robin conditions can be written as
\begin{equation}
  \boldsymbol{u}^{n+1}-\boldsymbol{\xi}^{n+1} = \frac{1}{L}\left(\tilde{\mathcal{F}}^n+\tilde{\mathcal{S}}^n-\mathcal{F}^{n+1}-\mathcal{S}^{n+1}\right).
  \label{appen:important}
\end{equation}
Thus, by applying the polarization identity \eqref{polar} in $L^2(\Gamma)$ to \eqref{app:main}, and using \eqref{appen:important}, the two squared Robin residuals are canceled exactly. We obtain the temporal pointwise estimates
\begin{equation}
\begin{aligned}
    &\frac{d}{dt}\mathbb E^{n+1}(t)+2\mu_f\Vert\mathbb{D}(\boldsymbol{u}^{n+1})\Vert^2_f+\frac{L}{2}\left(\Vert\boldsymbol{u}^{n+1}\Vert^2_\Gamma + \Vert\boldsymbol{\xi}^{n+1}\Vert^2_\Gamma\right)\\
    &+\frac{1}{2L}\left(\Vert\mathcal{F}^{n+1}\Vert_\Gamma^2+\Vert\mathcal{S}^{n+1}\Vert_\Gamma^2\right)=\frac{L}{4}\Vert\tilde{\boldsymbol{u}}^{n}+\tilde{\boldsymbol{\xi}}^n\Vert^2_\Gamma+\frac{1}{4L}\Vert\tilde{\mathcal{F}}^{n}+\tilde{\mathcal S}^n\Vert^2_\Gamma\\
    &\le \frac{L}{2}\left(\Vert\tilde{\boldsymbol{u}}^n\Vert^2_\Gamma+\Vert\tilde{\boldsymbol{\xi}}^n\Vert^2_\Gamma\right)+\frac{1}{2L}\left(\Vert\tilde{\mathcal F}^n\Vert^2_\Gamma+\Vert\tilde{\mathcal S}^n\Vert^2_\Gamma\right).
    \end{aligned}
    \label{pointwise}
\end{equation}
Integrating \eqref{pointwise} over $[t_n, t_{n+1})]$ gives
\begin{equation}
    \mathbb{E}^{n+1} + \mathbb{H}^{n+1} + \mathbb{I}^{n+1} \le \mathbb{E}^n + \mathbb{I}^n.
    \label{cfinal}
\end{equation}
Here, for $n\ge 1$, the time-shift definition \eqref{tilde} implies
$$
\int^{t_{n+1}}_{t_n}\Vert\tilde{v}^n\Vert^2_\Gamma dt=\int^{t_{n}}_{t_{n-1}}\Vert v^n\Vert^2_\Gamma dt, v=\{\boldsymbol{u}, \boldsymbol{\xi}, \mathcal F, \mathcal S\}.
$$
Summing \eqref{cfinal} for $n=0, \cdots, M-1$ yields \eqref{csta}.
\end{proof}

\color{black}

\bibliographystyle{unsrtnat} 
\bibliography{sample}
\end{document}